\documentclass[12pt,letterpaper]{amsart}
\usepackage{euler, epic,eepic,latexsym, amssymb, amscd, amsfonts, xypic, floatflt}
\input xy
\xyoption{all}
 \newlength{\baseunit}               % the basic unit length
 \newcount{\numlines}                % depth of picture (in number of lines)
\newcommand{\point}{\vspace{3mm}\par \noindent \refstepcounter{subsection}{\bf \thesubsection.} }
\newcommand{\tpoint}[1]{\vspace{3mm}\par \noindent \refstepcounter{subsection}{\bf \thesubsection.} 
  {\em #1. ---} }
\newcommand{\epoint}[1]{\vspace{3mm}\par \noindent \refstepcounter{subsection}{\bf \thesubsection.} 
  {\em #1.} }
\newcommand{\bpoint}[1]{\vspace{3mm}\par \noindent \refstepcounter{subsection}{\bf \thesubsection.} 
  {\bf #1.} }
\newcommand{\bpf}{\noindent {\em Proof.  }}
\newcommand{\epf}{\qed \vspace{+10pt}}

\newcommand{\Z}{\mathbb{Z}}

\newcommand{\E}{\mathbb{E}}
\renewcommand{\L}{\mathbb{L}}
\newcommand{\Q}{\mathbb{Q}}

\newcommand{\C}{\mathbb{C}}
\newcommand{\F}{\mathbb{F}}

\renewcommand{\H}{\mathbb{H}}

\newcommand{\W}{\mathbb{W}}
\newcommand{\ka}{\kappa}
\newcommand{\proj}{\mathbb P}

\newcommand{\cC}{{\mathcal{C}}}
\newcommand{\cF}{{\mathcal{F}}}
\newcommand{\oh}{{\mathcal{O}}}

\newcommand{\cK}{{\mathcal{K}}}

\renewcommand{\cL}{{\mathcal{L}}}
\newcommand{\cm}{{\mathcal{M}}}
\newcommand{\cmbar}{\overline{\cm}}

\newcommand{\fr}{\mathfrak{r}}
\newcommand{\al}{\alpha}

\newcommand{\be}{\beta}
\newcommand{\ga}{\gamma}
\newcommand{\Ga}{\Gamma}
\newcommand{\de}{\delta}
\newcommand{\De}{\Delta}
\newcommand{\si}{\sigma}
\newcommand{\la}{\lambda}

\newcommand{\Hom}{\operatorname{Hom}}

\newcommand{\Pic}{\operatorname{Pic}}

\newcommand{\Sym}{\operatorname{Sym}}
\newcommand{\Aut}{\operatorname{Aut}}

\newcommand{\dis}{\bullet}
\newcommand{\virt}{\rm{virt}}
\newcommand{\fixed}{\rm{fixed}}
\renewcommand{\top}{\rm{top}}

\newcommand{\cited}{}

\newcommand{\lremind}[1]{{\bf[label:  #1]}}
\newcommand{\notation}[1]{}
\renewcommand{\lremind}[1]{{}}

\newcommand{\secretnote}[1]{}
\newcommand{\cut}[1]{}

\begin{document}
\pagestyle{plain}
\title{{\large {The moduli space of curves and Gromov-Witten theory}}
% was \Large
}
\author{Ravi Vakil}
\address{Dept. of Mathematics, Stanford University, Stanford CA~94305--2125}
\email{vakil@math.stanford.edu}
\thanks{Partially supported by NSF CAREER/PECASE Grant DMS--0228011, and 
an Alfred P. Sloan Research Fellowship.
\newline
\indent
2000 
Mathematics Subject Classification:  Primary 
14H10, 14H81, 14N35,
% 14N35 Gromov-Witten invariants, quantum cohomology
% 14H81 Curves, relationships with physics
% 14H10 Curves, families, moduli (algebraic)
Secondary  14N10, 53D45, 14H15.
% 14N10 enumerative problems (combinatorial problems)
% 53D45 Gromov-Witten invariants, quantum cohomology, Frobenius manifolds
% 14H15 Curves, families, moduli (analytic)
% not used:
% 32G15 Deformations of analytic structures, Moduli of Riemann surfaces, teichmuller theory
}
\date{Sunday, February 19, 2006.}
%\subjclass{Primary ??? Secondary ???. }
\begin{abstract}
  The goal of this article is to motivate and describe how
  Gromov-Witten theory can and has provided tools to understand
  the moduli space of curves.  For example, ideas and methods from
  Gromov-Witten theory have led to both conjectures and theorems
  showing that the tautological part of the cohomology ring has a
  remarkable and profound structure.  As an illustration, we describe
  a new approach to Faber's intersection number conjecture via
  branched covers of the projective line (work with I.P. Goulden and
  D.M. Jackson, based on work with T.  Graber).  En route we 
  review the work of a large number of mathematicians.
\end{abstract}
\maketitle
\tableofcontents

{\parskip=12pt % closing bracket is just before the bibliography 

\section{Introduction}
These notes are intended to explain how Gromov-Witten theory has been
useful in understanding the moduli space of complex curves.  
We will focus on the
moduli space of smooth curves and how much of the recent progress in
understanding it has come through ``enumerative'' invariants in
Gromov-Witten theory, something which we take for granted these days,
but should really be seen as surprising.  There is one sense in which
it should not be surprising --- in many circumstances, modern
arguments can be loosely interpreted as the fact that we can
understand curves in general by studying branched covers of the
complex projective line, as all curves can be so expressed.  We will
see this theme throughout the notes, from a Riemann-style parameter
count in \S \ref{riemann} to the tool of relative virtual localization
in Gromov Witten theory in \S \ref{s:rvl}.

These notes culminate in an approach to Faber's
intersection number conjecture using relative Gromov-Witten theory
(joint work with Goulden and Jackson \cite{gjv}).  One motivation for
this article is to convince the reader that our approach is natural and
straightforward.

We first introduce the {\em moduli space of curves}, both the moduli
space of smooth curves, and the Deligne-Mumford compactification,
which we will see is something forced upon us by nature, not
arbitrarily imposed by man.  We will then define certain geometrically
natural cohomology classes on the moduli space of smooth curves (the
tautological subring of the cohomology ring), and discuss Faber's
foundational conjectures on this subring.  We will then extend these
notions to the moduli space of stable curves, and discuss Faber-type
conjectures in this context.  A key example is Witten's conjecture,
which really preceded (and motivated) Faber's conjectures, and opened
the floodgates to the last decade's flurry of developments.  We will
then discuss other relations in the tautological ring (both known and
conjectural).  We will describe Theorem $\star$
(Theorem~\ref{thmstarhere}), a blunt tool for proving many statements,
and Y.-P. Lee's Invariance conjecture, which may give all relations in
the tautological ring.  In order to discuss the proof of Theorem
$\star$, we will be finally drawn into Gromov-Witten theory, and we
will quickly review the necessary background.  In particular, we will
need the notion of ``relative Gromov-Witten theory'', including Jun Li's
degeneration formula \cite{li1, li2} and the relative virtual
localization formula \cite{thmstar}.  Finally, we will use these ideas
to tackle Faber's intersection number conjecture.

Because the audience has a diverse background, this article is 
intended to be read at many different levels, with as much rigor as
the reader is able to bring to it. Unless the reader has a solid knowledge of
the foundations of algebraic geometry, which is most likely not the
case, he or she will have to be willing to take a few notions on
faith, and to ask a local expert a few questions.

We will cover a lot of ground, but hopefully this article will include
enough background that the reader can make explicit computations to
see that he or she can actively manipulate the ideas involved.  You
are strongly encouraged to try these ideas out via the exercises.
They are of varying difficulty, and the amount of rigor required for
their solution should depend on your
background.

Here are some suggestions for further reading.  For a gentle and quick
introduction to the moduli space of curves and its tautological ring,
see \cite{notices}.  For a pleasant and very detailed discussion of
moduli of curves, see Harris and Morrison's foundational book
\cite{hm}.  An on-line resource discussing curves and links to topology
(including a glossary of important terms) is available at \cite{aim}.
For more on curves, Gromov-Witten theory, and localization, see
\cite[Chapter 22--27]{mirsym}, which is intended for both physicists
and mathematicians.  Cox and Katz' wonderful book \cite{coxkatz} gives
an excellent mathematical approach to mirror symmetry.  There is as of
yet no ideal book introducing (Deligne-Mumford) stacks, but Fantechi's
\cite{fantechi} and Edidin's \cite{edidin} both give an excellent idea
of how to think about them and work with them, and the appendix to
Vistoli's paper \cite{vistoli} lays out the foundations directly,
elegantly, and quickly, although this is necessarily a more serious
read.

\noindent {\bf Acknowledgments.} I am grateful to the organizers of
the June 2005 conference in Cetraro, Italy on ``Enumerative invariants
in algebraic geometry and string theory'' (Kai Behrend, Barbara
Fantechi, and Marco Manetti), to Fondazione C.I.M.E. (Centro
Internazionale Matematico Estivo), and to the Hotel San Michele. I
learned this material from my co-authors Graber, Goulden, and Jackson,
and from the other experts in the field, including Carel Faber, Rahul
Pandharipande, Y.-P. Lee, \dots, whose names are mentioned throughout
this article.  I thank Carel Faber, Soren Galatius, Tom Graber,
Y.-P. Lee and Rahul Pandharipande for improving the manuscript.

\section{The moduli space of curves}

We begin with some conventions and terminology.
We will work over $\C$, although these questions
remain interesting over arbitrary fields. 
We will work algebraically, and hence only briefly mention 
other important approaches to the subjects, such as the construction
of the moduli space of curves as a quotient of Teichmuller space.

By {\em smooth curve}, we mean a compact (also known as proper or
complete), smooth (also known as nonsingular) complex curve, i.e.\ a
Riemann surface, see Figure~\ref{curve}.  Our curves will be connected
unless we especially describe them as ``possibly disconnected''.  In
general our {\em dimensions} will be algebraic or complex, which is why
we refer to a Riemann surface as a curve --- they have
algebraic/complex dimension $1$.  Algebraic geometers tend to draw
``half-dimensional'' cartoons of curves (see also Figure~\ref{curve}).

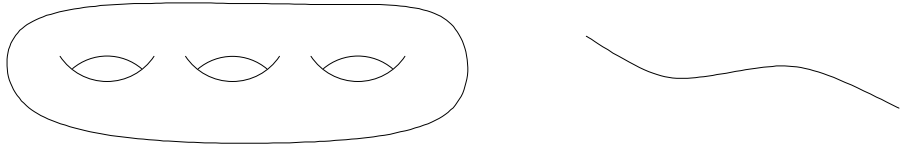
\begin{figure}
\begin{center}
\setlength{\unitlength}{0.00083333in}
\begingroup\makeatletter\ifx\SetFigFont\undefined%
\gdef\SetFigFont#1#2#3#4#5{%
  \reset@font\fontsize{#1}{#2pt}%
  \fontfamily{#3}\fontseries{#4}\fontshape{#5}%
  \selectfont}%
\fi\endgroup%
{\renewcommand{\dashlinestretch}{30}
\begin{picture}(5588,915)(0,-10)
\put(635.500,748.627){\arc{704.335}{0.5855}{2.5560}}
\put(636.000,221.000){\arc{667.083}{4.0119}{5.4129}}
\put(1418.500,748.627){\arc{704.335}{0.5855}{2.5560}}
\put(1418.500,219.744){\arc{669.649}{4.0132}{5.4116}}
\put(2200.500,748.500){\arc{704.195}{0.5852}{2.5563}}
\put(2201.000,221.250){\arc{666.701}{4.0114}{5.4134}}
\path(2025,880)(2060,880)(2096,880)
	(2131,880)(2165,880)(2200,880)
	(2234,880)(2269,880)(2303,879)
	(2337,878)(2370,876)(2403,875)
	(2436,872)(2467,869)(2498,866)
	(2528,862)(2556,857)(2583,852)
	(2609,845)(2634,839)(2657,831)
	(2679,823)(2699,814)(2718,804)
	(2736,793)(2753,782)(2769,769)
	(2783,755)(2797,741)(2810,725)
	(2822,708)(2833,689)(2844,670)
	(2853,650)(2861,628)(2869,606)
	(2875,583)(2879,559)(2883,535)
	(2885,511)(2886,486)(2886,462)
	(2884,438)(2881,414)(2876,391)
	(2870,368)(2863,346)(2855,325)
	(2845,305)(2834,286)(2822,268)
	(2809,250)(2795,234)(2780,220)
	(2764,206)(2747,193)(2728,180)
	(2708,168)(2687,157)(2664,146)
	(2639,135)(2614,125)(2587,116)
	(2558,107)(2529,98)(2498,90)
	(2466,83)(2434,76)(2401,69)
	(2367,63)(2333,58)(2299,53)
	(2265,48)(2230,44)(2196,41)
	(2162,38)(2128,35)(2094,32)
	(2060,30)(2026,27)(1992,25)
	(1963,24)(1932,22)(1902,21)
	(1871,20)(1839,18)(1807,17)
	(1775,16)(1741,15)(1707,14)
	(1673,14)(1638,13)(1603,13)
	(1567,12)(1531,12)(1495,12)
	(1458,12)(1422,13)(1386,13)
	(1350,13)(1314,14)(1279,15)
	(1245,16)(1210,17)(1177,18)
	(1144,19)(1111,21)(1080,23)
	(1049,24)(1018,26)(988,28)
	(958,30)(929,32)(896,35)
	(863,38)(830,41)(797,44)
	(763,48)(730,52)(696,57)
	(663,61)(629,67)(596,72)
	(563,79)(530,85)(497,92)
	(465,100)(434,108)(404,116)
	(375,125)(347,134)(320,144)
	(294,154)(269,164)(246,175)
	(224,186)(204,198)(184,210)
	(166,222)(149,234)(134,248)
	(118,263)(103,278)(89,294)
	(76,311)(64,329)(53,348)
	(44,367)(35,387)(28,408)
	(22,429)(17,450)(14,472)
	(12,494)(12,515)(12,537)
	(15,558)(18,578)(23,598)
	(30,618)(37,636)(46,654)
	(56,671)(67,687)(80,702)
	(93,717)(108,730)(122,742)
	(138,753)(155,764)(173,774)
	(193,784)(214,793)(237,802)
	(260,810)(286,818)(312,826)
	(340,833)(369,839)(400,845)
	(431,851)(463,856)(496,860)
	(529,865)(563,868)(597,872)
	(632,875)(666,877)(701,879)
	(736,881)(770,883)(805,884)
	(840,885)(875,886)(910,886)
	(941,887)(972,887)(1004,888)
	(1037,888)(1070,888)(1104,888)
	(1138,888)(1174,888)(1210,888)
	(1246,888)(1283,888)(1321,887)
	(1359,887)(1397,886)(1436,886)
	(1474,886)(1513,885)(1551,885)
	(1589,884)(1627,884)(1664,883)
	(1700,883)(1736,882)(1771,882)
	(1806,881)(1839,881)(1872,881)
	(1904,880)(1935,880)(1966,880)
	(1995,880)(2025,880)
\path(3626,681)(3628,679)(3634,676)
	(3643,670)(3657,661)(3677,648)
	(3701,633)(3729,615)(3760,596)
	(3793,576)(3827,556)(3861,537)
	(3895,519)(3927,501)(3958,486)
	(3987,472)(4015,460)(4042,450)
	(4067,441)(4092,434)(4116,428)
	(4140,424)(4164,420)(4189,418)
	(4211,418)(4234,418)(4258,418)
	(4282,420)(4307,422)(4333,425)
	(4359,428)(4386,432)(4414,436)
	(4442,441)(4471,446)(4500,451)
	(4529,456)(4559,461)(4588,466)
	(4618,471)(4647,476)(4676,480)
	(4704,484)(4732,487)(4759,490)
	(4786,492)(4812,494)(4838,494)
	(4863,494)(4889,493)(4914,492)
	(4939,489)(4964,485)(4990,480)
	(5017,474)(5044,466)(5073,458)
	(5103,447)(5135,436)(5169,423)
	(5204,408)(5241,392)(5280,375)
	(5319,357)(5358,339)(5397,320)
	(5434,302)(5468,286)(5498,271)
	(5524,258)(5544,248)(5559,240)
	(5568,235)(5574,232)(5576,231)
\end{picture}
}
\end{center}
\caption{A complex curve, and its real ``cartoon'' \lremind{curve} }
\label{curve}
\end{figure}

The reader likely needs no motivation to be interested in Riemann
surfaces.  A  natural question when you first hear of such
objects is: what are the Riemann surfaces?  How many of
them are there?  In other words, this question asks for a
classification of curves.

\epoint{Genus} A\label{genus}\lremind{genus} first invariant is the
{\em genus} of the smooth curve, which can be interpreted in three
ways: (i) the number of holes ({\em topological genus}; for example,
the genus of the curve in Figure~\ref{curve} is $3$), (ii) dimension
of space of space of differentials ($=h^0(C, \Omega_C)$, {\em
  geometric genus}), and (iii) the first cohomology group of the sheaf
of algebraic functions ($h^1(C, \oh_C)$, {\em arithmetic genus}).
These three notions are the same.  Notions (ii) and (iii) are related
by {\em Serre duality}\lremind{serreduality}
\begin{equation}\label{serreduality}
\boxed{H^0(C, \cF) \times
  H^1(C, \cK \otimes \cF^*) \rightarrow H^1(C, \cK) \cong \C}
\end{equation} where
$\cK$ is the canonical line bundle, which for smooth curves is the
sheaf of differentials $\Omega_C$.  Here $\cF$ can be any finite rank
vector bundle; $H^i$ refers to sheaf cohomology.  Serre duality
implies that $h^0(C, \cF) = h^1(C, \cK \otimes \cF^*)$, hence (taking
$\cF = \cK$).  $h^0(C, \Omega_C) = h^1(C, \oh_C)$.  (We will use these
important facts in the future!)

As we are working purely algebraically, we will not discuss why (i) is
the same as (ii) and (iii).

\bpoint{There is a $(3g-3)$-dimensional family of genus $g$ curves}

Remarkably, it was already known to  Riemann \cite[p.~134]{riemannref}
that there is a ``$3g-3$-dimensional family of genus $g$
curves''.\label{riemann}\lremind{riemann} You will notice that this
can't possibly be right if $g=0$, and you may know that this isn't
right if $g=1$, as you may have heard that elliptic curves are
parametrized by the $j$-line, which is one-dimensional.  So we will
take $g>1$, although there is a way to extend to $g=0$ and $g=1$  by making general enough
definitions.  (Thus there is a ``$(-3)$-dimensional moduli space'' of
genus $0$ curves, if you define moduli space appropriately --- in this
case as an Artin stack.  But that is another story.)

Let us now convince ourselves (informally) that there is a
$(3g-3)$-dimensional family of genus $g$ curves.  This will give me a
chance to introduce some useful facts that we will use later.  I will
use the same notation for vector bundles and their sheaves of
sections.  The sheaf of sections of a line bundle is called an {\em
  invertible sheaf}.

We will use five ingredients.

\noindent
{\bf (1)}  {\em Serre duality} \eqref{serreduality}.  
This is a hard fact.

\noindent
{\bf (2)}   {\em The Riemann-Roch formula.}
If $\cF$ is any coherent sheaf (for example, a finite rank vector bundle)
then $$
\boxed{h^0(C, \cF) - h^1(C, \cF) = \deg \cF - g + 1.}$$
This is an easy fact, although I will not explain why it is true.

\noindent
{\bf (3)} Line bundles of negative degree have no non-zero sections:
if $\cL$ is a line bundle of negative degree, then $\boxed{h^0(C, \cL) = 0}$.
Here is why: the degree of a line bundle $\cL$ can be defined as
follows.  Let $s$ be any non-zero meromorphic section of $\cL$.  Then
the degree of $\cL$ is the number of zeros of $s$ minus the number
of poles of $s$.  Thus if $\cL$ has an honest non-zero section (with no poles),
then the degree of $s$ is at least $0$. 

\noindent 
{\em Exercise.}  If $\cL$ is a degree $0$ line bundle with
a non-zero section $s$, show that $\cL$ is isomorphic to the trivial
bundle (the sheaf of functions) $\oh$.

\noindent
{\bf (4)}  Hence if $\cL$ is a line bundle with $\deg \cL > \deg \cK$,
then $h^1(C, \cL) = 0$ by Serre duality, from which
$\boxed{h^0(C, \cL) = \deg \cL - g + 1}$ by Riemann-Roch.

\noindent
{\bf (5)} {\em The Riemann-Hurwitz formula.}  Suppose $C \rightarrow
\proj^1$ is a degree $d$ cover of the complex projective line by a
genus $g$ curve $C$, with ramification $r_1$, \dots, $r_n$ at the
ramification points on $C$.  Then
$$
\boxed{\chi_{\top}(C) = d \chi_{\top} (\proj^1) -
  \sum (r_i-1),}$$
 where $\chi_{\top}$ is the topological Euler
characteristic, i.e. \lremind{rh} 
\begin{equation}
\boxed{2-2g = 2d - \sum(r_i-1).}
\label{rh}\end{equation}

We quickly review the language of divisors and line bundles on smooth
curves.  A {\em divisor} is a formal linear combination of points on
$C$, with integer co-efficients, finitely many non-zero.  A divisor is
{\em effective} if the co-efficients are non-negative.  The {\em
  degree} of a divisor is the sum of its co-efficients. Given a
divisor $D = \sum n_i p_i$ (where the $p_i$ form a finite set), we
obtain a line bundle $\oh(D)$ by ``twisting the trivial bundle $n_i$
times at the point $p_i$''.  This is best understood in terms of the
sheaf of sections.  Sections of the sheaf $\oh(D)$ (over some open
set) correspond to meromorphic functions that are holomorphic away
from the $p_i$; and if $n_i>0$, have a pole of order at most $n_i$ at
$p_i$; and if $n_i<0$, have a zero of order at least $-n_i$ at $p_i$.
Each divisor yields a line bundle along with a meromorphic section
(obtained by taking the function $1$ in the previous sentence's
description).  Conversely, each line bundle with a non-zero
meromorphic section yields a divisor, by taking the ``divisor of zeros
and poles'': if $s$ is a non-zero meromorphic section, we take the
divisor which is the sum of the zeros of $s$ (with multiplicity) minus
the sum of the poles of $s$ (with multiplicity).  These two
constructions are inverse to each other.  In short, line bundles with
the additional data of a non-zero {\em meromorphic} section correspond
to divisors.  This identification is actually quite subtle the first
few times you see it, and it is worth thinking through it carefully if
you have not done so before.  Similarly, line bundles with the
additional data of a non-zero {\em holomorphic} section correspond to
{\em effective} divisors.

We now begin our dimension count.  We do it in three steps.

\noindent 
{\bf Step 1.}
Fix a  curve $C$, and a degree $d$.  Let $\Pic^d C$ be the set of
degree $d$ line bundles on $C$. Pick a point $p \in C$.  Then there is
an bijection $\Pic^0 C \rightarrow \Pic^d C$ given by $\cF \rightarrow
\cF(d p)$.  (By $\cF(d p)$, we mean the ``twist of $\cF$ at $p$, $d$
times'', which is the same construction sketched two paragraphs previously.
 In terms of sheaves, if $d>0$, this means the sheaf of
meromorphic sections of $\cF$, that are required to be holomorphic
away from $p$, but may have a pole of order at most $d$ at $p$.  If
$d<0$, this means the sheaf of holomorphic sections of $\cF$ that are
required to have a zero of order at least $-d$ at $p$.)  If we believe
$\Pic^d C$ has some nice structure, which is indeed the case, then we
would expect that this would be an isomorphism.  In fact, $\Pic^d$ can
be given the structure of a complex manifold or complex variety, and
this gives an isomorphism of manifolds or varieties.

\noindent {\bf Step 2: ``$\dim \Pic^d C = g$.''}  There are quotes
around this equation because so far, $\Pic^d C$ is simply a set, so
this will just be a plausibility argument.  By Step 1, it suffices to
consider any $d > \deg \cK$.  Say $\dim \Pic^d C = h$.  We ask: how
many degree $d$ {\em effective divisors} are there (i.e.\ what is the
dimension of this family)?  The answer is clearly $d$, and $C^d$
surjects onto this set (and is usually $d!$-to-$1$).

But we can count effective divisors in a different way.
There is an $h$-dimensional family of line bundles by hypothesis,
and each one of these has a $(d-g+1)$-dimensional family of
non-zero sections, each of which gives a divisor of zeros.  But two
sections yield the same divisor if one is a multiple of the other.
Hence we get: $h + (d-g+1) - 1 = h+d-g$.

Thus $d=h+d-g$, from which $h=g$ as desired.

Note that we get a bit more: if we believe that $\Pic^d$ has an
algebraic structure, we have a fibration $(C^d) / S_d \rightarrow
\Pic^d$, where the fibers are isomorphic to $\proj^{d-g}$.  In
particular, $\Pic^d$ is reduced (I won't define this!), and
irreducible.  (In fact, as many of you know, it is isomorphic to the
dimension $g$ abelian variety $\Pic^0 C$.)

\noindent {\bf Step 3.}
Say $\cm_g$ has dimension $p$.  By fact {\bf (4)} above, if 
$d \gg 0$, and $D$ is a
divisor of degree $d$, then $h^0(C, \oh(D)) = d-g+1$.  
If we take two general sections $s$, $t$ of the line bundle $\oh(D)$,
we get a map to $\proj^1$ (given by $p \rightarrow [s(p);t(p)]$
--- note that this is well-defined), and this map is degree
$d$ (the preimage of $[0;1]$ is precisely $\operatorname{div} s$,
which has $d$ points counted with multiplicity).
Conversely, any degree $d$ cover  $f: C \rightarrow \proj^1$ arises
from two linearly independent sections of a degree $d$
line bundle.  (To get the divisor associated to one of them,
consider $f^{-1}([0;1])$, where points are counted with multiplicities;
to get the divisor associated to the other, consider $f^{-1}([1;0])$.)
Note that $(s,t)$ gives the same map to $\proj^1$ as $(s', t')$
if and only $(s,t)$ is a scalar multiple of $(s', t')$.
Hence the number of maps to $\proj^1$ arising from a fixed
curve $C$ and a fixed line bundle $\cL$ correspond
to the choices of two sections ($2(d-g+1)$ by fact {\bf (4)}),
minus $1$ to forget the scalar multiple, for a total
of $2d-2g+1$.  If we let the the line bundle vary, the
number of maps from a fixed curve is $2d-2g+1+\dim \Pic^d(C) =
2d-g+1$.  If we let the curve also vary, we see that the number
of degree $d$ genus $g$ covers of $\proj^1$ is 
$\boxed{p+2d-g+1}$.

But we can also count this number using the Riemann-Hurwitz formula
\eqref{rh}.  By that formula, there will be a total of $2g+2d-2$ branch points
(including multiplicity). 
Given the branch points (again, with multiplicity), there is a finite amount
of possible monodromy data around the branch points.
The Riemann Existence Theorem tells us that given any such monodromy
data, we can uniquely reconstruct the cover, so we have
$$
p+2d-g+1 = 2g+2d-2,$$
from which
$\boxed{p=3g-3}$.

Thus there is a $3g-3$-dimensional family of genus $g$ curves!  (By
showing that the space of branched covers is reduced and irreducible,
we could again ``show'' that the moduli space is reduced and irreducible.)

\bpoint{The moduli space of smooth curves}

It is time to actually define the moduli space of genus $g$ smooth curves,
denoted $\cm_g$, or at least to come close to it.  By ``moduli space
of curves'' we mean a ``parameter space for curves''.  As a first
approximation, we mean the set of curves, but we want to endow this
set with further structure (ideally that of a manifold, or even of a
smooth complex variety).  This structure should be given by nature,
not arbitrarily defined.

Certainly if there were such a space $\cm_g$, we would expect a
universal curve over it $\cC_g \rightarrow \cm_g$, so that the fiber
above the point $[C]$ representing a curve $C$ would be that same $C$.
Moreover, whenever we had a family of curves parametrized by some base
$B$, say $\cC_B \rightarrow B$ (where the fiber above any point $b \in
B$ is some smooth genus $g$ curve $C_b$), there should be a map $f: B
\rightarrow \cm_g$ (at the level of sets sending $b \in B$ to $[C_b]
\in \cm_g$), and then $f^* \cC_g$ should be isomorphic to $\cC_B$.

We can turn this into a precise definition.  The families we should
consider should be ``nice'' (``fibrations'' in the sense of
differential geometry).  It turns out that the corresponding algebraic
notion of ``nice'' is {\em flat}, which I will not define here.  We
can {\em define} $\cm_g$ to be the scheme such that the maps from any
scheme $B$ to it are in natural bijection with nice (flat) families of
genus $g$ curves over $B$. 
(Henceforth all families will be
assumed to be ``nice''=flat.)
 Some thought will convince you that only
one space (up to isomorphism) exists with this property.  This
``abstract nonsense'' is called {\em Yoneda's Lemma}.  The argument is
general, and  applies to nice families of any sort of
thing.  Categorical translation: we are saying that this contravariant
functor of families is {\em represented} by the functor $\Hom(\cdot,
\cm_g)$.  Translation: if such a space exists, then it is unique, up to
unique isomorphism.

If there is such a moduli space $\cm_g$, we gain some additional
information: cohomology classes on $\cm_g$ are ``characteristic
classes'' for families of genus $g$ curves.  More precisely, given any
family of genus $g$ curves $\cC_B \rightarrow B$, and any cohomology class $\al
\in H^*(\cm_g)$, we have a cohomology class on $B$:
if $f: B \rightarrow \cm_g$ is the moduli map, take $f^*
\al$.  These characteristic classes behave well with respect to pullback:
if $\cC_{B'} \rightarrow B'$ is a family obtained by pullback from 
$\cC_B \rightarrow B$, then the cohomology class on $B'$ induced by $\al$
is the pullback of the cohomology class on $B$ induced by $\al$.
The converse turns out to be true:  any such ``universal cohomology class'',
defined for all families and well-behaved under pullback, arises
from a cohomology class on $\cm_g$.  (The argument is actually quite
tautological, and the reader is invited to think it through.)
More generally, statements about the geometry of $\cm_g$ correspond
to ``universal statements about all families''.

Here is an example of a consequence.  A curve is {\em hyperelliptic}
if it admits a $2$-to-$1$ cover of $\proj^1$.  In the space of smooth
genus $3$ curves $\cm_3$, there is a Cartier divisor of hyperelliptic
curves, which means that the locus of hyperelliptic curves is locally
cut out by a single equation.  Hence in {\em any} family of genus $3$
curves over an arbitrarily horrible base, the hyperelliptic locus are
cut out by a single equation.  (For scheme-theoretic experts: for any
family $\cC_B \rightarrow B$ of genus $3$ curves, there is then a
closed subscheme of $B$ corresponding to the hyperelliptic locus.
What is an intrinsic scheme-theoretic definition of this locus?)

Hence all we have to do is show that there is such a scheme $\cm_g$.
Sadly, there is no such scheme!  We could just throw up our hands and
end these notes here.  There are two patches to this problem. One
solution is to relax the definition of moduli space (to get the notion
of {\em coarse moduli space}), which doesn't quite parametrize all
families of curves.  A second option is to extend the notion of {\em
  space}.  The first choice is the more traditional one, but it is
becoming increasingly clear that the second one is the better one.

This leads us to the notion of a {\em stack}, or in this case, the
especially nice stack known as a {\em Deligne-Mumford stack}.  This is
an extension of the idea of an idea of a scheme.  Defining a
Deligne-Mumford stack correctly takes some time, and is rather tiring
and uninspiring, but dealing with Deligne-Mumford stacks on a
day-to-day basis is not so bad --- you just pretend it is a scheme.
One might compare it to driving a car without knowing how the engine
works, but really it is more like driving a car while having only the
vaguest idea of what a car is.

Thus I will content myself with giving you a few cautions about where
your informal notion of Deligne-Mumford stack should differ with your
notion of scheme.  (I feel less guilty about this knowing that many
analytic readers will be similarly uncomfortable with the notion of a
scheme.)  The main issue is that when considering
cohomology rings (or the algebraic analog, Chow rings), we will
take $\Q$-co-efficients in order to avoid subtle technical issues.
The foundations of intersection theory for Deligne-Mumford stacks
were laid by Vistoli in \cite{vistoli}
(However, thanks to work of Andrew Kresch \cite{kresch}, it
is possible to take integral co-efficients using the Chow ring.
Then we have to accept the fact that cohomology groups can be non-zero
even in degree higher than the dimension of the space.  This is actually
something that for various reasons we {\em want} to be true, but such
a discussion is not appropriate in these notes.) 

A smooth (or nonsingular) Deligne-Mumford stack (over $\C$) is
essentially the same thing as a complex orbifold.  The main caution
about saying that they are the same thing is that there are actually
three different definitions of orbifold in use, and many users are
convinced that their version is the only version in use, causing
confusion for readers such as myself.

Hence for the rest of these notes, we will take for granted
that there is a moduli space of smooth curves $\cm_g$ (and we 
will make similar assumptions about other moduli spaces).

Here are some {\em facts} about the moduli space of curves.  The space
$\cm_g$ has (complex) dimension $3g-3$.  It is smooth (as a stack), so
it is an orbifold (given the appropriate definition), and we will
imagine that it is a manifold.  We have informally seen that it is
irreducible.

We make a brief brief excursion outside of algebraic geometry to show
that this space has some interesting structure.  In the analytic
setting, $\cm_g$ can be expressed as the quotient of {\em Teichmuller
  space} (a subset of $\C^{3g-3}$ homeomorphic to a ball) by a
discrete group, known as the {\em mapping class group}.  Hence the
cohomology of the quotient $\cm_g$ is the group cohomology of the
mapping class group.  (Here it is essential that we take the quotient
as an orbifold/stack.)  Here is a fact suggesting that the topology of
this space has some elegant structure: \lremind{hzeq}
\begin{equation}\label{hzeq}
\chi(\cm_{g}) =
B_{2g} / 2g (2g-2)
\end{equation} (due to Harer and Zagier \cite{hz}), where
$B_{2g}$ denotes the $2g$th Bernoulli number.

Other exciting recent work showing the attractive structure of the
cohomology ring is Madsen and Weiss' proof of Madsen's generalization
of Mumford's conjecture \cite{mw}.  We briefly give the statement.
There is a natural isomorphism between $H^*(\cm_g;\Q)$ and $H^*(\cm_{g+1};\Q)$
for $*<(g-1)/2$ (due to Harer and Ivanov).
Hence we can define the ring we could informally denote by
$H^*(\cm_{\infty}; \Q)$.  Mumford conjectured that this is a free
polynomial ring generated by certain cohomology classes
($\ka$-classes, to be defined in \S \ref{taut1}).  Madsen and Weiss proved this,
and a good deal more.  (See \cite{tillmann} for an overview of the
topological approach to the Mumford conjecture, and \cite{mt} for a
more technical discussion.)

\bpoint{Pointed nodal curves, and the moduli space of stable
pointed curves}

As our moduli space $\cm_g$ is a smooth orbifold of dimension $3g-3$,
it is wonderful in all ways but one:  it is not compact.  It would 
be useful to have a good compactification, one that is still smooth,
and also has good geometric meaning.
This leads us to extend our notion of smooth curves slightly.

A {\em node} of a curve is a singularity analytically isomorphic to
$xy=0$ in $\C^2$.  A {\em nodal curve} is a curve (compact, connected)
smooth away from finite number of points (possibly zero), which are
nodes.  An example is sketched in Figure~\ref{nodalcurve}, in both
``real'' and ``cartoon'' form.  One caution with the ``real''
picture: the two branches at the node are not tangent; this optical
illusion arises from the need of our limited brains to represent the
picture in three-dimensional space.  A {\em pointed nodal curve}
is a nodal curve with the additional data of $n$ distinct smooth points 
labeled $1$ through $n$ (or $n$ distinct labels of your choice,
such as $p_1$ through $p_n$).

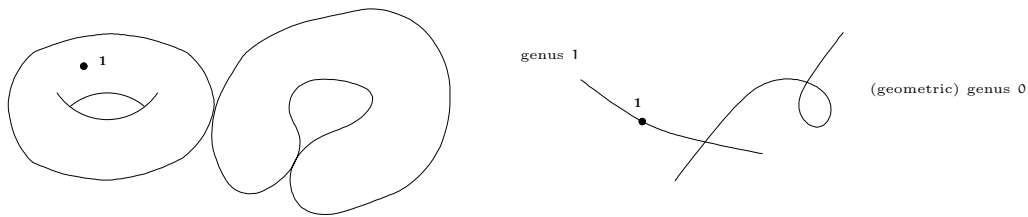
\begin{figure}
\begin{center}
\setlength{\unitlength}{0.00083333in}
\begingroup\makeatletter\ifx\SetFigFont\undefined%
\gdef\SetFigFont#1#2#3#4#5{%
  \reset@font\fontsize{#1}{#2pt}%
  \fontfamily{#3}\fontseries{#4}\fontshape{#5}%
  \selectfont}%
\fi\endgroup%
{\renewcommand{\dashlinestretch}{30}
\begin{picture}(5387,1323)(0,-10)
\put(630.500,980.255){\arc{757.731}{0.5915}{2.5501}}
\put(630.500,410.223){\arc{718.841}{4.0162}{5.4086}}
\path(3585,853)(3587,851)(3593,847)
	(3603,840)(3617,829)(3636,814)
	(3660,797)(3687,777)(3716,756)
	(3746,734)(3776,713)(3806,692)
	(3835,673)(3862,655)(3889,638)
	(3914,623)(3938,609)(3961,596)
	(3984,584)(4007,573)(4030,563)
	(4053,553)(4073,545)(4094,537)
	(4115,530)(4138,522)(4161,515)
	(4186,508)(4213,500)(4242,493)
	(4272,485)(4305,477)(4340,469)
	(4376,461)(4415,453)(4454,445)
	(4494,436)(4533,428)(4571,421)
	(4606,414)(4637,408)(4663,402)
	(4684,398)(4699,395)(4709,394)
	(4715,392)(4717,392)
\path(4172,224)(4174,226)(4177,231)
	(4183,239)(4193,252)(4206,269)
	(4224,292)(4244,319)(4268,349)
	(4294,382)(4322,418)(4350,454)
	(4379,490)(4408,525)(4437,559)
	(4464,592)(4490,622)(4515,650)
	(4538,676)(4561,699)(4582,720)
	(4602,740)(4621,757)(4640,772)
	(4658,786)(4675,798)(4693,809)
	(4710,818)(4730,828)(4750,836)
	(4770,843)(4791,849)(4811,853)
	(4832,856)(4852,858)(4873,859)
	(4893,858)(4914,856)(4934,854)
	(4953,850)(4972,845)(4990,839)
	(5007,833)(5023,826)(5039,818)
	(5053,810)(5066,801)(5078,793)
	(5089,783)(5099,774)(5107,765)
	(5115,755)(5122,745)(5129,734)
	(5134,723)(5138,712)(5141,700)
	(5143,689)(5144,677)(5145,666)
	(5144,654)(5142,643)(5140,632)
	(5136,622)(5132,612)(5127,603)
	(5122,594)(5116,587)(5109,580)
	(5103,574)(5096,569)(5088,565)
	(5081,562)(5073,560)(5065,559)
	(5057,558)(5048,559)(5040,560)
	(5031,563)(5022,566)(5013,571)
	(5004,576)(4996,583)(4987,590)
	(4980,598)(4972,607)(4966,617)
	(4960,628)(4955,639)(4951,650)
	(4948,663)(4945,675)(4944,688)
	(4944,701)(4945,714)(4947,728)
	(4950,742)(4955,757)(4961,773)
	(4968,789)(4978,808)(4989,827)
	(5002,849)(5017,872)(5034,897)
	(5053,925)(5073,953)(5095,983)
	(5117,1013)(5139,1042)(5160,1070)
	(5179,1094)(5194,1114)(5206,1129)
	(5214,1139)(5218,1144)(5220,1147)
\put(484,937){\blacken\ellipse{42}{42}}
\put(484,937){\ellipse{42}{42}}
\put(3967,591){\blacken\ellipse{42}{42}}
\put(3967,591){\ellipse{42}{42}}
\path(1783,1188)(1825,1202)(1868,1215)
	(1909,1226)(1948,1235)(1984,1243)
	(2018,1250)(2048,1256)(2076,1262)
	(2101,1267)(2123,1271)(2143,1276)
	(2161,1280)(2178,1283)(2194,1286)
	(2209,1289)(2224,1292)(2240,1294)
	(2256,1295)(2273,1296)(2292,1296)
	(2313,1294)(2335,1292)(2360,1287)
	(2387,1281)(2416,1273)(2447,1262)
	(2480,1248)(2514,1232)(2547,1211)
	(2579,1188)(2610,1159)(2638,1129)
	(2662,1097)(2682,1066)(2699,1037)
	(2714,1009)(2726,983)(2735,959)
	(2743,937)(2749,916)(2754,898)
	(2758,880)(2762,863)(2764,846)
	(2766,829)(2768,811)(2769,793)
	(2769,773)(2769,751)(2769,726)
	(2767,699)(2765,669)(2761,635)
	(2755,599)(2747,560)(2736,519)
	(2722,476)(2705,434)(2685,396)
	(2663,360)(2640,327)(2616,297)
	(2592,270)(2569,247)(2547,226)
	(2525,209)(2505,193)(2486,180)
	(2468,169)(2451,160)(2434,151)
	(2418,144)(2401,137)(2385,131)
	(2368,124)(2351,118)(2333,111)
	(2313,104)(2292,96)(2270,87)
	(2246,78)(2219,69)(2191,58)
	(2162,48)(2131,38)(2099,29)
	(2066,21)(2035,15)(2003,12)
	(1974,12)(1948,15)(1925,19)
	(1904,25)(1887,31)(1872,38)
	(1859,46)(1848,53)(1838,61)
	(1830,68)(1823,76)(1817,84)
	(1811,92)(1806,100)(1801,109)
	(1797,118)(1792,128)(1788,139)
	(1783,152)(1779,165)(1776,181)
	(1773,198)(1771,216)(1770,237)
	(1772,259)(1776,283)(1783,308)
	(1796,338)(1813,367)(1833,394)
	(1855,419)(1879,442)(1903,461)
	(1928,478)(1953,492)(1978,505)
	(2003,516)(2027,526)(2052,535)
	(2076,544)(2101,554)(2125,564)
	(2149,576)(2173,589)(2196,604)
	(2218,621)(2238,640)(2256,660)
	(2271,683)(2281,705)(2286,727)
	(2285,747)(2278,764)(2267,779)
	(2253,792)(2236,803)(2218,812)
	(2200,820)(2180,827)(2161,832)
	(2141,837)(2121,842)(2101,846)
	(2080,849)(2060,852)(2038,855)
	(2016,856)(1994,857)(1970,857)
	(1946,856)(1921,852)(1896,847)
	(1871,838)(1847,826)(1825,811)
	(1804,790)(1787,766)(1775,741)
	(1768,717)(1764,693)(1764,670)
	(1767,650)(1772,630)(1780,612)
	(1788,594)(1797,577)(1806,560)
	(1815,541)(1822,522)(1828,502)
	(1832,479)(1833,455)(1830,428)
	(1824,399)(1814,369)(1800,338)
	(1783,308)(1765,283)(1745,261)
	(1726,242)(1707,226)(1689,212)
	(1673,200)(1658,190)(1644,182)
	(1632,175)(1620,169)(1609,163)
	(1598,158)(1587,154)(1576,150)
	(1564,147)(1551,144)(1538,143)
	(1522,142)(1505,142)(1487,145)
	(1467,149)(1447,157)(1426,168)
	(1406,183)(1388,202)(1372,225)
	(1358,248)(1346,270)(1336,292)
	(1327,313)(1319,331)(1313,348)
	(1307,364)(1302,379)(1298,392)
	(1294,406)(1291,420)(1288,435)
	(1285,451)(1283,469)(1282,490)
	(1282,514)(1283,541)(1286,572)
	(1290,607)(1298,645)(1308,685)
	(1322,727)(1338,765)(1355,801)
	(1373,835)(1390,865)(1406,891)
	(1420,914)(1433,934)(1443,950)
	(1452,964)(1460,976)(1467,986)
	(1473,995)(1479,1003)(1486,1011)
	(1493,1019)(1502,1028)(1513,1037)
	(1526,1048)(1543,1061)(1564,1075)
	(1589,1092)(1618,1110)(1653,1129)
	(1692,1149)(1736,1169)(1783,1188)
\path(400,1105)(438,1114)(475,1121)
	(510,1126)(540,1130)(566,1134)
	(587,1136)(605,1138)(620,1139)
	(632,1139)(642,1140)(652,1140)
	(661,1140)(671,1139)(683,1139)
	(698,1138)(716,1136)(737,1134)
	(763,1130)(793,1126)(828,1121)
	(865,1114)(903,1105)(944,1094)
	(980,1082)(1012,1071)(1039,1061)
	(1060,1053)(1076,1046)(1089,1041)
	(1098,1037)(1106,1034)(1112,1031)
	(1118,1028)(1125,1024)(1132,1019)
	(1141,1011)(1153,1001)(1167,987)
	(1183,970)(1201,948)(1220,923)
	(1238,895)(1253,865)(1265,836)
	(1274,810)(1281,787)(1286,767)
	(1290,751)(1293,738)(1295,728)
	(1296,718)(1297,710)(1298,702)
	(1297,692)(1296,681)(1294,668)
	(1291,652)(1286,632)(1278,608)
	(1268,580)(1255,550)(1238,518)
	(1218,487)(1196,460)(1176,436)
	(1157,416)(1141,400)(1128,388)
	(1116,378)(1107,371)(1099,365)
	(1091,361)(1083,356)(1074,351)
	(1063,346)(1048,339)(1030,330)
	(1006,320)(977,307)(943,294)
	(904,280)(861,266)(821,256)
	(782,247)(747,241)(716,236)
	(689,232)(667,230)(649,228)
	(635,227)(622,227)(612,227)
	(602,228)(593,228)(583,229)
	(571,230)(557,232)(539,234)
	(518,237)(493,240)(463,245)
	(430,250)(395,257)(358,266)
	(320,277)(285,289)(256,300)
	(232,310)(213,317)(199,323)
	(188,328)(180,331)(174,334)
	(170,336)(165,339)(161,342)
	(155,347)(148,355)(138,365)
	(127,379)(113,397)(97,419)
	(81,446)(65,476)(53,505)
	(43,534)(35,561)(29,585)
	(24,607)(21,625)(18,641)
	(16,654)(14,666)(13,676)
	(13,686)(12,695)(12,705)
	(13,717)(14,730)(17,746)
	(20,764)(25,786)(31,810)
	(40,837)(51,866)(65,895)
	(83,925)(102,952)(120,974)
	(136,992)(150,1006)(162,1016)
	(171,1024)(178,1029)(185,1032)
	(191,1035)(197,1037)(205,1040)
	(214,1043)(227,1048)(243,1054)
	(264,1061)(291,1071)(323,1082)
	(359,1094)(400,1105)
\put(3208,979){\makebox(0,0)[lb]{\smash{{{\SetFigFont{5}{6.0}{\rmdefault}{\mddefault}{\updefault}genus $1$}}}}}
\put(5387,769){\makebox(0,0)[lb]{\smash{{{\SetFigFont{5}{6.0}{\rmdefault}{\mddefault}{\updefault}(geometric) genus $0$}}}}}
\put(3920,686){\makebox(0,0)[lb]{\smash{{{\SetFigFont{5}{6.0}{\rmdefault}{\mddefault}{\updefault}{\bf 1}}}}}}
\put(581,948){\makebox(0,0)[lb]{\smash{{{\SetFigFont{5}{6.0}{\rmdefault}{\mddefault}{\updefault}{\bf 1}}}}}}
\end{picture}
}
\end{center}
\caption{A pointed nodal curve, and its real ``cartoon'' \lremind{nodalcurve}}
\label{nodalcurve}
\end{figure}

The {\em geometric genus}  of an {\em irreducible} 
curve is its genus once all of the nodes are ``unglued''.  For example,
the components of the curve in Figure~\ref{nodalcurve}
have genus $1$ and $0$.

We define the {\em (arithmetic) genus} of a pointed nodal curve
informally as the genus of a ``smoothing'' of the curve,
which is indicated in Figure~\ref{smoothing}.
More formally, we define it as  $h^1(C, \oh_C)$.  
This notion behaves well with respect to deformations.
(More formally, it is locally constant in flat families.)

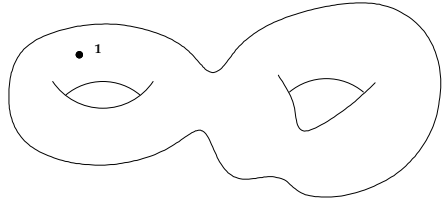
\begin{figure}
\begin{center}
\setlength{\unitlength}{0.00083333in}
\begingroup\makeatletter\ifx\SetFigFont\undefined%
\gdef\SetFigFont#1#2#3#4#5{%
  \reset@font\fontsize{#1}{#2pt}%
  \fontfamily{#3}\fontseries{#4}\fontshape{#5}%
  \selectfont}%
\fi\endgroup%
{\renewcommand{\dashlinestretch}{30}
\begin{picture}(2718,1251)(0,-10)
\put(596.500,942.255){\arc{757.731}{0.5915}{2.5501}}
\put(596.500,372.223){\arc{718.841}{4.0162}{5.4086}}
\put(1992.500,390.223){\arc{718.841}{4.0162}{5.4086}}
\put(450,899){\blacken\ellipse{42}{42}}
\put(450,899){\ellipse{42}{42}}
\path(973,1021)(960,1026)(946,1032)
	(931,1037)(916,1042)(899,1047)
	(881,1053)(863,1058)(845,1063)
	(826,1067)(807,1071)(788,1075)
	(770,1079)(751,1082)(734,1085)
	(716,1087)(699,1089)(682,1090)
	(665,1091)(647,1092)(629,1093)
	(610,1093)(591,1092)(571,1092)
	(552,1091)(532,1089)(513,1087)
	(494,1085)(476,1083)(458,1080)
	(441,1078)(424,1075)(408,1071)
	(392,1068)(376,1064)(360,1060)
	(343,1056)(326,1051)(309,1046)
	(292,1040)(275,1034)(259,1028)
	(243,1021)(227,1015)(213,1008)
	(199,1001)(186,995)(174,988)
	(163,981)(151,973)(139,964)
	(128,955)(116,945)(105,934)
	(95,922)(85,910)(76,897)
	(67,883)(59,870)(52,856)
	(46,842)(41,828)(36,813)
	(32,800)(28,786)(25,771)
	(22,755)(19,739)(17,722)
	(15,704)(13,686)(12,667)
	(12,649)(12,631)(12,614)
	(13,597)(14,581)(16,566)
	(18,551)(21,534)(25,517)
	(30,501)(35,484)(42,467)
	(49,450)(58,434)(67,418)
	(77,403)(88,388)(99,375)
	(111,362)(123,351)(136,340)
	(147,331)(160,322)(174,313)
	(188,304)(204,296)(220,287)
	(238,279)(256,272)(274,264)
	(292,258)(311,251)(329,246)
	(347,241)(365,236)(382,232)
	(399,228)(416,225)(434,222)
	(452,220)(470,218)(490,216)
	(509,214)(529,213)(550,212)
	(570,212)(590,212)(610,212)
	(629,213)(648,215)(666,216)
	(684,218)(702,220)(719,223)
	(737,226)(755,229)(773,233)
	(792,238)(811,242)(831,248)
	(850,253)(869,259)(888,265)
	(906,272)(923,278)(939,285)
	(954,292)(969,298)(983,305)
	(998,313)(1014,322)(1029,331)
	(1043,340)(1058,350)(1073,359)
	(1087,369)(1100,378)(1113,387)
	(1124,395)(1135,402)(1145,408)
	(1154,414)(1163,418)(1170,422)
	(1176,425)(1183,427)(1189,428)
	(1196,429)(1202,429)(1208,429)
	(1214,427)(1219,425)(1225,422)
	(1230,418)(1235,414)(1239,409)
	(1244,403)(1248,397)(1252,390)
	(1258,380)(1264,368)(1270,355)
	(1276,341)(1283,325)(1290,309)
	(1298,292)(1305,276)(1312,260)
	(1320,246)(1327,232)(1335,219)
	(1343,207)(1351,196)(1360,185)
	(1370,174)(1380,165)(1391,156)
	(1403,148)(1414,141)(1426,136)
	(1437,131)(1448,128)(1459,126)
	(1471,124)(1483,124)(1496,124)
	(1509,124)(1523,126)(1537,127)
	(1551,129)(1565,131)(1579,133)
	(1591,135)(1604,137)(1615,138)
	(1629,139)(1642,139)(1656,139)
	(1670,137)(1684,135)(1697,132)
	(1710,129)(1721,124)(1731,119)
	(1740,114)(1749,107)(1757,100)
	(1766,91)(1776,82)(1786,73)
	(1798,64)(1810,55)(1822,46)
	(1836,39)(1851,32)(1862,28)
	(1875,24)(1890,21)(1905,18)
	(1923,16)(1941,14)(1961,13)
	(1981,12)(2002,12)(2023,13)
	(2044,15)(2065,17)(2085,19)
	(2106,23)(2122,26)(2139,30)
	(2156,34)(2174,39)(2193,44)
	(2212,50)(2231,57)(2251,64)
	(2270,72)(2290,80)(2309,89)
	(2328,98)(2347,107)(2365,117)
	(2382,127)(2398,137)(2414,147)
	(2429,158)(2444,169)(2458,180)
	(2473,193)(2487,206)(2502,220)
	(2516,234)(2530,250)(2544,266)
	(2558,282)(2571,299)(2583,316)
	(2595,334)(2606,351)(2617,368)
	(2626,385)(2635,402)(2643,419)
	(2650,436)(2657,453)(2664,471)
	(2670,489)(2675,508)(2681,528)
	(2686,548)(2690,569)(2694,590)
	(2697,611)(2700,632)(2703,653)
	(2704,674)(2705,694)(2706,713)
	(2706,732)(2705,750)(2704,768)
	(2703,785)(2701,804)(2698,822)
	(2694,841)(2690,860)(2685,879)
	(2679,899)(2672,918)(2665,937)
	(2657,955)(2648,973)(2639,990)
	(2630,1006)(2620,1022)(2610,1036)
	(2600,1050)(2589,1062)(2578,1075)
	(2566,1087)(2553,1099)(2540,1111)
	(2525,1123)(2510,1134)(2494,1145)
	(2477,1155)(2460,1165)(2442,1174)
	(2425,1183)(2407,1190)(2390,1197)
	(2372,1202)(2355,1207)(2338,1211)
	(2323,1215)(2307,1217)(2290,1220)
	(2273,1222)(2256,1223)(2237,1224)
	(2218,1224)(2199,1224)(2179,1224)
	(2160,1223)(2140,1222)(2120,1220)
	(2101,1218)(2082,1216)(2063,1214)
	(2045,1211)(2027,1208)(2010,1205)
	(1992,1202)(1974,1198)(1956,1194)
	(1937,1190)(1918,1186)(1898,1181)
	(1878,1176)(1858,1170)(1838,1165)
	(1818,1159)(1799,1153)(1780,1147)
	(1762,1140)(1744,1134)(1728,1128)
	(1712,1122)(1696,1116)(1682,1109)
	(1666,1102)(1650,1095)(1634,1087)
	(1618,1079)(1602,1070)(1587,1061)
	(1571,1052)(1556,1042)(1541,1032)
	(1528,1022)(1514,1012)(1502,1002)
	(1490,992)(1479,983)(1469,973)
	(1460,964)(1449,953)(1439,941)
	(1428,929)(1418,917)(1408,904)
	(1399,892)(1389,879)(1380,867)
	(1371,856)(1363,845)(1355,836)
	(1348,827)(1341,820)(1335,813)
	(1328,807)(1322,803)(1315,798)
	(1308,795)(1302,793)(1295,791)
	(1288,790)(1281,790)(1274,791)
	(1268,793)(1262,795)(1256,798)
	(1250,802)(1244,806)(1237,812)
	(1230,819)(1222,826)(1215,835)
	(1206,844)(1197,854)(1188,865)
	(1179,875)(1170,886)(1161,896)
	(1152,905)(1142,915)(1134,923)
	(1125,931)(1116,939)(1106,947)
	(1094,955)(1083,963)(1070,971)
	(1057,979)(1044,987)(1030,995)
	(1016,1002)(1002,1008)(988,1015)(973,1021)
\path(1693,774)(1693,773)(1695,770)
	(1700,762)(1708,750)(1719,733)
	(1730,713)(1742,693)(1753,674)
	(1762,656)(1770,640)(1777,624)
	(1783,610)(1787,596)(1791,582)
	(1794,567)(1796,552)(1798,537)
	(1801,522)(1803,506)(1806,491)
	(1809,477)(1813,464)(1817,452)
	(1823,442)(1830,435)(1838,429)
	(1847,425)(1857,424)(1868,425)
	(1880,427)(1894,432)(1909,437)
	(1926,445)(1943,453)(1962,463)
	(1982,474)(2001,486)(2021,499)
	(2040,511)(2059,524)(2077,537)
	(2094,549)(2110,560)(2125,572)
	(2141,584)(2156,596)(2171,608)
	(2185,621)(2200,634)(2216,649)
	(2233,664)(2249,679)(2265,694)
	(2278,707)(2288,717)(2294,723)(2297,726)
\put(547,910){\makebox(0,0)[lb]{\smash{{{\SetFigFont{5}{6.0}{\rmdefault}{\mddefault}{\updefault}{\bf 1}}}}}}
\end{picture}
}
\end{center}
  \caption{By smoothing the curve of Figure~\ref{nodalcurve}, we see
    that the its genus is $2$ \lremind{smoothing}}
\label{smoothing}\end{figure}

\noindent 
{\em Exercise (for those with enough background):}
If $C$ has $\de$ nodes, and its irreducible components
have geometric genus $g_1$, \dots, $g_k$ respectively, 
show that $\sum_{i=1}^k (g_i-1) + 1 + \de$.

We define the {\em dual graph} of a a pointed nodal curve as follows.
It consists of vertices, edges, and ``half-edges''.  The vertices
correspond to the irreducible components of the curve, and are labeled
with the {\em geometric}  genus of the component.  When the genus is $0$, 
the label will be
omitted for convenience.  The edges correspond to the nodes, and join
the corresponding vertices.  (Note that an edge can contain a vertex
to itself.)  The half-edges correspond to the labeled points.  The
dual graph corresponding to Figure~\ref{nodalcurve} is given in
Figure~\ref{dualexample}.

\begin{figure}
\begin{center}
\setlength{\unitlength}{0.00083333in}
\begingroup\makeatletter\ifx\SetFigFont\undefined%
\gdef\SetFigFont#1#2#3#4#5{%
  \reset@font\fontsize{#1}{#2pt}%
  \fontfamily{#3}\fontseries{#4}\fontshape{#5}%
  \selectfont}%
\fi\endgroup%
{\renewcommand{\dashlinestretch}{30}
\begin{picture}(2146,633)(0,-10)
\put(786,346){\ellipse{150}{150}}
\put(1536,346){\blacken\ellipse{74}{74}}
\put(1536,346){\ellipse{74}{74}}
\put(1836,309){\ellipse{604}{604}}
\path(711,346)(186,346)
\path(711,346)(186,346)
\path(861,346)(1536,346)
\path(861,346)(1536,346)
\put(762,323){\makebox(0,0)[lb]{\smash{{{\SetFigFont{5}{6.0}{\rmdefault}{\mddefault}{\updefault}1}}}}}
\put(0,319){\makebox(0,0)[lb]{\smash{{{\SetFigFont{5}{6.0}{\rmdefault}{\mddefault}{\updefault}{\bf 1}}}}}}
\end{picture}
}
\end{center}
\caption{The dual graph to the pointed nodal curve of Figure~\ref{nodalcurve}
(unlabeled vertices are genus $0$)}
\label{dualexample}
\end{figure}
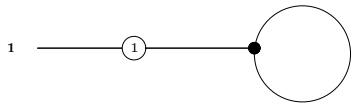

A nodal curve is said to be {\em stable} if it has finite automorphism
group.  This is equivalent to a combinatorial condition: (i) each
genus $0$ vertex of the dual graph has valence at least three, and
(iii) each genus $1$ vertex has valence at least one.  

\noindent 
{\em Exercise.}
Prove this.  You may use the fact that a genus $g \geq 2$ curve has
finite automorphism group, and that an elliptic curve (i.e.\ a $1$-pointed
genus $1$ curve) has finite automorphism group.  While you are proving
this, you may as well show that the automorphism group of a stable
genus $0$ curve is trivial.

\epoint{Exercise} Draw\label{m05}\lremind{m05} all possible stable
dual graphs for $g=0$ and $n \leq 5$; also for $g=1$ and $n \leq 2$.
In particular, show there are no stable dual graphs if $(g,n) =
(0,0)$, $(0,1)$, $(0,2)$, $(1,0)$.

\noindent {\bf Fact.}  There is a moduli space
of stable nodal curves of genus $g$ with $n$ marked points, denoted 
$\cmbar_{g,n}$.  There is an open subset
corresponding to smooth curves, denoted $\cm_{g,n}$.
The space $\cmbar_{g,n}$ is irreducible, of dimension $3g-3+n$,
and smooth.

(For Gromov-Witten experts:  you can interpret this space
as the moduli space of stable maps to a point.  But this in some sense
backwards, both historically, and in terms of the importance of both spaces.)

\noindent {\em Exercise.}  Show that $\chi(\cm_{g,n}) = (-1)^n \frac
{(2g+n-3)! B_{2g}} { 2g (2g-2)!}$, using the Harer-Zagier fact earlier
\eqref{hzeq}.

\epoint{Strata} To\label{stratadef}\lremind{stratadef} each stable
graph $\Ga$ of genus $g$ with $n$ points, we associate the subset
$\cm_{\Ga} \subset \cmbar_{g,n}$ of curves with that dual graph.  This
translates to the space of curves of a given topological type.  Notice
that if $\Ga$ is the dual graph given in Figure~\ref{dualexample}, we
can obtain any curve in $\cm_{\Ga}$ by taking a genus $0$ curve
with three marked points and gluing two of the points together, and
gluing the result to a genus $1$ curve with two marked points.  (This
is most clear in Figure~\ref{nodalcurve}.)  Thus each $\cm_{\Ga}$ is
naturally the quotient of a product of $\cm_{g',n'}$'s by some
symmetric group.  For example, if $\Ga$ is as in
Figure~\ref{dualexample}, $\cm_\Ga = (\cm_{0,3} \times \cm_{1,2}) /
S_2$.

These $\cm_{\Ga}$ give a stratification of $\cmbar_{g,n}$, and this
stratification is essentially as nice as one could hope.  For example,
the divisors (the closure of the codimension one strata) meet
transversely along smaller strata.  The dense open set $\cm_{g,n}$ is
one stratum; the rest are called {\em boundary strata}.
The codimension $1$ strata are called {\em boundary divisors}.

Notice that even if we were initially interested only in unpointed
Riemann surfaces, i.e.\ in the moduli space $\cm_g$, then
this compactification forces us to consider $\cm_{\Ga}$, which
in turn forces us to consider pointed nodal curves.

\noindent
{\em Exercise.}  By computing $\dim \cm_{\Ga}$, check that the 
codimension of the boundary stratum corresponding to a dual graph $\Ga$ is
precisely the number of edges of the dual graph.  (Do this first in some
easy case!)  

\epoint{Important exercise} Convince\label{m04}\lremind{m04} yourself
that $\cmbar_{0,4} \cong \proj^1$.  The isomorphism is given as
follows.  Given four distinct points $p_1$, $p_2$, $p_3$, $p_4$ on a
genus $0$ curve (isomorphic to $\proj^1$), we may take their
cross-ratio $\la = (p_4-p_1)(p_2-p_3) / (p_4-p_3)(p_2-p_1)$, and in
turn the cross-ratio determines the points $p_1$, \dots, $p_4$ up to
automorphisms of $\proj^1$.  The cross-ratio can take on any value in
$\proj^1- \{ 0, 1, \infty \}$.  The three $0$-dimensional strata
correspond to these three missing points --- figure out which
stratum corresponds to which of these three points.

\noindent
{\em Exercise.}  Write down the strata of $\cmbar_{0,5}$, along with which 
stratum is in the closure of which other stratum (cf.\ Exercise~\ref{m05}).

\bpoint{Natural morphisms among these moduli spaces}

We next describe some natural maps between these moduli spaces.  For
example, given any $n$-pointed genus $g$ curve (where $(g,n) \neq
(0,3), (1,1)$, $n >0$), we can forget the $n$th point, to obtain an
$(n-1)$-pointed nodal curve of genus $g$.  This curve may not be
stable, but it can be ``stabilized'' by contracting all components
that are $2$-pointed genus $0$ curves.  This gives us a map
$\cmbar_{g,n} \rightarrow \cmbar_{g,n-1}$, which we dub the {\em
  forgetful morphism}.

\noindent
{\em Exercise.}  Create an example of a dual graph where
stabilization is necessary.  Also, explain why we excluded the
cases $(g,n)=(0,3), (1,1)$.

\epoint{Important exercise}
Interpret\label{forgetfuluniversal}\lremind{forgetfuluniversal}
$\cmbar_{g,n+1} \rightarrow \cmbar_{g,n}$ as the universal curve over
$\cmbar_{g,n}$.  (This is a bit subtle.  Suppose $C$ is a nodal curve,
with node $p$. Which {\em stable} pointed curve with $1$ marked
point corresponds to $p$?  Similarly, suppose $(C,p)$ is a pointed
curve.  Which stable $2$-pointed curve corresponds to $p$?)

Given an $(n_1+1)$-pointed curve of genus $g_1$, and
an $(n_2+1)$-pointed curve of genus $g_2$, we 
can glue the first curve to the second along the last point of each,
resulting in an $(n_1+n_2)$-pointed curve of genus $g_1+g_2$.
This gives a map
$$
\cmbar_{g_1, n_1+1} \times \cmbar_{g_2, n_2+1} \rightarrow \cmbar_{g_1+g_2, n_1+n_2}.$$
Similarly, we could take a single $(n+2)$-pointed curve of genus $g$,
and glue its last two points together to get an $n$-pointed curve
of genus $g+1$; this gives a map
$$
\cmbar_{g,n+2} \rightarrow \cmbar_{g+1, n}.$$
We call these last two types of maps {\em gluing morphisms}.

We call the forgetful and gluing morphisms the {\em natural morphisms}
between moduli spaces of curves.

\section{Tautological cohomology classes on  moduli spaces of curves,
and their structure}

We now define some cohomology classes on these two sorts of moduli
spaces of curves, $\cm_g$ and $\cmbar_{g,n}$.  Clearly by Harer and
Zagier's Euler-characteristic calculation \eqref{hzeq}, we should
expect some interesting classes, and it is a challenge to name some.
Inside the cohomology ring, there is a subring, called the {\em
  tautological (sub)ring} of the cohomology ring, that consists
informally of the geometrically natural classes.  An equally informal
definition of the tautological ring is: all the classes you can easily
think of.  (Of course, this isn't a mathematical statement.  But we do
not know of a single algebraic class in $H^*(\cm_g)$ that can be
explicitly written down, that is provably not tautological, even
though we expect that they exist.)  Hence we care very much about this
subring.

The reader may work in cohomology, or in the Chow ring (the algebraic
analogue of cohomology).  The tautological elements will live
naturally in either, and the reader can choose what he or she is most
comfortable with.    In order to emphasize that one can work algebraically,
and also that our dimensions and codimensions are algebraic, I will
use the notation of the Chow ring $A^i$, but most readers
will prefer to interpret all statements in the cohomology ring.
There is a natural map $A^i \rightarrow H^{2i}$, and the reader
should be conscious of that doubling of the index.

If $\al$ is a $0$-cycle on a compact orbifold $X$, then $\int_X
\al$ is defined to be its degree.

\bpoint{Tautological classes on $\cm_g$, take one} \label{taut1}\lremind{taut1}

A good way of producing cohomology classes on $\cm_g$ is to take
Chern classes of some naturally defined vector bundles.

On the universal curve $\pi: \cC_g \rightarrow \cm_g$ over $\cm_g$,
there is a natural line bundle; on the fiber $C$ of $\cC_g$, it is the
line bundle of differentials $\cL$ of $C$.  Define $\psi := c_1(\cL)$, which
lies in $A^1(\cC_g)$ (or $H^2(\cm_g)$ --- but again, we will stick to the
language of $A^*$).  Then $\psi^{i+1} \in A^{i+1}(\cC_g)$, and as $\pi$
is a proper map, we can push this class forward to
$\cm_g$, to get the {\em Mumford-Morita-Miller $\ka$-class}
$$
\boxed{\ka_i := \pi_* \psi^{i+1}, \quad i= 0, 1, \dots.}$$

Another natural vector bundle is the following.
Each genus $g$ curve (i.e.\ each point of $\cm_g$)  has a $g$-dimensional space of
differentials (\S \ref{genus}), and the corresponding  rank $g$ vector bundle
on $\cm_g$ is called the {\em Hodge bundle}, denoted $\E$.
(It can also be defined by $\E := \pi_* \cL$.)  
We define the $\la$-classes by
$$
\boxed{\la_i := c_i(\E), \quad i = 0, \dots, g.}
$$

We define the {\em tautological ring} as the subring of the Chow ring
generated by the
$\ka$-classes.  (We will have another definition in \S \ref{tautdef}.)  This
ring is denoted $R^*(\cm_g)
\subset A^*(\cm_g)$ (or $R^*(\cm_g)
\subset H^{2*}(\cm_g)$).

It is a miraculous ``fact'' that everything else you can think of
seems to lie in this subring.  For example, the following generating
function identity determines the $\la$-classes from the $\ka$-classes
in an attractive way, and incidentally serves as an advertisement
for the fact that generating functions (with coefficients
in the Chow ring) are a good way to package information
\cite[p.~111]{faber}:
$$
\sum_{i=0}^\infty \la_i t^i = \exp \left( \sum_{i=1}^{\infty} 
  \frac { B_{2i} \ka_{2i-1} } {2i (2i-1)} t^{2i-1} \right).$$

\bpoint{Faber's conjectures} 

The study of the tautological ring was
begin in Mumford's fundamental paper \cite{mumford}, but there was no
reason to think that it was particularly well-behaved.  But just over
a decade ago, Carel Faber proposed a remarkable constellation of
conjectures (first in print in \cite{faber}), suggesting that the
tautological ring has a beautiful combinatorial structure.  It is
reasonable to state that Faber's conjectures have motivated a great
deal of the remarkable progress in understanding the topology of the
moduli space of curves over the last
decade.\label{s:fconj}\lremind{s:fconj}

Although Faber's conjectures deal just with the moduli of smooth 
curves, their creation required knowledge of the compactification, and even
of Gromov-Witten theory, as we will later see.

A good portion of Faber's conjectures can be informally summarized
as:  ``$R^*(\cm_g)$ behaves  like the ($(p,p)$-part of
the) cohomology ring of a $(g-2)$-dimensional complex projective manifold.''
We now describe (most of) Faber's conjectures more precisely.
I have chosen to cut them into three pieces.

\noindent
{\bf I. ``Vanishing/socle'' conjecture.}
$R^i(\cm_g) = 0$ for $i>g-2$, and $R^{g-2}(\cm_g) \cong \Q$.  This was
proved by Looijenga \cite{looijenga} and Faber \cite[Thm.~2]{faber}.
(Looijenga's theorem will be stated explicitly below, see
Theorem~\ref{lthm}.)  We will prove the ``vanishing'' part
$R^i(\cm_g)=0$ for $i>g-2$ in \S \ref{thmstarloo}, and show that
$R^{g-2}(\cm_g)$ is generated by a single element as a consequence of
Theorem~\ref{Zd1}.  These statements comprise Looijenga's theorem
(Theorem~\ref{lthm}).  The remaining part (that this generator
$R^{g-2}(\cm_g)$ is non-zero) is a theorem of Faber's, and we omit
its proof.

\noindent
{\bf II.  Perfect pairing conjecture.}    
The analog of Poincar\'{e} duality holds: for $0 \leq i \leq g-2$, the natural
product $R^i(\cm_g) \times R^{g-2-i}(\cm_g) \rightarrow R^{g-2}(\cm_g)
\cong \Q$ is a perfect pairing.
This conjecture is currently completely open, and is only known in 
special cases.

We call a ring satisfying {\bf I} and {\bf II} a {\em Poincar\'e
duality ring of dimension $g-2$}.

A little thought will convince you that thanks to {\bf II} if we knew
the ``top intersections'' (i.e.\ the products of $\ka$-classes of
total degree $g-2$, as a multiple of the generator of
$R^{g-2}(\cm_g)$), then we would know the complete structure of the
tautological ring.    Faber predicts the answer to this as well.

\noindent
{\bf III. Intersection number conjecture (take one).}
(We will give a better statement in Conjecture~\ref{taketwo}, in terms
of a partial compactification of $\cm_{g,n}$.)\lremind{fconj1}
For any $n$-tuple of non-negative integers $(d_1, \dots, d_n)$,
\begin{equation}\label{fconj1}
\frac { (2g-3+n)! (2g-1)!!} { (2g-1)! \prod_{j=1}^n (2d_j+1)!!}
\ka_{g-2} = \sum_{\si \in S_n} \ka_{\si}\end{equation}
where 
if $\si=(a_{1,1} \cdots a_{1, i_1})(a_{2,1}
\cdots a_{2, i_2}) \cdots$ is the cycle decomposition of $\si$,
then $\ka_{\si}$ is defined to be $\prod_j (d_{a_{j,1}}+d_{a_{j,2}}+
\cdots + d_{a_{j,i_j}})$.
 Recall that $(2k-1)!! = 1 \times 3 \times \cdots \times (2k-1) = (2k)! / 2^k k!$.

For example, we have
$$\ka_{i-1} \ka_{g-i-1} + \ka_{g-2} = \frac { (2g-1)!!} { (2i-1)!! (2g-2i-1)!!} \ka_{g-2}$$
and
$$
\ka_1^{g-2} = \frac 1 {g-1} 2^{2g-5}  (g-2)!^2 \ka_{g-2}.$$

Remarkably, Faber was able to deduce this elegant conjecture from
a very limited amount of experimental data.

Faber's intersection number conjecture begs an obvious question: why
is this formula so combinatorial?  What is the combinatorial structure
behind this ring?  Faber's alternate description of the intersection
number conjecture (Conjecture~\ref{taketwo}) will be even more
patently combinatorial.

Faber's intersection number conjecture is now a theorem.  Getzler and
Pandharipande showed that it is a formal consequence of the {\em
  Virasoro conjecture} for the projective plane \cite{gp}.  The
Virasoro conjecture is due to the physicists Eguchi, Hori, Xiong, and
also the mathematician Sheldon Katz, and deals with the Gromov-Witten
invariants to some space $X$.  (See \cite[Sect.~ 10.1.4]{coxkatz} for a
statement.)  Getzler and Pandharipande show that the Virasoro
conjecture in $\proj^2$ implies a recursion among the intersection
numbers on the (compact) moduli space of stable curves, which in turn
is equivalent to a recursion for the top intersections in Faber's
conjecture.  They then show that the recursions have a unique
solution, and that Faber's prediction is a solution.

Givental has announced a proof of Virasoro conjecture for projective space
(and more generally Fano toric varieties)
\cite{givental}.  The details of the proof have not appeared, but
Y.-P. Lee and Pandharipande are writing a book \cite{lp} giving the
details.  This theorem is really a tour-de-force, and the most
important result in Gromov-Witten theory in some time.
However, it seems  a round-about and high-powered way of proving
Faber's intersection number conjecture.  For example, by its nature, it
cannot shed light on the combinatorial structure behind the intersection
numbers.  For this reason, it seems worthwhile giving a more direct
argument.  At the end of these notes, I will outline a program
for tackling this conjecture (joint with the combinatorialists
I.P. Goulden and D.M. Jackson), and a proof in a large class of cases.

(There are two other conjectures in this constellation worth
mentioning.  Faber conjectures that $\ka_1$, \dots, $\ka_{[g/3]}$
generate the tautological ring, with no relations in degrees $\leq
[g/3]$.  Both Morita \cite{morita1} and Ionel \cite{i} have given
proofs of the first part of this conjecture a few years ago.
Faber also conjectures that 
$R^*(\cm_g)$ satisfies the Hard Lefschetz and Hodge Positivity  properties 
with respect to the
class $\ka_1$ \cite[Conj.~ 1(bis)]{faber}.

As evidence, Faber has checked that his conjectures hold true in genus
up to $21$ \cite{faberpc}. I should emphasize that this check is very
difficult to do --- the rings in question are quite large and
complicated!  Faber's verification involves some clever constructions,
and computer-aided computations.

Morita has recently announced a conjectural form of the tautological
ring, based on the representation theory of the symplectic group
$\operatorname{Sp}(2g,\Q)$ \cite[Conj.~1]{morita2}.  This is a new and
explicit (and attractive) proposed description of the tautological
ring.  One might hope that his conjecture may imply Faber's
conjecture, and may also be provable.

\bpoint{Tautological classes on $\cmbar_{g,n}$} \label{taut2}\lremind{taut2}

We can similarly define a tautological ring on the compact moduli
space of stable pointed curves, $\cmbar_{g,n}$.  In fact here the
definition is cleaner, and even sheds new light on the tautological ring of
$\cm_g$.  As before, this ring includes ``all classes one can easily
think of'', and as before, it will be most cleanly described in terms
of Chern classes of natural vector bundles.  Before we give a formal
definition, we begin by discussing some natural classes on
$\cmbar_{g,n}$.

\epoint{Strata} We note first that we have some obvious (co)homology
classes on $\cmbar_{g,n}$, that we didn't have on $\cm_g$: the
fundamental classes of the (closure of the) strata.  We will discuss
these classes and their relations at some length before moving on.

In genus $0$ (i.e., on $\cmbar_{0,n}$), the cohomology (and Chow) ring
is generated by these classes.  (The reason is that each stratum of
the boundary stratification is by (Zariski-)open subsets of affine
space.)  We will see why the tautological groups are generated by
strata in Exercise~\ref{g0g1}.

We thus have generators of the cohomology groups; it remains to find
the relations.  On $\cmbar_{0,4}$, the situation is especially nice.
We have checked that $\cmbar_{0,4}$ is isomorphic to $\proj^1$
(Exercise~\ref{m04}), and there are three boundary points.  They are
homotopic (as any two points on $\proj^1$ are homotopic) --- and even
{\em rationally equivalent}, the algebraic version of homotopic in the
theory of Chow groups.

By pulling back these relations by forgetful morphisms, and pushing
forward by gluing morphisms, we get many other relations for various
$\cmbar_{0,n}$.  We dub these {\em cross-ratio relations}, although
they go by many other names in the literature.  Keel has shown that
these are {\em all} the relations
\cite{keel}.\label{g0relns}\lremind{g0relns}

In genus $1$, the tautological ring (although not the cohomology or
Chow rings!) are again generated by strata.  (We will see why in
Exercise~\ref{g1boundary}, and again in Exercise~\ref{g0g1}.)\label{g1promise}\lremind{g1promise} We
again have cross-ratio relations, induced by a single
(algebraic/complex) codimension $1$ relation on $\cmbar_{0,4}$.
Getzler proved a new (codimension $2$) relation on $\cmbar_{1,4}$
\cite[Thm.~1.8]{getzler} (now known as {\em Getzler's relation}).  (It is
remarkable that this relation, on an important compact smooth
fourfold, parametrizing four points on elliptic curves, was discovered
so late.)  Via the natural morphisms, this induces relations on
$\cmbar_{1,n}$ for all $n$.  Some time ago, Getzler announced that
these two sorts of relations were the only relations among the strata
\cite[par.~2]{getzler}.

In genus $2$, there are very natural cohomology classes that are not
combination of strata, so it is now time to describe other
tautological classes.

\epoint{Other tautological classes}
Once again, we can define classes as 
Chern classes of natural vector bundles.

On $\cmbar_{g,n}$, for $1 \leq i \leq n$, we define the line bundle
$\L_i$ as follows.  On the universal curve $\cC_{g,n}
\rightarrow \cmbar_{g,n}$, the cotangent space at the fiber above
$[(C, p_1, \dots, p_n)] \in \cmbar_{g,n}$ at point $p_i$ is a
one-dimensional vector space, and this vector space varies smoothly
with $[(C, p_1, \dots, p_n)]$.  This is $\L_i$.  More precisely, if
$s_i: \cmbar_{g,n} \rightarrow \cC_{g,n}$ is the section of $\pi$
corresponding to the $i$th marked point, then $\cL_i$ is the pullback
by $s_i$ of the sheaf of relative differentials or the relative
dualizing sheaf (it doesn't matter which, as the section meets only
the smooth locus).  Define $\psi_i = c_1(\L_i) \in A^1(\cmbar_{g,n})$.

A genus $g$ nodal curve has a $g$-dimensional vector space
of sections of the dualizing line bundle.  These vector spaces vary
smoothly, yielding the {\em Hodge bundle} $\E_{g,n}$ on
$\cmbar_{g,n}$.  (More precisely, if $\pi$ is the universal curve over
$\cmbar_{g,n}$, and $\cK_{\pi}$ is the relative dualizing line bundle
on the universal curve, then $\E_{g,n} := \pi_* \cK_{\pi}$.)
Define $\la_i := c_i(\E_{g,n})$ on $\cmbar_{g,n}$.
Clearly the restriction of  the Hodge bundle
and $\la$-classes from $\cmbar_g$ to $\cm_g$ are the same
notions defined earlier.

Similarly, there is a more general definition of $\ka$-classes,
due to Arbarello and Cornalba \cite{ac}.

One might reasonably hope that these notions should behave well under
the forgetful morphism $\pi: 
\cmbar_{g,n+1} \rightarrow \cmbar_{g,n}$ 
(which we can interpret as the universal
curve by Exercise~\ref{forgetfuluniversal}).

\noindent {\em Exercise.}   Show that there is
a natural isomorphism $\pi^* \E_{g,n} \cong  \E_{g,n+1}$, and hence
that $\pi^* \la_k = \la_k$.

The behavior of the $\psi$-classes under pullback by the forgetful morphism
has a slight twist.

\tpoint{Comparison lemma}  {\em $\psi_1 = \pi^* \psi_1 + D_{0,\{ 1, n+1 \}}$.}
\label{comparison}\lremind{comparison}

(Caution:  the two $\psi_1$'s in the comparison lemma
are classes on two different spaces!)
Here $D_{0, \{ 1, n+1 \}}$ means the boundary divisor
corresponding to reducible curves with one node, where
one component is genus $0$ and contains only the marked points $p_1$
and $p_{n+1}$.
The analogous statement applies  with $1$ replaced by any number up to $n$
of course.

\noindent {\em Exercise (for people with more background).}  Prove the Comparison
lemma~\ref{comparison}.  (Hint:  First show that we have equality
away from $D_{0, \{ 1, n+1 \}}$.  Hence 
$\psi_1 = \pi^* \psi_1 + k D_{0,\{ 1, n+1 \}}$ for some integer $k$,
and this integer $k$ can be computed on a single test family.)

As an application:

\epoint{Exercise} Show that $\psi_1$ on $\cmbar_{0,4}$ is $\oh(1)$
(where $\cmbar_{0,4} \cong \proj^1$,
Exercise~\ref{m04}). \label{psi1}\lremind{psi1}

\noindent
{\em Exercise.} Express $\psi_1$
explicitly as a sum of boundary divisors on $\cmbar_{0,n}$.

We are now ready to define the tautological ring of $\cmbar_{g,n}$.
We do this by defining the rings for all $g$ and $n$ at once.

\epoint{Definition} The\label{tautdef}\lremind{tautdef} system of
tautological rings $( R^*(\cmbar_{g,n}) \subset
A^*(\cmbar_{g,n})_{g,n})$ is the smallest system of $\Q$-algebras
closed under pushforwards by the natural morphisms.

This elegant definition is due to Faber and Pandharipande \cite[\S 0.1]{fp3}.

Define the tautological ring of any open subset of $\cmbar_{g,n}$ by its
restriction from $\cmbar_{g,n}$.
In particular, we can recover our original definition of the
tautological ring of $\cm_g$ (\S \ref{taut1}). 

It is a surprising fact that everything else you can think of (such as
$\psi$-classes, $\la$-classes and $\ka$-classes) will lie in this ring.  (It is
immediate that fundamental classes of strata lie in this ring: they
are pushforwards of the fundamental classes of their ``component
spaces'', cf.\ \S \ref{stratadef}.)

We next give an equivalent description of the tautological {\em
  groups}, which will be convenient for many of our arguments, because
we do not need to make use of the multiplicative structure.  
In this description, the $\psi$-classes play a central role.

\epoint{Definition \cite[Defn.\ 4.2]{thmstar}} The system of
tautological rings $( R^*(\cmbar_{g,n}) \subset
A^*(\cmbar_{g,n})_{g,n})$ is the smallest system of $\Q$-vector spaces
closed under pushforwards by the natural morphisms, such that all
monomials in $\psi_1$, \dots, $\psi_n$ lie in $R^*(\cmbar_{g,n})$.
\label{ourtautdef}\lremind{ourtautdef}

The equivalence of Definition ~\ref{tautdef} and Definition~\ref{ourtautdef} is not difficult (see for
example \cite{thmstar}).

\bpoint{Faber-type conjectures
for $\cmbar_{g,n}$, and the conjecture of Hain-Looijenga-Faber-Pandharipande}

In analogy with Faber's conjecture, we have the following.

\bpoint{Conjecture}
$R^*(\cmbar_{g,n})$\label{hlfpconjecture}\lremind{hlfpconjecture} is a
Poincar\'e-duality ring of dimension $3g-3+n$.

This was first asked as a question by Hain and Looijenga
\cite[Question~5.5]{hl}, first stated as a speculation by Faber and
Pandharipande \cite[Speculation~3]{fplog} (in the case $n=0$), and first stated as a
conjecture by Pandharipande \cite[Conjecture~1]{icm}.  In analogy with
Faber's conjecture, we break this into two parts.

\noindent
{\bf I. ``Socle'' conjecture.}
$R^{3g-3+n}(\cmbar_{g,n}) \cong \Q$.
This is obvious if we define the tautological
ring in terms of  cohomology:  $H^{2(3g-3+n)}(\cmbar_{g,n}) \cong \Q$,
and the zero-dimensional strata show that the tautological zero-cycles
are not all zero.
However, in the tautological Chow ring, the socle conjecture is not at all obvious.
Moreover, the conjecture is not true in the full Chow ring ---
$A_0(\cmbar_{1,11})$ is uncountably generated, while
the conjecture states that $R_0(\cmbar_{1,11})$ has a
single generator. (By $R_0$, we of course mean $R^{3g-3+n}$.)

We will prove the  vanishing conjecture in \S \ref{soclepf}.

\noindent
{\bf II.  Perfect pairing conjecture}    
For $0\leq i \leq 3g-3+n$, the natural
product $$R^i(\cmbar_{g,n}) \times R^{3g-3+n-i}(\cmbar_{g,n}) \rightarrow R^{3g-3+n}(\cmbar_{g,n})
\cong \Q$$ is a perfect pairing.
(We currently have no idea why this should be true.)

Hence, in analogy with Faber's conjecture, if this conjecture were
true, then we could recover the entire ring by knowing the top
intersections.  This begs the question of how to compute all top
intersections.

\bpoint{Fact/recipe (Mumford and Faber)} \label{psithing}\lremind{psithing}If we
knew the top intersections of $\psi$-classes, we would know all top
intersections.  In other words, there is an algorithm to compute all
top intersections if we knew the numbers
\begin{equation}
\label{gtop}
\int_{\cmbar_{g,n}} \psi_1^{a_1} \cdots \psi_n^{a_n}, \quad \quad \sum a_i =
3g-3+n.\end{equation}

(This is a worthwhile exercise for people with some familiarity with
the moduli space of curves.)  
This is the basis of Faber's wonderful computer program \cite{fprog}
computing top intersections of various tautological classes.
For more information, see \cite{falg}.  This construction
is useful in understanding the definition (Defn.\ \ref{ourtautdef}) of
the tautological group in terms of the $\psi$-classes.

Until a key insight of Witten's, there was no a priori reason to
expect that these numbers should behave nicely.  We will survey
three methods of computing these numbers: (i) partial results in low
genus; (ii) Witten's conjecture; and (iii) via the ELSV formula.  A
fourth (attractive) method was given in Kevin Costello's thesis
\cite{costello}.

\epoint{Top intersections on $\cmbar_{g,n}$:  partial results in low
genus}
Here are two crucial relations among top intersections.

\noindent
{\bf Dilaton equation.} If $\cmbar_{g,n}$ exists
(i.e.\ there are stable $n$-pointed genus $g$ curve, or equivalently
$2g-2+n>0$), then
$$\int_{\cmbar_{g,n+1}} \psi_1^{\be_1} \cdots \psi_2^{\be_2}
\cdots \psi_n^{\be_n} \psi_{n+1} = (2g-2+n) \int_{\cmbar_{g,n}} \psi_1^{\be_1}
\cdots \psi_n^{\be_n}.$$ 

\noindent
{\bf String equation.} If $2g-2+n>0$, then 
$$\int_{\cmbar_{g,n+1}} \psi_1^{\be_1} \psi_2^{\be_2}
\cdots \psi_n^{\be_n}  = \sum_{i=1}^n \int_{\cmbar_{g,n}} \psi_1^{\be_1}  \psi_2^{\be_2}
\cdots \psi_i^{\be_i-1} \cdots \psi_n^{\be_n}$$
(where you ignore terms where you see negative exponents).

\noindent {\em Exercise (for those with more experience).}  Prove these using the
Comparison lemma~\ref{comparison}.

Equipped with the string equation alone, we can compute
all top intersections in genus $0$, i.e.
$\int_{\cmbar_{0,n}} \psi_1^{\be_1} \cdots \psi_n^{\be_n}$
where $\sum \be_i = n-3$.  (In any such expression, some $\be_i$ must
be $0$, so the string equation may be used.)
Thus we can recursively solve for these numbers, starting
from the base case
$\int_{\cmbar_{0,3}} \emptyset = 1$.

\noindent {\em Exercise.}
Show that $$\int_{\cmbar_{0,n}} \psi_1^{a_1} \cdots \psi_n^{a_n}  = \binom {n-3}
{a_1, \cdots, a_n}.$$

In genus $1$, the story is similar.  In this case, we need both the
string and dilaton equation.

\noindent {\em Exercise.}  Show that any integral
$$
\int_{\cmbar_{1,n}} \psi_1^{\be_1} \cdots \psi_n^{\be_n}$$ can be
computed using the string and dilaton equation from the base case
$\int_{\cmbar_{1,1}} \psi_1 = 1/24$.

We now sketch why the base case $\int_{\cmbar_{1,1}} \psi_1 = 1/24$ is
true.  We calculate this by choosing a finite cover $\proj^1
\rightarrow \cmbar_{1,1}$.  Consider a general pencil of cubics in the
projective plane.  In other words, take two general homogeneous cubic polynomials
$f$ and $g$ in three variables, and consider the linear combinations
of $f$ and $g$.  The non-zero  linear combinations modulo scalars are
parametrized by a $\proj^1$.  Thus we get a family of cubics parametrized
by 
$\proj^1$, i.e.\ $\cC \rightarrow \proj^1$.  

You can verify that in
this family, there will be twelve singular fibers, that are cubics
with one node.  One way of verifying this is as follows:  $f=g=0$
consists of nine points $p_1$, \dots, $p_9$ (basically by Bezout's
theorem --- you expect two cubics to meet at nine points).  There is a
map $\proj^2 - \{ p_1, \dots, p_9 \} \rightarrow \proj^1$.  If $\cC$
is the blow-up of $\proj^2$ at the nine points, then this map extends
to $\cC \rightarrow \proj^1$, and this is the total space of the
family.
The (topological) Euler characteristic of $\cC$ is the Euler
characteristic of $\proj^2$ (which is $3$) plus $9$ (as each blow-up
replaces a point by a $\proj^1$), i.e.\ $\chi(\cC) = 12$.  Considering
$\cC$ as a fibration over $\proj^1$, most fibers are elliptic curves,
which have Euler characteristic $0$.  Hence $\chi(\cC)$ is the sum of
the Euler characteristics of the singular fibers.  Each singular fiber
is a nodal cubic, which is isomorphic to $\proj^1$ with two points
glued together (depicted in Figure~\ref{delta0}); this is the union of $\C^*$ (which has Euler
characteristic $0$) with a point, so $\chi(\cC)$ is the number of
singular fibers.  (This argument needs further justification at every
point!)

We have a section of $\cC \rightarrow \proj^1$, given by the
exceptional fiber $E$ of the blow-up of $p_1$.  Hence we have a moduli map
$\mu: \proj^1 \rightarrow \cmbar_{1,1}$ of smooth curves.  Clearly it
doesn't map $\proj^1$ to a point, as some of the fibers are smooth,
and twelve are singular.  Thus the moduli map $\mu$ is surjective (as the image
is an irreducible closed set that is not a point).
You might suspect that $\mu$ has degree $12$, as the preimage of the
boundary divisor $\De \in \cmbar_{1,1}$ has $12$ preimages, and one
can check that $\mu$ is nonsingular here.  However, we come to one of
the twists of stack theory --- each point of $\cmbar_{1,1}$, including
$\De$, has degree $1/2$ --- each point should be counted with
multiplicity one over the size its automorphism group, and each
$1$-pointed genus $1$ stable curve has precisely one nontrivial
automorphism.

Thus $24 \int_{\cmbar_{1,1}} \psi_1 = \int_{\proj^1} \mu^* \psi_1$, so
we wish to show that $\int_{\proj^1} \mu^* \psi_1 = 1$.  This is an
explicit computation on $\cC \rightarrow \proj^1$.  You may check that
on the blow-up to $\cC$, the dualizing sheaf to the fiber at $p_1$ is
given by $-\oh(E)|_E$.  As $E^2 = -1$, we have $\int_{\proj^1} \mu^*
\psi_1 = -E^2=1$ as desired.

In higher genus, the string and dilaton equation are also very useful.

\noindent {\em Exercise.}  Fix $g$.  Show that using
the string and dilaton equation, all of the numbers \eqref{gtop}
(for all $n$) can be computed from a finite number of base cases.
The number of base cases required is the number of partitions
of $3g-3$.
(It is useful to describe this more precisely, by explicitly
describing the generating function for
\eqref{gtop} in terms of these base cases.)

\epoint{Witten's conjecture}
So how do we get at these remaining base cases?  The answer was
given by Witten \cite{w}.
(This presentation is not chronological --- Witten's conjecture
came first, and motivated most of what followed.  In particular, 
it predates Faber's conjectures, and was used to generate the data that
led Faber to his conjectures.)

\noindent
{\bf Witten's conjecture (Kontsevich's theorem).} 
Let $$F_g = \sum_{n \geq 0} \frac  1 {n!} \sum_{k_1, \dots, k_n} 
\left( \int_{\cmbar_{g,n}} \psi_1^{k_1} \cdots \psi_n^{k_n} \right)
t_{k_1} \cdots t_{k_n}$$
be the generating function for the genus $g$ numbers
\eqref{gtop}, and 
and let $$F = \sum F_g \hbar^{2g-2}$$
be the generating function for all genus.
(This is {\em Witten's free energy}, or the {\em Gromov-Witten
potential of a point}.)
Then
$$
(2n+1) \frac {\partial^3} {\partial t_n \partial t_0^2} F =
\left(  \frac {\partial^2} {\partial t_{n-1} \partial t_0} F \right)
\left( \frac {\partial^3} {\partial t^3_0} F \right)  +
2 \left(  \frac {\partial^3} {\partial t_{n-1} \partial t_0^2} F
\right)  \left( \frac {\partial^2 }  {\partial t_0^2} F \right)
+ \frac 1 4 \frac { \partial^5} {\partial t_{n-1} \partial t_0^4} F.
$$

Witten's conjecture now has many proofs, by Kontsevich \cite{kontsevich},
Okounkov-Pandharipande \cite{op}, Mirzakhani \cite{mirzakhani}, and Kim-Liu \cite{kimliu}.
It is a sign of the richness of this conjecture that these proofs
are all very different, and all very enlightening in different ways.

The reader should not worry about the details of this formula, and
should just look at its shape.  Those familiar with integrable systems
will recognize this as the Korteweg-de Vries (KdV) equation, in some guise.  There was a
later reformulation due to Dijkgraaf, Verlinde, and
Verlinde \cite{dvv}, in terms of the
Virasoro algebra.  Once again, the reader should not worry about the
precise statement, and concentrate on the form of the conjecture.
Define differential operators ($n \geq -1$)
\begin{eqnarray*}
L_{-1} &=& - \frac {\partial} {\partial t_0} + \frac {\hbar^{-2}} 2 t_0^2
+ \sum_{i=0}^{\infty}  t_{i+1} \frac \partial {\partial t_i}
\\
L_0 &=& - \frac 3 2 \frac \partial {\partial t_1} + \sum_{i=0}^{\infty}
\frac {2i+1} 2 t_i \frac {\partial} {\partial t_i} + \frac 1 {16}
\\
L_n &=&
\sum_{k=0}^{\infty} \frac { \Ga(m+n+ \frac 3 2) }{ \Ga(k+ \frac 1 2)} (t_k - \de_{k,1})
\delta_{n+k} + \frac {\hbar^2} 2 \sum_{k=1}^{n-1} (-1)^{k+1} \frac {\Ga(n-k+ \frac 1 2)}
{\Ga( -k-\frac 1 2)} \de_k \de_{n-k-1} \quad \quad (n>0)
\end{eqnarray*}
These operators satisfy
$[L_m, L_n] = (m-n) L_{m+n}$. 

% Caution!  L_0 has a 1/2 instead of a 3/2 in Getzler-Pandharipande p. 2. 
% GeP also uses
% h instead of  h^2.
% fror mirror book p. 516, we had 1/2.  I'm sticking with ours

\noindent {\em Exercise.}  Show that $L_{-1} e^F = 0$ is equivalent to
the string equation.  Show that $L_{0} e^F = 0$ is equivalent to
the dilaton equation. 

Witten's conjecture is equivalent to: $L_n e^F = 0$ for all $n$.
These equations let you inductively solve for the co-efficients of
$F$, and hence compute all these numbers.

\epoint{The Virasoro conjecture} 
The Virasoro formulation of Witten's conjecture has a far-reaching
generalization, the {\em Virasoro conjecture} described earlier.
Instead of top intersections on the moduli space of curves,
it addresses top (virtual) intersections on the moduli space
of maps of curves to some space $X$.  Givental's proof (to be
explicated by Lee and Pandharipande) for the case of projective
space (and more generally Fano toric varieties) was mentioned earlier.  It is also worth mentioning
Okounkov and Pandharipande's proof in the case where $X$ is a curve;
this is also a tour-de-force.

\epoint{Hurwitz numbers and the ELSV formula}
We can also recover these top intersections via
the old-fashioned theme of branched covers of the projective line,
the very technique that let us compute the dimension of the moduli
space of curves, and of the Picard variety \S \ref{riemann}.

Fix a genus $g$, a degree $d$, and a partition of $d$ into $n$ parts,
$\al_1 + \cdots + \al_n=d$, which we write as $\al \vdash d$.  Let\lremind{Hr}
\begin{equation}
r:=2g+d+n-2.\label{Hr}
\end{equation}  Fix $r+1$ points $p_1, \dots, p_r, \infty \in \proj^1$.
Define the {\em Hurwitz number} $H^g_{\al}$ to be the number of branched covers of $\proj^1$ by
a Riemann surface, that are unbranched away from $p_1, \dots, p_r,
\infty$, such that the branching over $\infty$ is given by $\al_1,
\dots, \al_n$ (i.e.\ there are $n$ preimages of $\infty$, and the
branching at the $i$th preimage is of order $\al_i$, i.e.\ the map is
analytically locally given by $t \mapsto t^{\al_i}$), and there is the
simplest possible branching over each $p_i$, i.e.\ the branching is
given by $2+1+\cdots+1=d$. 
(To describe this {\em simple branching} more explicitly:
above any such branch point, $d-2$ of the sheets are unbranched, and the remaining two sheets come together.  The analytic picture of the two sheets is the projection
of the parabola $y^2=x$ to the $x$-axis in $\C^2$.)
We consider the $n$ preimages of $\infty$
to be labeled.  Caution: in the literature, sometimes the preimages of
$\infty$ are {\em not labeled}; that definition of Hurwitz number will
be smaller than ours by a factor of $\# \Aut \al$, where $\Aut \al$ is
the subgroup of $S_n$ fixing the $n$-tuple $(\al_1, \dots, \al_n)$ (e.g.\ if $\al =
(2,2,2,5,5)$, then $\# \Aut \al = 3! 2!$).

One technical point:  each cover is counted with multiplicity
$1$ over the size of the automorphism group of the cover.

\noindent
{\em Exercise.}  Use the Riemann-Hurwitz formula \eqref{rh} to show
that if the cover is connected, then it has genus $g$.

Experts will recognize these as relative descendant Gromov-Witten
invariants of $\proj^1$; we will discuss relative Gromov-Witten
invariants of $\proj^1$ in Section~\ref{s:rvl}.  However, they are
something much more down-to-earth.  The following result shows that
this number is a purely combinatorial object.  In particular, there
are a finite number of such covers.

\tpoint{Proposition}
{\em 
$$H^g_{\al} = \# \left\{ (\si_1, \dots, \si_r) : \si_i \text{ transpositions generating
  $S_d$}, \prod_{i=1}^r \si_i \in C(\al) \right\} \# \Aut \al / d!,$$ where the
$\si_i$ are transpositions generating the symmetric group $S_d$, and
$\cC(\al)$ is the conjugacy class in $S_d$ corresponding to partition
$\al$.  } \label{hc}\lremind{hc}

Before we give the proof, we make some preliminary comments.  As an
example, consider $d=2$, $\al=2$, $g$ arbitrary, so $r=2g+1$.  The
above formula gives $H^g_{\al} = 1/2$, which at first blush seems like
nonsense --- how can we count covers and get a non-integer?  Remember
however the combinatorial/stack-theoretic principal that objects
should be counted with multiplicity  $1$ over the size of
their automorphism group.  Any double cover of this sort always has a
non-trivial involution (the ``hyperelliptic involution'').  Hence
there is indeed one cover, but it is counted as ``half a cover''.
Fortunately, this is the only case of Hurwitz numbers for which this
is an issue.  The reader may want to follow this particular case
through in the proof.

\noindent {\em Proof of Proposition~\ref{hc}.}
Pick another point $0 \in \proj^1$ distinct from $p_1, \dots,
p_r, \infty$.  Choose branch cuts from $0$ to $p_1, \dots, p_r,
\infty$ (non-intersecting paths from $0$ to $p_1$, $0$ to $p_2$, \dots, $0$
to $\infty$) such that their cyclic order around $0$ is $p_1, \dots,
p_r, \infty$.  Suppose $C \rightarrow \proj^1$ is one of the branched
covers counted by $H^g_{\al}$.  Then label the $d$ preimages of $0$ with  $1$
through $d$ in some way.  We will count these labeled covers, and
divide by $d!$ at the end.  Now cut along the preimages of the
branch-cuts.  As $\proj^1$ minus the branch-cuts is homeomorphic to a
disc, which is simply connected, its preimage must be $d$ copies of
the disc, labelled $1$ through $d$ according to the label on the
preimage of $0$.  We may reconstruct $C \rightarrow \proj^1$ by
determining how to glue these sheets together along the branch cuts.
The monodromy of the cover $C \rightarrow \proj^1$ around $p_i$ is an
element $\si_i$ of $S_d$, and this element will be a transposition,
corresponding to the two sheets being interchanged above that branch
point.  Similarly, the monodromy around $\infty$ is also an element
$\si_{\infty}$ of $S_d$, with cycle type $\al$.  The cover has the
additional data of the bijection of the cycles with the parts of
$\al$. 
 In $\pi_1( \proj^1 - \{ p_1, \dots, p_r, \infty \})$, the loops around $p_1$, \dots, $p_r$, $\infty$ multiply
to the identity, so $\si_1 \si_2 \cdots \si_r \si_{\infty} = e$.
(Here we use the fact that the branch cuts meet $0$ in this particular order.)
Thus $\si_{\infty}^{-1} = \si_1 \cdots  \si_r$.  This is the only
relation among these generators of $\pi_1(\proj^1 - \{ p_1, \dots, p_r, \infty \}$.
Furthermore,
the cover $C$ is connected, meaning that we can travel from any one of the $d$
sheets to any of the others, necessarily by travelling around the branch points.
This implies that the $\si_1$, \dots, $\si_r$, $\si_{\infty}$ (and hence
just the $\si_1$, \dots, $\si_r$) generate a transitive subgroup of $S_d$.
But the only transitive subgroup of $S_d$ containing a transposition $\si_1$
is all of $S_d$.

Conversely, given the data of transposition $\si_1$, \dots, $\si_r$ generating
$S_d$, with product of cycle type $\al$, along with a labelling
of the parts of the product (of which there are $\# \Aut \al$),
we can construct a connected cover $C \rightarrow \proj^1$, by the Riemann
existence theorem.  Thus, upon forgetting the labels of the $d$
sheets, we obtain the desired equality. \epf

The above proof clearly extends to deal with more general
Hurwitz numbers, where arbitrary branching is specified over
each of a number of points.

Proposition~\ref{hc} shows that any Hurwitz number may be readily
computed by hand or by computer.  What is interesting is the structure
behind them.  In 1891, Hurwitz \cite{hurwitz} showed that \lremind{hformula}
\begin{equation}
\label{hformula}
H^0_\al = r! d^{n-3} \prod \left(  \frac{ \al_i^{\al_i}}
  { \al_i!} \right).
\end{equation} By modern standards, he provided an outline of
a proof.  His work was forgotten by a large portion of the mathematics
community, and later people proved special cases, including D\'enes
\cite{denes} in the case $n=1$, Arnol'd \cite{arnold}
in the case $n=2$.  In the case $n=d$ (so $\al=1^d$) was stated by the
physicists Crescimanno and Taylor \cite{ct}, who apparently asked the
combinatorialist Richard Stanley about it, who in turn asked Goulden
and Jackson. Goulden Jackson independently discovered and proved Hurwitz'
original theorem in the mid-nineties \cite{gj}.  Since then, many
proofs have been given, including one by myself using moduli of curves \cite{g1}.

Goulden and Jackson studied the problem for higher genus, and
conjectured a structural formula for Hurwitz numbers in general.
Their polynomiality conjecture \cite[Conj.~1.2]{gjconj} implies the 
following.

\tpoint{Goulden-Jackson Polynomiality Conjecture (one version)} {\em For each $g,n$,
  there is a symmetric polynomial $P_{g,n}$ in $n$ variables, with
  monomials of homogeneous degree between $2g-3+n$ and $3g-3+n$,
   such that\lremind{pconj} \label{pconj}
$$
H^g_{\al}= r! \prod_{i=1}^n  
\left(  \frac{ \al_i^{\al_i}}
{ \al_i!} \right)
P_{g,n}(\al_1, \dots, \al_n).
$$
}

The reason this conjecture (and the original version) is true is an
amazing theorem of Ekedahl, Lando, M.~Shapiro, and Vainshtein.

\tpoint{Theorem (ELSV formula, by Ekedahl, Lando, M.~Shapiro, and Vainshtein
\cite{elsv1, elsv2})}\lremind{elsv}
{\em \begin{equation}\label{elsv}
H^g_{\al}= r! \prod_{i=1}^n  
\left(  \frac{ \al_i^{\al_i}}
{ \al_i!} \right)
\int_{\cmbar_{g,n}} \frac { 1 - \la_1 + \cdots + (-1)^g \la_g}
{ ( 1 - \al_1 \psi_1) \cdots (1 - \al_n \psi_n)}
\end{equation}
(if $\cmbar_{g,n}$ exists).}

We will give a proof in \S \ref{elsvproof}.

Here is how to interpret the right side of the equation.
Note that the $\al_i$ are integers, and the $\psi_i$'s and $\la_k$'s
are cohomology (or Chow) classes.
Formally invert the denominator, e.g.\
$$
\frac 1 {1 - \al_1 \psi_1}= 1 + \al_1 \psi_1 + \al_1^2 \psi_1^2 +
\cdots.$$ Then multiply everything out inside the integral sign, and
discard all but the summands of total codimension $3g-3+n$ (i.e.\
dimension $0$).  Then take the degree of this cohomology class.

For example, if $g=0$ and $n=4$, we get
\begin{eqnarray*}
H^g_{\al} &=&
 r! 
\prod_{i=1}^4  
\left(  \frac{ \al_i^{\al_i}}
{ \al_i!} \right)
\int_{\cmbar_{0,4}} \frac { 1 - \la_1 + \cdots \pm \la_g}
{ ( 1 - \al_1 \psi_1) \cdots (1 - \al_4 \psi_4)}\\
&=& 
r! 
\prod_{i=1}^4  
\left(  \frac{ \al_i^{\al_i}}
{ \al_i!} \right)
\int_{\cmbar_{0,4}} 
 ( 1 + \al_1 \psi_1+ \cdots ) \cdots (1 + \al_4 \psi_4 + \cdots)\\
&=&
r! 
\prod_{i=1}^4  
\left(  \frac{ \al_i^{\al_i}}
{ \al_i!} \right)
\int_{\cmbar_{0,4}} 
 ( \al_1 \psi_1+ \cdots + \al_4 \psi_4 ) \\
&=&
r! 
\prod_{i=1}^4  
\left(  \frac{ \al_i^{\al_i}}
{ \al_i!} \right)
(\al_1 + \cdots + \al_4)  \quad \quad \quad \text{(Exercise~\ref{psi1})}\\
&=&
r! 
\prod_{i=1}^4  
\left(  \frac{ \al_i^{\al_i}}
{ \al_i!} \right) d.
\end{eqnarray*}

\noindent
{\em Exercise.}  Recover Hurwitz' original formula \eqref{hformula} from the
ELSV-formula, at least if $n\geq 3$.

More generally, expanding the integrand  of \eqref{elsv} yields
\begin{equation}\label{Pgn}
 \sum_{a_1 + \cdots + a_n + k = 3g-3+n} \left( (-1)^k  
\left( \int_{\cmbar_{g,n}} \psi_1^{a_1} 
  \cdots \psi_n^{a_n} \la_k\right) \right)
\left( \al_1^{a_1} \cdots \al_n^{a_n} \right).
\end{equation}
This is a polynomial in $\al_1$, \dots, $\al_n$ of homogeneous
degree between
$2g-3+n$ and $3g-3+n$.
Thus this explains the mystery polynomial in the Goulden-Jackson Polynomiality
Conjecture~\ref{pconj} --- and the coefficients turn out to be top intersections
on the moduli space of curves!  (The original polynomiality conjecture
was actually different, and some translation is necessary in order to make
the connection with the ELSV formula \cite{gjv1}.)

There are many other consequences of the ELSV formula; see \cite{elsv2, gjv1} for surveys.

We should take a step back to see how remarkable the ELSV formula is.
To any reasonable mathematician, Hurwitz numbers (as defined by
Proposition \ref{hc}) are purely discrete, combinatorial objects.  Yet
their structure is fundamentally determined by topology of the moduli
space of curves.  Put more strikingly ---
the combinatorics of transpositions in the symmetric
group lead inexorably to the tautological ring of the moduli
space of curves!

\point We return to our original motivation for discussing the ELSV
formula: computing top intersections of $\psi$-classes on the moduli
space of curves $\cmbar_{g,n}$.  Fix $g$ and $n$.  As stated earlier,
any given Hurwitz number may be readily computed (and this can
be formalized elegantly in the language of generating functions).
Thus any number of values of $P_{g,n}(\al_1, \dots, \al_n)$ may be
computed.  However, we know that $P_{g,n}$ is a symmetric polynomial
of known degree, and it is straightforward to show that one can
determine the co-efficients of a polynomial of known degree from
enough values.  In particular, from \eqref{Pgn}, the coefficients of
the highest-degree terms in $P_{g,n}$ are precisely the top
intersections of $\psi$-classes.\label{trick}\lremind{trick}

This is a powerful perspective.  As an important example, Okounkov and
Pandharipande used the ELSV formula to prove Witten's conjecture.

\bpoint{Back to Faber-type conjectures}

This concludes our discussion of Faber-type conjectures for
$\cmbar_{g,n}$.  I have two more remarks about Faber-type conjectures.
The first is important, the second a side-remark.

\epoint{Faber's intersection number conjecture on $\cm_g$, take two}
We define the moduli space of $n$-pointed genus $g$ curves with
``rational tails'', denoted $\cm_{g,n}^{rt}$, as follows.  
We define $\cm_{g,n}^{rt}$ as the dense open subset
of $\cmbar_{g,n}$ parametrizing pointed nodal curves where one component
is nonsingular of genus $g$ (and the remaining components form trees
of genus $0$ curves sprouting from it --- hence the phrase ``rational
tails'').  If $g>1$, then $\cm_{g,n}^{rt} = \pi^{-1}(\cm_g)$, 
where $\pi: \cmbar_{g,n} \rightarrow \cmbar_g$ is the forgetful
morphism.
Note that $\cm_{g}^{rt}=\cm_g$.

We may restate Faber's intersection number conjecture (for $\cm_g$) in
terms of this moduli space.  By our re-definition of the tautological
ring on $\cm_g$ in \S \ref{taut2} (Definition~\ref{ourtautdef},
using also Faber's constructions of \S \ref{psithing}), the ``top intersections'' are
determined by $\pi_* \psi_1^{a_1} \cdots \psi_n^{a_n}$ (where $\pi: 
\cm_{g,n}^{rt} \rightarrow \cm_g$) for $\sum a_i = g-2+n$.

Then Faber's intersection number conjecture translates to the following.

\tpoint{Faber's intersection number conjecture (take two)}
{\em   If all $\al_i >1$, then 
$$\psi_1^{\al_1} \cdots \psi_n^{\al_n} = 
\frac { (2g-3+k)! (2g-1)!!} { (2g-1)! \prod_{j=1}^k (2d_j+1)!!}
[\text{generator}] \quad \text{for $\sum \al_i = g-2+n$}$$ where
$[\text{generator}] = \ka_{g-2} = \pi_* \psi_1^{g-1}$.}\label{taketwo}\lremind{taketwo}

(This reformulation is also due to Faber.)
This description is certainly more beautiful than the original one
\eqref{fconj1},
which suggests that we are closer to the reason for it to be true.

\point The other conjectures of Faber were extended to $\cm_{g,n}^{rt}$ by
Pandharipande \cite[Conj.\ 1]{icm}. \label{fabertypert}\lremind{fabertypert}

\epoint{Remark: Faber-type conjectures for curves of compact type}
Based on the cases of the $\cm_g$ and $\cm_{g,n}$, Faber and
Pandharipande made another conjecture for curves of ``compact type''.
A curve is said to be of {\em compact type} if its Jacobian is
compact, or equivalently if its dual graph has no loops, or
equivalently, if the curve has no nondisconnecting nodes.  Define
$\cm_{g,n}^c \subset \cmbar_{g,n}$ to be the moduli space of curves of
compact type.  It is $\cmbar_{g,n}$ minus an irreducible divisor,
corresponding to singular curves with one irreducible component (called
$\Delta_0$, although we will not use this notation).

\tpoint{Conjecture (Faber-Pandharipande \cite[Spec.\ 2]{fplog},
  \cite[Conj.\ 1]{icm}} {\em $R^*(\cm_g^c)$ is a Poincar\'e duality
  ring of dimension $2g-3$. }\label{ctconj}\lremind{ctconj}

Again, this has a {\em vanishing/socle} part and a {\em
perfect pairing} part.  There is something that can be considered
the corresponding {\em intersection number} part, Pandharipande and Faber's
$\lambda_g$ theorem \cite{lambdag}.

We will later (\S \ref{ctpf})  give a proof of the vanishing/socle portion
of the conjecture, that 
$R^i(\cm_g^c)=0$ for $i>2g-3$, and is  $1$-dimensional 
if $i=2g-3$.  The perfect pairing part is essentially completely open.

\bpoint{Other relations in the tautological ring}

We have been concentrating on top intersections in the tautological
ring.  I wish to discuss more about other relations (in smaller
codimension) in the tautological ring.

In genus $0$, as stated earlier (\S \ref{g0relns}), all classes on
$\cmbar_{g,n}$ are generated by the strata, and the only relation among them
are the cross-ratio relations.  We have also determined the
$\psi$-classes in terms of the boundary classes.

In genus $1$, we can verify that $\psi_1$ can be expressible in terms
of boundary strata.  On $\cmbar_{1,1}$,
if the boundary point is denoted $\de_0$ (the class of the nodal
elliptic curve shown in Figure~\ref{delta0}), we have shown $\psi_1 =
\de_0 / 12$.  (Reason: we proved it was true on a finite cover, in the
course of showing that $\int_{\cmbar_{1,1}} \psi_1 = 1/24$.)  
We know how to pull back $\psi$-classes by forgetful morphisms, so we
can now verify the following.

\noindent {\em Exercise.}  Show that in the cohomology group of
$\cmbar_{1,n}$, $\psi_i$ is equivalent to a linear combination of
boundary divisors.  (Hint:  use the Comparison Lemma \ref{comparison}.)

\begin{figure}
\begin{center}
\setlength{\unitlength}{0.00083333in}
\begingroup\makeatletter\ifx\SetFigFont\undefined%
\gdef\SetFigFont#1#2#3#4#5{%
  \reset@font\fontsize{#1}{#2pt}%
  \fontfamily{#3}\fontseries{#4}\fontshape{#5}%
  \selectfont}%
\fi\endgroup%
{\renewcommand{\dashlinestretch}{30}
\begin{picture}(1511,1323)(0,-10)
\put(790,1098){\blacken\ellipse{42}{42}}
\put(790,1098){\ellipse{42}{42}}
\path(513,1188)(555,1202)(598,1215)
	(639,1226)(678,1235)(714,1243)
	(748,1250)(778,1256)(806,1262)
	(831,1267)(853,1271)(873,1276)
	(891,1280)(908,1283)(924,1286)
	(939,1289)(954,1292)(970,1294)
	(986,1295)(1003,1296)(1022,1296)
	(1043,1294)(1065,1292)(1090,1287)
	(1117,1281)(1146,1273)(1177,1262)
	(1210,1248)(1244,1232)(1277,1211)
	(1309,1188)(1340,1159)(1368,1129)
	(1392,1097)(1412,1066)(1429,1037)
	(1444,1009)(1456,983)(1465,959)
	(1473,937)(1479,916)(1484,898)
	(1488,880)(1492,863)(1494,846)
	(1496,829)(1498,811)(1499,793)
	(1499,773)(1499,751)(1499,726)
	(1497,699)(1495,669)(1491,635)
	(1485,599)(1477,560)(1466,519)
	(1452,476)(1435,434)(1415,396)
	(1393,360)(1370,327)(1346,297)
	(1322,270)(1299,247)(1277,226)
	(1255,209)(1235,193)(1216,180)
	(1198,169)(1181,160)(1164,151)
	(1148,144)(1131,137)(1115,131)
	(1098,124)(1081,118)(1063,111)
	(1043,104)(1022,96)(1000,87)
	(976,78)(949,69)(921,58)
	(892,48)(861,38)(829,29)
	(796,21)(765,15)(733,12)
	(704,12)(678,15)(655,19)
	(634,25)(617,31)(602,38)
	(589,46)(578,53)(568,61)
	(560,68)(553,76)(547,84)
	(541,92)(536,100)(531,109)
	(527,118)(522,128)(518,139)
	(513,152)(509,165)(506,181)
	(503,198)(501,216)(500,237)
	(502,259)(506,283)(513,308)
	(526,338)(543,367)(563,394)
	(585,419)(609,442)(633,461)
	(658,478)(683,492)(708,505)
	(733,516)(757,526)(782,535)
	(806,544)(831,554)(855,564)
	(879,576)(903,589)(926,604)
	(948,621)(968,640)(986,660)
	(1001,683)(1011,705)(1016,727)
	(1015,747)(1008,764)(997,779)
	(983,792)(966,803)(948,812)
	(930,820)(910,827)(891,832)
	(871,837)(851,842)(831,846)
	(810,849)(790,852)(768,855)
	(746,856)(724,857)(700,857)
	(676,856)(651,852)(626,847)
	(601,838)(577,826)(555,811)
	(534,790)(517,766)(505,741)
	(498,717)(494,693)(494,670)
	(497,650)(502,630)(510,612)
	(518,594)(527,577)(536,560)
	(545,541)(552,522)(558,502)
	(562,479)(563,455)(560,428)
	(554,399)(544,369)(530,338)
	(513,308)(495,283)(475,261)
	(456,242)(437,226)(419,212)
	(403,200)(388,190)(374,182)
	(362,175)(350,169)(339,163)
	(328,158)(317,154)(306,150)
	(294,147)(281,144)(268,143)
	(252,142)(235,142)(217,145)
	(197,149)(177,157)(156,168)
	(136,183)(118,202)(102,225)
	(88,248)(76,270)(66,292)
	(57,313)(49,331)(43,348)
	(37,364)(32,379)(28,392)
	(24,406)(21,420)(18,435)
	(15,451)(13,469)(12,490)
	(12,514)(13,541)(16,572)
	(20,607)(28,645)(38,685)
	(52,727)(68,765)(85,801)
	(103,835)(120,865)(136,891)
	(150,914)(163,934)(173,950)
	(182,964)(190,976)(197,986)
	(203,995)(209,1003)(216,1011)
	(223,1019)(232,1028)(243,1037)
	(256,1048)(273,1061)(294,1075)
	(319,1092)(348,1110)(383,1129)
	(422,1149)(466,1169)(513,1188)
\put(886,1098){\makebox(0,0)[lb]{\smash{{{\SetFigFont{5}{6.0}{\rmdefault}{\mddefault}{\updefault}{\bf 1}}}}}}
\end{picture}
}
\end{center}
\caption{The curve corresponding to the point $\de_0 \in \cmbar_{1,1}$\lremind{delta0}}\label{delta0}
\end{figure}

\epoint{Slightly trickier exercise}
Use the above to show that the tautological ring in genus $1$ is
generated (as a group) by boundary classes.  (This fact was promised
in \S \ref{g1promise}.)  \label{g1boundary}\lremind{g1boundary}

In genus $2$, this is no longer true: $\psi_1$ is not equivalent to a
linear combination of boundary strata on $\cmbar_{2,1}$.  However, in
1983, Mumford showed that $\psi_1^2$ (on $\cmbar_{2,1}$) is a
combination of boundary strata (\cite{mumford}, see also \cite[eqn.\
(4)]{g2}); in 1998, Getzler showed the same for $\psi_1 \psi_2$ (on
$\cmbar_{2,2}$) \cite{g2}.  These two results can be used to show that
on $\cmbar_{2,n}$, all tautological classes are linear combinations of
strata, and from classes ``constructed using $\psi_1$ on
$\cmbar_{2,1}$''.  Figure~\ref{m21} may help elucidate what classes we
mean --- they correspond to dual graphs, with at most one marking
$\psi$ on an edge incident to one genus $2$ component.  The class in
question is defined by gluing together the class of $\psi_i$ on
$\cmbar_{2,v}$ corresponding to that genus $2$ component (where $v$ is
the valence, and $i$ corresponds to the edge labeled by $\psi$) with
the fundamental classes of the $\cmbar_{0,v_j}$'s corresponding to the
other vertices.  The question then arises: what are the relations
among these classes?  On top of the cross-ratio and Getzler relation,
there is a new relation due to Belorousski and Pandharipande, in
codimension $2$ on $\cmbar_{2,3}$ \cite{bp}.  We do not know if these
three relations generate all the relations.  (All the genus $2$
relations mentioned in this paragraph are given by explicit formulas,
although they are not pretty to look at.)
\begin{figure}
\begin{center}
\setlength{\unitlength}{0.00083333in}
\begingroup\makeatletter\ifx\SetFigFont\undefined%
\gdef\SetFigFont#1#2#3#4#5{%
  \reset@font\fontsize{#1}{#2pt}%
  \fontfamily{#3}\fontseries{#4}\fontshape{#5}%
  \selectfont}%
\fi\endgroup%
{\renewcommand{\dashlinestretch}{30}
\begin{picture}(2409,1409)(0,-10)
\put(1220,588){\ellipse{150}{150}}
\put(1820,1038){\blacken\ellipse{74}{74}}
\put(1820,1038){\ellipse{74}{74}}
\put(620,888){\blacken\ellipse{74}{74}}
\put(620,888){\ellipse{74}{74}}
\path(624,885)(1149,619)
\path(624,885)(1149,619)
\path(1825,1037)(1282,630)
\path(1825,1037)(1282,630)
\path(1238,511)(1304,154)
\path(1238,511)(1304,154)
\path(623,886)(280,424)
\path(623,886)(280,424)
\path(626,883)(128,1249)
\path(626,883)(128,1249)
\path(626,889)(1082,1251)
\path(626,889)(1082,1251)
\path(1819,1033)(2004,510)
\path(1819,1033)(2004,510)
\path(1821,1030)(2355,1144)
\path(1821,1030)(2355,1144)
\put(1196,565){\makebox(0,0)[lb]{\smash{{{\SetFigFont{5}{6.0}{\rmdefault}{\mddefault}{\updefault}2}}}}}
\put(1293,737){\makebox(0,0)[lb]{\smash{{{\SetFigFont{5}{6.0}{\rmdefault}{\mddefault}{\updefault}$\psi$}}}}}
\put(2007,380){\makebox(0,0)[lb]{\smash{{{\SetFigFont{5}{6.0}{\rmdefault}{\mddefault}{\updefault}{\bf 5}}}}}}
\put(2409,1135){\makebox(0,0)[lb]{\smash{{{\SetFigFont{5}{6.0}{\rmdefault}{\mddefault}{\updefault}{\bf 4}}}}}}
\put(1290,31){\makebox(0,0)[lb]{\smash{{{\SetFigFont{5}{6.0}{\rmdefault}{\mddefault}{\updefault}{\bf 2}}}}}}
\put(1100,1303){\makebox(0,0)[lb]{\smash{{{\SetFigFont{5}{6.0}{\rmdefault}{\mddefault}{\updefault}{\bf 6}}}}}}
\put(141,310){\makebox(0,0)[lb]{\smash{{{\SetFigFont{5}{6.0}{\rmdefault}{\mddefault}{\updefault}{\bf 1}}}}}}
\put(0,1310){\makebox(0,0)[lb]{\smash{{{\SetFigFont{5}{6.0}{\rmdefault}{\mddefault}{\updefault}{\bf 3}}}}}}
\end{picture}
}
\end{center}
\caption{A class on $\cmbar_{2,6}$ --- a codimension $1$ class
on a boundary stratum, constructed using $\psi_1$ on $\cmbar_{2,3}$ and gluing
morphisms\lremind{m21}}\label{m21}
\end{figure}
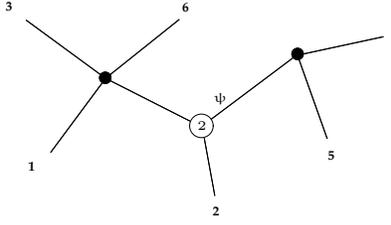

In general genus, the situation should get asymptotically worse
as $g$ grows.  However, there is a general statement that can be made:

\tpoint{Getzler's conjecture \cite[footnote~1]{g2} (Ionel's theorem
  \cite{i1})} {\em If $g>0$, all degree $g$ polynomials in
  $\psi$-classes vanish on $\cm_{g,n}$ (hence live on the boundary on
  $\cmbar_{g,n}$).} \label{getzlerconj}\lremind{getzlerconj}

We will interpret this result as a special case of a more general
result (Theorem $\star$), in  \S \ref{thmstargetzler}.  In keeping with the theme of this article,
the proof will  be Gromov-Witten theoretic.

\epoint{Y.-P. Lee's Invariance Conjecture}
There is another general statement that may well give {\em all} the
relations in {\em every} genus: Y.-P. Lee's {Invariance
  conjecture}.  It is certainly currently
beyond our current ability to either prove it.
Lee's conjecture is strongly motivated by
Gromov-Witten theory.\label{sic}\lremind{sic}

Before we state the conjecture, we discuss the
consequences and evidence.  All of the known relations in the
tautological rings are consequences of the conjecture.  For example,
the genus $2$ implications are shown by Arcara and Lee in \cite{al2}.
They then predicted a {\em new} relation in $\cmbar_{3,1}$ in
\cite{al3}.  Simultaneously and independently, this relation was
proved by Kimura and X. Liu \cite{kl}. This seems to be good evidence
for the conjecture being true.

More recently, the methods behind the conjecture have allowed Lee to
turn these predictions into proofs, {\em not conditional on the truth
  of the conjecture} \cite{lee2}.  Thus for example Arcara and Lee's
work yields a proof of the new relation on $\cmbar_{3,1}$.

We now give the statement.  The conjecture is most naturally expressed
in terms of the tautological rings of {\em possibly-disconnected}
curves.  The definition of a stable possibly-disconnected curve is the
same as that of a stable curve, except the curve is not required to be
connected.  We denote the moduli space of $n$-pointed genus $g$
possibly-disconnected curves by $\cmbar_{g,n}^{\dis}$.  The reader can
quickly verify that our discussion of the moduli space of curves
carries over without change if we consider possibly-disconnected
curves.  For example, $\cmbar_{g,n}^{\dis}$ is nonsingular and
pure-dimensional of dimension $3g-3+n$ (although not in general
irreducible).  It contains $\cmbar_{g,n}$ as a component, so
any statements about $\cmbar_{g,n}^{\dis}$ will imply statements about
$\cmbar_{g,n}$.  Note that the disjoint union of two curves of
arithmetic genus $g$ and $h$ is a curve of arithmetic genus $g+h-1$:
Euler characteristics add under disjoint unions.  Note also that a
possibly-disconnected marked curve is stable if and only if all of its
connected components are stable.

\noindent {\em Exercise.}
Show that $\cmbar_{-1,6}^{\dis}$ is a union of $\binom 6 3 / 2$ points
--- any $6$-pointed genus $-1$ stable curve must be the disjoint union
of two $\proj^1$'s, with $3$ of the $6$ labeled points on each component.

\noindent {\em Exercise.} Show that any component of $\cmbar_{g,n}^{\dis}$
is the quotient of a product of $\cmbar_{g',n'}$'s by a finite group.

Tautological classes are generated by classes corresponding to a dual
graph, with each vertex (of genus $g$ and valence $n$, say) labeled by
some cohomology class on $\cmbar^{\dis}_{g,n}$ (possibly the
fundamental class); call this a decorated dual graph.  (We saw an
example of a decorated dual graph in Figure~\ref{m21}.  Note that
$\psi$-classes will always be associated to some half edge.)  Decorated dual
graphs are not required to be connected.  If $\Ga$ is a decorated dual
graph (of genus $g$ with $n$ tails, say), let $\dim \Ga$ be the
dimension of the corresponding class in $A_*(\cmbar_{g,n}^{\dis})$.

For each positive integer $l$, we will describe a linear operator $\fr_l$
that sends formal linear combinations of decorated dual graphs to
formal linear combinations of decorated dual graphs. 
It is homogeneous of degree $-l$: it sends (dual graphs corresponding
to) dimension $k$ classes to (dual graphs corresponding to) dimension
$k-l$ classes.  

We now describe its action on a single decorated dual graph $\Ga$ of
genus $g$ with $n$ marked points (or half-edges), labeled $1$ through
$n$.  Then $\fr_l (\Ga)$ will be a formal linear combination of other graphs,
each of genus $g-1$ with $n+2$ marked points.

There are three types of contributions to $\fr_l \Ga$.
(In each case, we discard any graph that is not stable.)

\noindent {\em 1. Edge-cutting.}
There are two contributions for each {\em directed} edge, i.e.\
an edge with chosen starting and ending point.
(Caution:  there are two possible directions for each edge in general,
except for those edges that are ``loops'', connecting a
single vertex to itself.  In this case, both directions are considered
the same.)
We cut the edge, regarding the two half-edges as ``tails'', or
marked points.  The starting half-edge
is labeled $n+1$, and the ending half-edge is labeled $n+2$.
One summand will correspond to adding an extra
decoration of $\psi^l$ to point $n+1$. 
(In other words, $\psi_{n+1}^l$ is multiplied by whatever
cohomology class is already decorating that vertex.)
A second summand will correspond to the adding an extra
decoration of $\psi^l$ to point $n+2$, and this summand
appears with multiplicity $(-1)^{l-1}$.

\noindent {\em 2. Genus reduction} 
  For each vertex we produce $l$ graphs as follows.  We reduce the
  genus of the vertex by $1$, and add two new tails to this vertex, 
  labelled $n+1$ and $n+2$; we decorate them with 
  $\psi^m$ and $\psi^{l-1-m}$ respectively, where $0 \leq m \leq
  l-1$. Each such graph is taken with 
  multiplicity $(-1)^{m+1}$.

\noindent {\em 3. Vertex-splitting.}
  For each vertex, we produce a number of graphs as follows.
We split the vertex into two, and the first new vertex is
given the tail $n+1$, and the second is given the tail $n+2$.
The two new tails are decorated by 
$\psi^m$ and $\psi^{l-1-m}$ respectively, where
$0 \leq m \leq l-1$.  We then take one such graph for each choice
of splitting of the genus $g=g_1+g_2$  and partitioning of  
the other incident edges.
Each such graph is taken with multiplicity  $(-1)^{m+1}$.

Then $\fr_l(\Ga)$ is the sum of the above summands.  Observe that
when $l$ is odd (resp.\ even), the result is symmetric (resp.\ anti-symmetric)
in labels $n+1$ and $n+2$.

By linearity, this defines the action of $\fr_l$ on any linear combination of
directed graphs.

\tpoint{Y.-P. Lee's Invariance Conjecture \cite[Conj.~1--2]{yplee}} {\em 
\label{ic}\lremind{ic}
\begin{enumerate}
\item[(a)] If
  $\sum c_i \Ga_i =0$ holds in 
$A^*(\cmbar_{g,n}^{\dis})$, then
$\fr_l(\sum c_i \Ga_i) = 0$ in $A^*(\cmbar_{g-1,n+2})$.  
\item[(b)] Conversely, if $\sum c_i \Ga_i$ has {\bf positive} pure dimension, and
$\fr_l(\sum c_i \Ga_i) = 0$ in $A^*(\cmbar_{g-1,n+2})$, then
  $\sum c_i \Ga_i =0$ holds in 
$A^*(\cmbar_{g,n}^{\dis})$.
\end{enumerate}}

This can be used to produce tautological equations inductively!  The
base case is when $\dim \cmbar_{g,n}^{\dis} = 0$, which is known: we
will soon show (\S \ref{soclepf}) that $R_0(\cmbar_{g,n}) \cong \Z$,
and hence dimension $0$ tautological classes on $\cmbar_{g,n}$ are
determined by their degree (and dimension $0$ tautological classes on
$\cmbar_{g,n}^{\dis}$ are determined by their degree on each connected
component).  Note that the algorithm is a finite process: the
dimension $l$ relations on $\cmbar_{g,n}$ or $\cmbar_{g,n}^{\dis}$
produced by this algorithm are produced after a finite number of
steps.

Even more remarkably,
this seems to produce {\em all}  tautological relations:

\tpoint{Y.-P. Lee's Invariance Conjecture, continued \cite[Conj.~3]{yplee}}
{\em Conjecture~\ref{ic}(b) will produce {\bf all} tautological equations
inductively.}

A couple of remarks are in order.  Clearly this is a very
combinatorial description.  It was dictated by Gromov-Witten theory,
as explained in \cite{yplee}.  In particular, it uses the fact that all
tautological equations are invariant under the action of lower
triangular subgroups of the twisted loop groups,
and proposes that they are the only equations invariant in this way.

In order to see the magic of this conjecture in action, and to get
experience with the $\fr_l$ operators, it is best to work out an
example.  The simplest dimension $1$ relation is the following.

\noindent {\em Exercise.}
Show that the pullback of the (dimension $0$) cross-ratio relation 
(\S \ref{g0relns}) on
$\cmbar_{0,4}$ to a (dimension $1$) relation on $\cmbar_{0,5}$ is
implied by the Invariance Conjecture.  (Some rather beautiful cancellation
happens.)

\section{A blunt tool:  Theorem $\star$ and consequences}

We now describe a blunt tool from which much  of the previously
described structure of the tautological ring follows.
Although it is statement purely about the stratification of the moduli
space of curves, we will see (\S \ref{thmstarpf}) that it is proved via Gromov-Witten theory.

\tpoint{Theorem $\star$ \cite{thmstar}} 
{\em  Any tautological class of codimension $i$
is trivial away from strata satisfying\label{thmstarhere}\lremind{thmstarhere}
$$
\boxed{{\text {\# genus $0$ vertices  }} \geq i-g+1.}
$$}

(Recall that the genus $0$ vertices correspond to components
of the curve with geometric genus $0$.)

More precisely, any tautological class is zero upon restriction to the
(large) open set corresponding to 
the  open set corresponding to
$${\text {\# genus $0$ vertices  }} < i-g+1.$$
Put another way:  given any tautological class of codimension $i$,
you can move it into the set of curves with at least $i-g+1$ genus $0$ components.
A third formulation is that the tautological classes of codimension $i$
are pushed forward from classes on the locus of curves with
at least $i-g+1$ genus $0$ components.

We remark that this is false for the Chow ring as a whole ---
this is fundamentally a statement about tautological classes.

We will discuss the proof in \S \ref{thmstarpf}, but first we give consequences.
There are in some sense four morals of this result.

First, this is the fundamental geometry behind many of the theorems
we have been discussing.  We will see that  they follow from Theorem $\star$ by
straightforward combinatorics.  As a sign of this, we will
often get strengthenings of what was known or conjectured previously.

Second, this suggests the potential importance of a filtration of the
moduli space by number of genus $0$ curves.  It would be interesting
to see if this filtration really is fundamental, for example if it
ends up being relevant in understanding the moduli space of curves in
another way.  So far this has not been the case.

Third, as we will see from the proof, once one knows
a clean statement of what one wants to prove, the proof is relatively
straightforward,  at least in outline.

And fourth, the proof will once again show the centrality of
Gromov-Witten theory to the study of the moduli of curves.

\bpoint{Consequences of Theorem $\star$}

We begin with a warm-up example.

\epoint{Theorem $\star$ implies Getzler's conjecture ~\ref{getzlerconj}
  (Ionel's theorem)}  Any degree $g$ monomial is a codimension $g$
tautological class, which vanishes on the open set of $\cmbar_{g,n}$
corresponding to curves with no genus $0$ components.  If $g>0$, this
is non-empty and includes $\cm_{g,n}$.\label{thmstargetzler}\lremind{thmstargetzler}

In particular:  (1)  we get a proof of Getzler's conjecture;
(2) we see that more classes vanish on this set --- all tautological
classes of degree at least $g$, not just polynomials in the $\psi$-classes;
(3) we observe that the classes vanish on a bigger set than
$\cm_{g,n}$, and that what is relevant is not the smoothness of the
curves, but the fact that they have no genus $0$ components.
(4)  This gives a moral reason for Getzler's conjecture
not to hold in genus $0$.

\epoint{Theorem $\star$ implies the first part of  Looijenga's Theorem (Faber's
vanishing conjecture)}
Recall (\S \ref{s:fconj}) that Looijenga's Theorem is
part of the ``vanishing'' part of Faber's conjectures:
\label{thmstarloo}\lremind{thmstarloo}

\tpoint{Theorem \cite{looijenga}} {\em We have $R^i(\cm_g) = 0$ for $i>g-2$,
  $\dim R^{g-2}(\cm_g) \leq 1$.}\label{lthm}\lremind{lthm}

We will show that Theorem $\star$ implies the first part now;
we will show the second part as a consequence of Theorem~\ref{Zd1}.

First, if the codimension of a tautological class is greater than or
equal to $g$, then we get vanishing on the open set where there are no
genus $0$ components, so we get vanishing for the same reason as
Getzler's conjecture. 

The case of codimension $g-1$ is more subtle.  
From the definition of the tautological ring,
tautological classes are obtained by taking $\psi$-classes,
and multiplying, gluing, and pushing forward by forgetful morphisms.
Now on $\cm_g = \cm_{g,0}$, there are no $\psi$-classes, no boundary 
strata, and no tautological classes of codimension less than $g-1$.
hence all codimension $g-1$ tautological classes on $\cm_{g,0}$ are pushed
forward from tautological classes on $\cm_{g,1}$, which are necessarily of
codimension $g$.  These also vanish by Theorem $\star$ by the same argument as before.

As before, one can say more: 

\noindent {\em Exercise.}  Extend this argument to the moduli space of
curves with rational tails we can extend Looijenga's theorem to the
moduli space of curves with rational tails $\cm_{g,n}^{rt}$.  (First
determine the dimension of the conjectural Poincar\'e duality ring!)

The Faber-type conjecture for this space was mentioned in \S
\ref{fabertypert}.  I should point out that I expect that Looijenga's
proof extends to this case without problem, but I haven't checked.

\epoint{Theorem $\star$ implies the socle part of
  Hain-Looijenga-Faber-Pandharipande conjecture~\ref{hlfpconjecture}
  on $\cmbar_{g,n}$} Recall the socle part of the Hain-Looijenga-Faber-Pandharipande
conjecture ~\ref{hlfpconjecture}, that $R_0(\cmbar_{g,n}) \cong
\Q$.\label{soclepf}\lremind{soclepf} (We write
$R^{3g-3+n}(\cmbar_{g,n})$ as $R_0(\cmbar_{g,n})$ to remind the reader
that the statement is about tautological $0$-cycles.)

We show how this is implied by Theorem $\star$.  
This was first shown in 
\cite{socle}, which can be seen as a first step toward the
statement and proof of Theorem $\star$.

Our goal is to show
that all tautological $0$-cycles are commensurate, and that one of
them is non-zero.  Clearly the latter is true, as the class of a
$0$-dimensional stratum (a point) is tautological, and is non-zero as
it has non-zero degree, so we concentrate on the first statement.

By Theorem $\star$, any dimension $0$
tautological class is pushed forward from the locus  of curves
with at least $(3g-3+n)-g+1=2g-2+n$ genus $0$ components.

\noindent {\em Exercise.}
Show that  the only stable dual graphs with $2g-2+n$ genus $0$
components has all vertices genus $0$ and trivalent. 
Show that these are the $0$-dimensional strata. (See Figure~\ref{m12}
for the $0$-dimensional strata on $\cmbar_{1,2}$.)
\begin{figure}
\begin{center}
\setlength{\unitlength}{0.00083333in}
\begingroup\makeatletter\ifx\SetFigFont\undefined%
\gdef\SetFigFont#1#2#3#4#5{%
  \reset@font\fontsize{#1}{#2pt}%
  \fontfamily{#3}\fontseries{#4}\fontshape{#5}%
  \selectfont}%
\fi\endgroup%
{\renewcommand{\dashlinestretch}{30}
\begin{picture}(3724,629)(0,-10)
\put(636,307){\blacken\ellipse{74}{74}}
\put(636,307){\ellipse{74}{74}}
\put(1086,307){\blacken\ellipse{74}{74}}
\put(1086,307){\ellipse{74}{74}}
\put(3186,307){\blacken\ellipse{74}{74}}
\put(3186,307){\ellipse{74}{74}}
\put(2736,307){\blacken\ellipse{74}{74}}
\put(2736,307){\ellipse{74}{74}}
\put(2436,307){\ellipse{600}{600}}
\path(186,307)(1536,307)
\path(186,307)(1536,307)
\path(2736,307)(3186,307)
\path(3186,307)(3636,457)
\path(3186,307)(3636,157)
\put(3704,490){\makebox(0,0)[lb]{\smash{{{\SetFigFont{5}{6.0}{\rmdefault}{\mddefault}{\updefault}{\bf 1}}}}}}
\put(3724,64){\makebox(0,0)[lb]{\smash{{{\SetFigFont{5}{6.0}{\rmdefault}{\mddefault}{\updefault}{\bf 2}}}}}}
\put(0,268){\makebox(0,0)[lb]{\smash{{{\SetFigFont{5}{6.0}{\rmdefault}{\mddefault}{\updefault}{\bf 1}}}}}}
\put(1600,274){\makebox(0,0)[lb]{\smash{{{\SetFigFont{5}{6.0}{\rmdefault}{\mddefault}{\updefault}{\bf 2}}}}}}
\end{picture}
}
\end{center}
\caption{The $0$-dimensional strata on $\cmbar_{1,2}$ --- notice
that all vertices are genus $0$ and trivalent, and that there are $2g-2+n$
of them\lremind{m12}}\label{m12}
\end{figure}
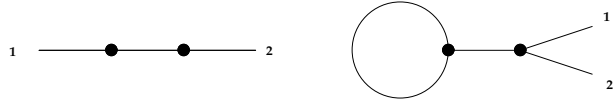

Hence $R_0(\cmbar_{g,n})$ is generated by these finite number of points.
It remains to show that any two of these points are equivalent
in the Chow ring.  A geometric way of showing this is by observing
that all points in $\cmbar_{0,N}$ are equivalent in the Chow ring, and
that our $0$-dimensional strata are in the image of $\cmbar_{0,2g+n}$ under
$2g$ gluing morphisms.  A more combinatorial way of showing this
is by showing that each $1$-dimensional stratum is isomorphic to $\proj^1$,
and that any two $0$-dimensional strata can be connected by a chain
of $1$-dimensional strata.

\noindent {\em Exercise.}
Complete one of these arguments.

As in the earlier applications of Theorem $\star$ too: we can verify
that the perfect pairing conjecture in codimension 1 and probably 2
(although Tom Graber and I haven't delved too deeply into 2).  This is
combinatorially more serious, but not technically hard.

\epoint{Theorem $\star$ implies the Faber-Pandharipande vanishing/socle
  conjecture on curves of compact type} We now show the ``vanishing/socle
part'' of the Faber-type conjecture for curves of compact type
(Faber-Pandharipande Conjecture~\ref{ctconj}).
\label{ctpf}\lremind{ctpf}

First, suppose that $i>2g-3$.  We will show that $R^i(\cm_g^c) = 0$.
By Theorem $\star$, any such tautological class vanishes on the
open set where there are at most $i-g+1 > g-2$ genus $0$ vertices.  
Then our goal follows from the next exercise.

\noindent {\em Exercise.}
Show that any genus $g$ ($0$-pointed) stable graph that is a
tree has at most $g-2$ genus $0$ vertices.
Moreover, if equality holds, then each vertex is either genus $1$
of valence $1$, or genus $0$ of valence $3$.  (Examples when
$g=6$ are given in Figure~\ref{g6}.)
\begin{figure}
\begin{center}
\setlength{\unitlength}{0.00083333in}
\begingroup\makeatletter\ifx\SetFigFont\undefined%
\gdef\SetFigFont#1#2#3#4#5{%
  \reset@font\fontsize{#1}{#2pt}%
  \fontfamily{#3}\fontseries{#4}\fontshape{#5}%
  \selectfont}%
\fi\endgroup%
{\renewcommand{\dashlinestretch}{30}
\begin{picture}(3766,1231)(0,-10)
\put(1508,308){\ellipse{150}{150}}
\put(1484,285){\makebox(0,0)[lb]{\smash{{{\SetFigFont{5}{6.0}{\rmdefault}{\mddefault}{\updefault}1}}}}}
\put(383,1133){\ellipse{150}{150}}
\put(359,1110){\makebox(0,0)[lb]{\smash{{{\SetFigFont{5}{6.0}{\rmdefault}{\mddefault}{\updefault}1}}}}}
\put(83,683){\ellipse{150}{150}}
\put(59,660){\makebox(0,0)[lb]{\smash{{{\SetFigFont{5}{6.0}{\rmdefault}{\mddefault}{\updefault}1}}}}}
\put(533,233){\ellipse{150}{150}}
\put(509,210){\makebox(0,0)[lb]{\smash{{{\SetFigFont{5}{6.0}{\rmdefault}{\mddefault}{\updefault}1}}}}}
\put(1058,1133){\ellipse{150}{150}}
\put(1034,1110){\makebox(0,0)[lb]{\smash{{{\SetFigFont{5}{6.0}{\rmdefault}{\mddefault}{\updefault}1}}}}}
\put(1733,683){\ellipse{150}{150}}
\put(1709,660){\makebox(0,0)[lb]{\smash{{{\SetFigFont{5}{6.0}{\rmdefault}{\mddefault}{\updefault}1}}}}}
\put(2633,1133){\ellipse{150}{150}}
\put(2609,1110){\makebox(0,0)[lb]{\smash{{{\SetFigFont{5}{6.0}{\rmdefault}{\mddefault}{\updefault}1}}}}}
\put(3383,1133){\ellipse{150}{150}}
\put(3359,1110){\makebox(0,0)[lb]{\smash{{{\SetFigFont{5}{6.0}{\rmdefault}{\mddefault}{\updefault}1}}}}}
\put(2333,683){\ellipse{150}{150}}
\put(2309,660){\makebox(0,0)[lb]{\smash{{{\SetFigFont{5}{6.0}{\rmdefault}{\mddefault}{\updefault}1}}}}}
\put(2783,83){\ellipse{150}{150}}
\put(2759,60){\makebox(0,0)[lb]{\smash{{{\SetFigFont{5}{6.0}{\rmdefault}{\mddefault}{\updefault}1}}}}}
\put(3233,83){\ellipse{150}{150}}
\put(3209,60){\makebox(0,0)[lb]{\smash{{{\SetFigFont{5}{6.0}{\rmdefault}{\mddefault}{\updefault}1}}}}}
\put(3683,683){\ellipse{150}{150}}
\put(3659,660){\makebox(0,0)[lb]{\smash{{{\SetFigFont{5}{6.0}{\rmdefault}{\mddefault}{\updefault}1}}}}}
\put(458,758){\blacken\ellipse{74}{74}}
\put(458,758){\ellipse{74}{74}}
\put(758,608){\blacken\ellipse{74}{74}}
\put(758,608){\ellipse{74}{74}}
\put(1058,758){\blacken\ellipse{74}{74}}
\put(1058,758){\ellipse{74}{74}}
\put(1358,608){\blacken\ellipse{74}{74}}
\put(1358,608){\ellipse{74}{74}}
\put(3008,608){\blacken\ellipse{74}{74}}
\put(3008,608){\ellipse{74}{74}}
\put(2708,458){\blacken\ellipse{74}{74}}
\put(2708,458){\ellipse{74}{74}}
\put(3308,458){\blacken\ellipse{74}{74}}
\put(3308,458){\ellipse{74}{74}}
\put(3008,908){\blacken\ellipse{74}{74}}
\put(3008,908){\ellipse{74}{74}}
\path(2708,458)(3008,608)
\path(3008,608)(3008,908)
\path(3008,608)(3308,458)
\path(458,758)(758,608)(1058,758)(1358,608)
\path(1356,603)(1472,368)
\path(1361,601)(1657,664)
\path(761,612)(578,293)
\path(462,753)(395,1054)
\path(460,753)(153,689)
\path(1060,753)(1060,1057)
\path(2693,1085)(3008,902)
\path(3011,902)(3321,1088)
\path(3316,462)(3623,645)
\path(3313,454)(3255,155)
\path(2707,451)(2765,155)
\path(2707,451)(2397,645)
\end{picture}
}
\end{center}
\caption{The two $0$-pointed genus $6$ stable trees with at least $4$
genus $0$ vertices\lremind{g6}}\label{g6}
\end{figure}
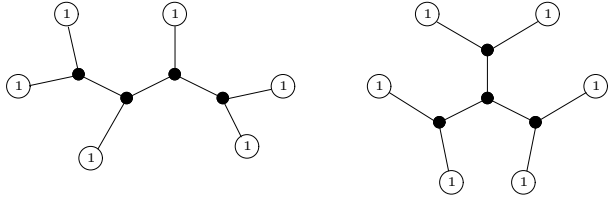

Next, if $i=2g-3$, then our codimension $2g-3$ (hence $g$) class is
pushed forward from strata of the form described in the previous
exercise.  But each stratum has dimension $g$,
so the tautological class must be a linear combination
of fundamental classes of such strata.

Furthermore, any two such strata are equivalent (in cohomology, or even in
the Chow ring) by arguments analogous to either of those we used for
$\cmbar_{g,n}$.

Thus we have shown $R^{2g-3}(\cm_g^c)$ is generated by the fundamental
class of a single such stratum.  It remains to show that this is
non-zero.  This argument is short, but requires a little more
background than we have presented.  (For the experts: it suffices to
show that $\la_g \neq 0$ on this stratum $\cm_{\Ga}$.  We have a cover
$\pi \cmbar_{1,1}^g \rightarrow \cm_{\Ga}$ via gluing morphisms, and
the pullback of the Hodge bundle splits into the Hodge bundles of each
of the $g$ elliptic curves. Thus $\pi^* \la_g$ is the product of the
$\la_1$-classes on each factor, so $\deg \pi^* \la_g =
(\int_{\cmbar_{1,1}} \la_1)^{g} = 1/24^g \neq 0$.)

As always, Theorem $\star$ gives extra information.  (1) This argument
extends to curves of compact type with points.  (2) We can now attack
part of the Poincar\'e duality portion of the conjecture.  (3) We get
an explicit generator of $R^{2g-3}(\cm_g^c)$ (a stratum of a particular form,
e.g.\ Figure~\ref{g6}).

\epoint{Theorem $\star$ helps determine the tautological ring
in low dimension}
In the course of proving $R_0(\cmbar_{g,n}) \cong \Q$, we showed that
$R_0(\cmbar_{g,n})$ was generated by 0-strata.  A similar argument
shows that $R_i(\cmbar_{g,n})$ generated by boundary strata for
$i=1,2$.  (We are already aware that this will not extend to $i=3$, as
$\psi_1$ on $\cmbar_{2,1}$ is the fundamental class of a stratum.)

In general, Theorem $\star$ implies that in order to understand
tautological classes in dimension up to $i$, you need only understand
curves of genus up to $(i+1)/2$, with not too many marked points.

The moral of this is that the ``top'' (lowest-codimension) part of the
tautological ring used to be considered the least mysterious (given
the definition of the tautological ring, it is easy to give
generators), and the bottom was therefore the most mysterious.  Now
the situation is the opposite. For example, in codimension 3, we can
describe the generators of the tautological ring, but we have no idea
what the relations are.  However, we know exactly what the
tautological ring looks like in {\em dimension} $3$.

\epoint{Exercise} Use Theorem $\star$ and a similar argument to show
that the tautological groups of $\cmbar_{0,n}$ and $\cmbar_{1,n}$ are
generated by boundary strata.
\label{g0g1}\lremind{g0g1}

\epoint{Additional consequences} For many additional consequences of
Theorem $\star$, see \cite{thmstar}.  For example, we recover {\em
  Diaz' theorem} ($\cm_g$ contains no complete subvarieties of $\dim >
g-2$), as well as generalizations and variations such as:
$\cmbar_{g,n}^{c}$ contains no complete subvarieties of $\dim >
2g-3+n$.

The idea behind the proof of Theorem $\star$ is rather naive.
But before we can discuss it, we will have to finally enter
the land of Gromov-Witten theory, and define stable relative maps
to $\proj^1$, which we will interpret as a generalization of the notion of 
a branched cover.

\section{Stable relative maps to $\proj^1$ and relative virtual localization}

We\label{s:rvl}\lremind{s:rvl}
 now discuss the theory of {\em stable relative maps}, and
``virtual'' localization on their moduli space ({\em relative virtual
  localization}).  We will follow J. Li's algebro-geometric definition
of stable relative maps \cite{li1}, and his description of their
obstruction theory \cite{li2}, but we point out earlier definitions of
stable relative maps in the differentiable category due to A.-M. Li
and Y. Ruan \cite{lr}, and Ionel and Parker \cite{ip1, ip2}, and
Gathmann's work in the algebraic category in genus 0 \cite{gathmann}.
We need the algebraic category for several reasons, most importantly
because we will want to apply virtual localization.

Stable relative maps are variations of the notion of stable maps, and
the reader may wish to become comfortable with that notion first.
(Stable maps are discussed in Abramovich's article in this volume, for
example.)

We are interested in the particular case of stable relative maps to
$\proj^1$, relative to at most two points, so we will define stable relative
maps only in this case.   
For concreteness, we define stable maps to $X=\proj^1$ relative to
one point $\infty$; the case of zero or two points is the obvious
variation on this theme.  Such a stable relative map to $(\proj^1,
\infty)$ is defined as follows.  We are given the data of a degree $d$
of the map, a genus $g$ of the source curve, a number $m$ of marked
points, and a partition $d = \al_1 + \cdots + \al_n$, which we write
$\al \vdash d$.

Then a relative map is the following data:
\begin{itemize}
\item
a morphism $f_1$ from a nodal 
$(m+n)$-pointed genus $g$ curve $(C, p_1, \dots, p_m, q_1, \dots, q_n)$ (where
as usual the $p_i$ and $q_j$ are distinct nonsingular points) to a
chain of $\proj^1$'s, $T = T_0 \cup T_1 \cup \cdots \cup T_t$ (where
$T_i$ and $T_{i+1}$ meet), with a point $\infty \in T_t - T_{t-1}$.
Unfortunately, there are two points named $\infty$.  We will call the one on $X$,  $\infty_X$, and the one on $T$, $\infty_T$, whenever there is any ambiguity.
\item  A projection $f_2: T \rightarrow X$ contracting $T_i$ 
to $\infty_X$ (for $i>0$)
and giving an isomorphism from $(T_0, T_0 \cap T_1)$ 
(resp.\ $(T_0, \infty)$) 
to $X$ if $t>0$  (resp.\ if $t=0$).  Denote $f_2 \circ f_1$ by $f$.
\item We have an equality of divisors on $C$: $f_1^* \infty_T = \sum \al_i q_i$.
In particular,  $f_1^{-1} \infty_T$ consists 
of nonsingular (marked) points of $C$.
\item The preimage of each node $n$ of $T$ is a union of nodes of $C$.
  At any such node $n'$ of $C$, the two branches map to the two
  branches of $n$, and their orders of branching are the same.  (This
  is called the {\em predeformability} or {\em kissing} condition.)
\end{itemize}
If follows that the degree of  $f_1$ is
$d$ on each $T_i$.
An {\em isomorphism} of two such maps is a commuting diagram
$$
\xymatrix{
(C, p_1, \dots, p_m, q_1, \dots, q_n) \ar[d]_{f_1} \ar[r]^{\sim} & 
(C', p'_1, \dots, p'_m, q'_1, \dots, q'_n) \ar[d]^{f_1} \\
(T, \infty_T) \ar[r]^{\sim} \ar[d]_{f_2} &  (T, \infty_T) \ar[d]^{f_2}\\
(X,\infty_X) \ar[r]^= & (X,\infty_X)}
$$
where all horizontal morphisms are isomorphisms, the bottom (although
not necessarily the middle!) is an equality, the top horizontal
isomorphism sends $p_i$ to $p'_i$ and $q_j$ to $q'_j$.  Note that the
middle isomorphism must  preserve the isomorphism of $T_0$ with $X$, and is
hence the identity on $T_0$, but for $i>0$, the isomorphism may not be
the identity on $T_i$.

This data of a relative map is often just denoted $f$, with the remaining information left
 implicit.

We say that $f$ is {\em stable} if it has finite automorphism group.
This corresponds to the following criteria. 
\begin{itemize}
\item Any $f_1$-contracted geometric genus $0$ component has at least
  $3$ ``special points'' (node branches or marked points).
\item Any $f_1$-contracted geometric genus $1$ component has at least $1$ ``special point''.
\item If $0<i<t$ (resp.\ $0<i=t$), then the preimage of $T_i - T_{i+1} - T_{i-1}$
(resp.\ $T_i - \{ \infty \} - T_{i-1}$) is not smooth  unmarked
curve.  In other words, not every component mapping to $T_i$ is of the
form $[x;y] \rightarrow [x^q;y^q]$, where the coordinates on the target
are given by  $[0;1] = T_i \cap T_{i-1}$
and $[1;0] = T_i \cap T_{i+1}$  (resp.\ $[1;0] = \infty$).
\end{itemize}
 (The first two conditions
are the same as for stable maps.  The third condition is new.)
A picture of a stable relative map is given in Figure~\ref{rsmfig}.
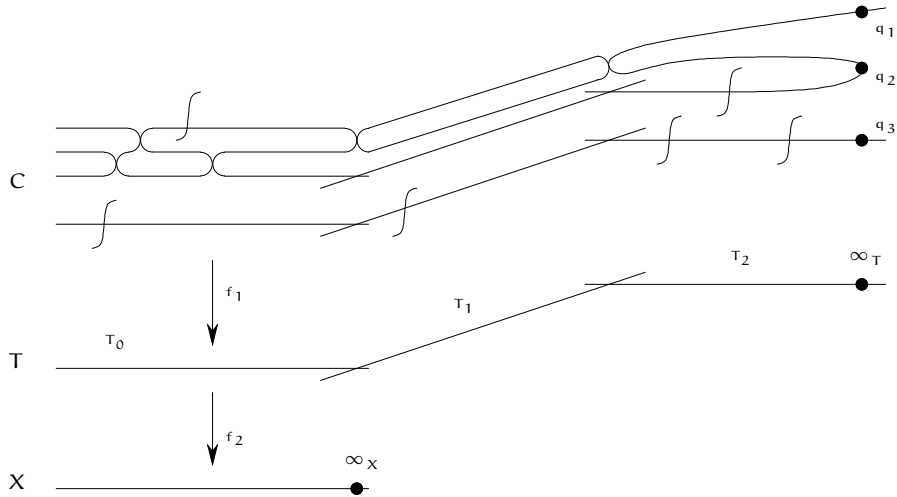
\begin{figure}
\begin{center}
\setlength{\unitlength}{0.00083333in}
\begingroup\makeatletter\ifx\SetFigFont\undefined%
\gdef\SetFigFont#1#2#3#4#5{%
  \reset@font\fontsize{#1}{#2pt}%
  \fontfamily{#3}\fontseries{#4}\fontshape{#5}%
  \selectfont}%
\fi\endgroup%
{\renewcommand{\dashlinestretch}{30}
\begin{picture}(5487,3080)(0,-10)
\put(1350.000,2070.000){\arc{150.000}{1.5708}{4.7124}}
\put(1200.000,2070.000){\arc{150.000}{4.7124}{7.8540}}
\put(2100.000,2220.000){\arc{150.000}{4.7124}{7.8540}}
\put(2250.000,2220.000){\arc{150.000}{1.5708}{4.7124}}
\put(900.000,2220.000){\arc{150.000}{1.5708}{4.7124}}
\put(600.000,2070.000){\arc{150.000}{4.7124}{7.8540}}
\put(750.000,2220.000){\arc{150.000}{4.7124}{7.8540}}
\put(750.000,2070.000){\arc{150.000}{1.5708}{4.7124}}
\put(3675.000,2670.000){\arc{150.000}{4.7124}{7.8540}}
\put(2175,45){\blacken\ellipse{74}{74}}
\put(2175,45){\ellipse{74}{74}}
\put(5325,1320){\blacken\ellipse{74}{74}}
\put(5325,1320){\ellipse{74}{74}}
\put(5325,2220){\blacken\ellipse{74}{74}}
\put(5325,2220){\ellipse{74}{74}}
\put(5325,3020){\blacken\ellipse{74}{74}}
\put(5325,3020){\ellipse{74}{74}}
\put(5325,2670){\blacken\ellipse{74}{74}}
\put(5325,2670){\ellipse{74}{74}}
\path(1275,1470)(1275,945)
\blacken\path(1245.000,1065.000)(1275.000,945.000)(1305.000,1065.000)(1275.000,1029.000)(1245.000,1065.000)
\path(1275,645)(1275,195)
\blacken\path(1245.000,315.000)(1275.000,195.000)(1305.000,315.000)(1275.000,279.000)(1245.000,315.000)
\path(300,795)(2250,795)
\path(300,45)(2250,45)
\path(300,1695)(2250,1695)
\path(2250,1995)(1350,1995)
\path(1350,2145)(2100,2145)
\path(300,1995)(600,1995)
\path(600,2145)(300,2145)
\path(300,2295)(750,2295)
\path(750,1995)(1200,1995)
\path(900,2145)(1200,2145)
\path(900,2295)(2100,2295)
\path(1950,720)(3975,1395)
\path(1950,1620)(3975,2295)
\path(1950,1920)(3975,2595)
\path(2250,2145)(3675,2595)
\path(2250,2295)(3675,2745)
\path(3600,1320)(5475,1320)
\path(3600,2220)(5475,2220)
\path(2550,1920)(2549,1920)(2546,1919)
	(2539,1918)(2530,1915)(2519,1911)
	(2510,1906)(2502,1900)(2496,1892)
	(2491,1882)(2488,1870)(2485,1859)
	(2483,1848)(2481,1834)(2480,1819)
	(2478,1804)(2476,1787)(2475,1770)
	(2474,1753)(2472,1736)(2470,1721)
	(2469,1706)(2467,1692)(2465,1681)
	(2463,1670)(2459,1658)(2454,1648)
	(2448,1640)(2440,1634)(2431,1629)
	(2420,1625)(2411,1622)(2404,1621)
	(2401,1620)(2400,1620)
\path(675,1845)(674,1845)(671,1844)
	(664,1843)(655,1840)(644,1836)
	(635,1831)(627,1825)(621,1817)
	(616,1807)(613,1795)(610,1784)
	(608,1773)(606,1759)(605,1744)
	(603,1729)(601,1712)(600,1695)
	(599,1678)(597,1661)(595,1646)
	(594,1631)(592,1617)(590,1606)
	(588,1595)(584,1583)(579,1573)
	(573,1565)(565,1559)(556,1554)
	(545,1550)(536,1547)(529,1546)
	(526,1545)(525,1545)
\path(1200,2520)(1199,2520)(1196,2519)
	(1189,2518)(1180,2515)(1169,2511)
	(1160,2506)(1152,2500)(1146,2492)
	(1141,2482)(1138,2470)(1135,2459)
	(1133,2448)(1131,2434)(1130,2419)
	(1128,2404)(1126,2387)(1125,2370)
	(1124,2353)(1122,2336)(1120,2321)
	(1119,2306)(1117,2292)(1115,2281)
	(1113,2270)(1109,2258)(1104,2248)
	(1098,2240)(1090,2234)(1081,2229)
	(1070,2225)(1061,2222)(1054,2221)
	(1051,2220)(1050,2220)
\path(3600,2520)(3603,2520)(3610,2520)
	(3622,2520)(3639,2520)(3659,2520)
	(3682,2520)(3705,2520)(3727,2520)
	(3749,2520)(3769,2520)(3788,2520)
	(3807,2520)(3825,2520)(3844,2520)
	(3863,2520)(3877,2520)(3892,2520)
	(3907,2520)(3924,2520)(3941,2520)
	(3959,2520)(3979,2520)(3999,2520)
	(4020,2520)(4043,2520)(4066,2520)
	(4089,2520)(4114,2520)(4139,2520)
	(4164,2520)(4190,2520)(4216,2520)
	(4242,2520)(4268,2520)(4295,2520)
	(4322,2520)(4350,2520)(4374,2520)
	(4399,2520)(4425,2520)(4451,2520)
	(4478,2520)(4506,2520)(4535,2521)
	(4565,2521)(4595,2521)(4626,2522)
	(4657,2522)(4688,2523)(4720,2523)
	(4751,2524)(4782,2525)(4812,2526)
	(4842,2528)(4870,2529)(4898,2531)
	(4925,2532)(4951,2534)(4975,2536)
	(4999,2538)(5021,2540)(5042,2542)
	(5063,2545)(5094,2550)(5123,2555)
	(5150,2560)(5177,2567)(5201,2574)
	(5224,2581)(5244,2589)(5263,2597)
	(5280,2605)(5294,2613)(5306,2621)
	(5316,2629)(5324,2637)(5330,2644)
	(5335,2651)(5338,2657)(5339,2663)
	(5339,2668)(5338,2672)(5337,2677)
	(5333,2682)(5329,2687)(5323,2691)
	(5315,2696)(5305,2700)(5294,2704)
	(5281,2708)(5266,2712)(5249,2715)
	(5231,2718)(5211,2721)(5190,2724)
	(5167,2726)(5142,2729)(5116,2731)
	(5088,2732)(5066,2734)(5043,2735)
	(5020,2736)(4994,2737)(4968,2738)
	(4940,2738)(4912,2739)(4882,2739)
	(4851,2740)(4819,2740)(4786,2740)
	(4752,2740)(4718,2740)(4684,2740)
	(4650,2739)(4615,2739)(4581,2738)
	(4548,2737)(4515,2736)(4483,2735)
	(4452,2733)(4421,2732)(4392,2730)
	(4364,2729)(4337,2727)(4311,2725)
	(4286,2722)(4263,2720)(4231,2716)
	(4200,2712)(4171,2708)(4143,2704)
	(4116,2699)(4091,2694)(4066,2688)
	(4042,2683)(4020,2678)(3998,2672)
	(3978,2667)(3960,2662)(3943,2657)
	(3927,2653)(3912,2649)(3899,2645)
	(3887,2642)(3876,2639)(3865,2637)
	(3855,2635)(3844,2634)(3832,2634)
	(3822,2634)(3811,2636)(3801,2638)
	(3792,2640)(3784,2644)(3776,2648)
	(3769,2653)(3763,2658)(3757,2664)
	(3753,2670)(3750,2676)(3748,2682)
	(3747,2688)(3746,2694)(3747,2700)
	(3748,2705)(3750,2711)(3753,2716)
	(3757,2723)(3763,2729)(3771,2736)
	(3781,2742)(3793,2750)(3807,2757)
	(3823,2765)(3842,2772)(3863,2780)
	(3887,2788)(3913,2796)(3941,2804)
	(3972,2812)(4006,2819)(4041,2827)
	(4080,2835)(4103,2840)(4128,2845)
	(4153,2850)(4180,2854)(4209,2859)
	(4240,2865)(4272,2870)(4307,2876)
	(4345,2882)(4384,2888)(4427,2895)
	(4472,2901)(4520,2909)(4571,2916)
	(4624,2924)(4681,2932)(4739,2941)
	(4800,2950)(4862,2959)(4926,2968)
	(4989,2977)(5052,2986)(5114,2994)
	(5173,3003)(5228,3010)(5279,3018)
	(5325,3024)(5364,3030)(5397,3034)
	(5424,3038)(5444,3041)(5459,3043)
	(5468,3044)(5473,3045)(5475,3045)
\path(4200,2370)(4199,2370)(4196,2369)
	(4189,2368)(4180,2365)(4169,2361)
	(4160,2356)(4152,2350)(4146,2342)
	(4141,2332)(4138,2320)(4135,2309)
	(4133,2298)(4131,2284)(4130,2269)
	(4128,2254)(4126,2237)(4125,2220)
	(4124,2203)(4122,2186)(4120,2171)
	(4119,2156)(4117,2142)(4115,2131)
	(4113,2120)(4109,2108)(4104,2098)
	(4098,2090)(4090,2084)(4081,2079)
	(4070,2075)(4061,2072)(4054,2071)
	(4051,2070)(4050,2070)
\path(4950,2370)(4949,2370)(4946,2369)
	(4939,2368)(4930,2365)(4919,2361)
	(4910,2356)(4902,2350)(4896,2342)
	(4891,2332)(4888,2320)(4885,2309)
	(4883,2298)(4881,2284)(4880,2269)
	(4878,2254)(4876,2237)(4875,2220)
	(4874,2203)(4872,2186)(4870,2171)
	(4869,2156)(4867,2142)(4865,2131)
	(4863,2120)(4859,2108)(4854,2098)
	(4848,2090)(4840,2084)(4831,2079)
	(4820,2075)(4811,2072)(4804,2071)
	(4801,2070)(4800,2070)
\path(4575,2670)(4574,2670)(4571,2669)
	(4564,2668)(4555,2665)(4544,2661)
	(4535,2656)(4527,2650)(4521,2642)
	(4516,2632)(4513,2620)(4510,2609)
	(4508,2598)(4506,2584)(4505,2569)
	(4503,2554)(4501,2537)(4500,2520)
	(4499,2503)(4497,2486)(4495,2471)
	(4494,2456)(4492,2442)(4490,2431)
	(4488,2420)(4484,2408)(4479,2398)
	(4473,2390)(4465,2384)(4456,2379)
	(4445,2375)(4436,2372)(4429,2371)
	(4426,2370)(4425,2370)
\put(5250,1470){\makebox(0,0)[lb]{\smash{{{\SetFigFont{5}{6.0}{\rmdefault}{\mddefault}{\updefault}$\infty_T$}}}}}
\put(600,945){\makebox(0,0)[lb]{\smash{{{\SetFigFont{5}{6.0}{\rmdefault}{\mddefault}{\updefault}$T_0$}}}}}
\put(0,45){\makebox(0,0)[lb]{\smash{{{\SetFigFont{8}{9.6}{\rmdefault}{\mddefault}{\updefault}$X$}}}}}
\put(0,795){\makebox(0,0)[lb]{\smash{{{\SetFigFont{8}{9.6}{\rmdefault}{\mddefault}{\updefault}$T$}}}}}
\put(0,1920){\makebox(0,0)[lb]{\smash{{{\SetFigFont{8}{9.6}{\rmdefault}{\mddefault}{\updefault}$C$}}}}}
\put(5400,2595){\makebox(0,0)[lb]{\smash{{{\SetFigFont{5}{6.0}{\rmdefault}{\mddefault}{\updefault}$q_2$}}}}}
\put(5400,2295){\makebox(0,0)[lb]{\smash{{{\SetFigFont{5}{6.0}{\rmdefault}{\mddefault}{\updefault}$q_3$}}}}}
\put(1350,1245){\makebox(0,0)[lb]{\smash{{{\SetFigFont{5}{6.0}{\rmdefault}{\mddefault}{\updefault}$f_1$}}}}}
\put(1350,345){\makebox(0,0)[lb]{\smash{{{\SetFigFont{5}{6.0}{\rmdefault}{\mddefault}{\updefault}$f_2$}}}}}
\put(5400,2895){\makebox(0,0)[lb]{\smash{{{\SetFigFont{5}{6.0}{\rmdefault}{\mddefault}{\updefault}$q_1$}}}}}
\put(2100,195){\makebox(0,0)[lb]{\smash{{{\SetFigFont{5}{6.0}{\rmdefault}{\mddefault}{\updefault}$\infty_X$}}}}}
\put(4500,1470){\makebox(0,0)[lb]{\smash{{{\SetFigFont{5}{6.0}{\rmdefault}{\mddefault}{\updefault}$T_2$}}}}}
\put(2775,1170){\makebox(0,0)[lb]{\smash{{{\SetFigFont{5}{6.0}{\rmdefault}{\mddefault}{\updefault}$T_1$}}}}}
\end{picture}
}
\end{center}
\caption{An example of a stable relative map \lremind{rsmfig}
}\label{rsmfig}
\end{figure}

Thus we have some behavior familiar from the theory 
of stable maps: we can have contracted components, so long as they are
``stable'', and don't map to any nodes of $T$, or to $\infty_T$.  We
also have some new behavior: the target $X$ can ``sprout'' a chain of
$\proj^1$'s at $\infty_X$.  Also, the action of $\C^*$ on the map via
the action on a component $T_i$ ($i>0$) that preserves the two
``special points'' of $T_i$ (the intersections with $T_{i-1}$ and
$T_{i+1}$ if $i<t$, and the intersection with $T_{i-1}$ and $\infty$
if $i=t$) is considered to preserve the stable map.  For example,
Figure~\ref{isorsm} shows two isomorphic stable maps.
\begin{figure}
\begin{center}
\setlength{\unitlength}{0.00083333in}
\begingroup\makeatletter\ifx\SetFigFont\undefined%
\gdef\SetFigFont#1#2#3#4#5{%
  \reset@font\fontsize{#1}{#2pt}%
  \fontfamily{#3}\fontseries{#4}\fontshape{#5}%
  \selectfont}%
\fi\endgroup%
{\renewcommand{\dashlinestretch}{30}
\begin{picture}(7824,2772)(0,-10)
\put(1062.000,2070.000){\arc{150.000}{1.5708}{4.7124}}
\put(912.000,2070.000){\arc{150.000}{4.7124}{7.8540}}
\put(1812.000,2220.000){\arc{150.000}{4.7124}{7.8540}}
\put(1962.000,2220.000){\arc{150.000}{1.5708}{4.7124}}
\put(612.000,2220.000){\arc{150.000}{1.5708}{4.7124}}
\put(312.000,2070.000){\arc{150.000}{4.7124}{7.8540}}
\put(462.000,2220.000){\arc{150.000}{4.7124}{7.8540}}
\put(462.000,2070.000){\arc{150.000}{1.5708}{4.7124}}
\put(3387.000,2670.000){\arc{150.000}{4.7124}{7.8540}}
\put(5187.000,2070.000){\arc{150.000}{1.5708}{4.7124}}
\put(5037.000,2070.000){\arc{150.000}{4.7124}{7.8540}}
\put(5937.000,2220.000){\arc{150.000}{4.7124}{7.8540}}
\put(6087.000,2220.000){\arc{150.000}{1.5708}{4.7124}}
\put(4737.000,2220.000){\arc{150.000}{1.5708}{4.7124}}
\put(4437.000,2070.000){\arc{150.000}{4.7124}{7.8540}}
\put(4587.000,2220.000){\arc{150.000}{4.7124}{7.8540}}
\put(4587.000,2070.000){\arc{150.000}{1.5708}{4.7124}}
\put(7512.000,2670.000){\arc{150.000}{4.7124}{7.8540}}
\put(1887,45){\blacken\ellipse{74}{74}}
\put(1887,45){\ellipse{74}{74}}
\put(3462,1320){\blacken\ellipse{74}{74}}
\put(3462,1320){\ellipse{74}{74}}
\put(6012,45){\blacken\ellipse{74}{74}}
\put(6012,45){\ellipse{74}{74}}
\put(7587,1320){\blacken\ellipse{74}{74}}
\put(7587,1320){\ellipse{74}{74}}
\path(987,1470)(987,945)
\blacken\path(957.000,1065.000)(987.000,945.000)(1017.000,1065.000)(987.000,1029.000)(957.000,1065.000)
\path(987,645)(987,195)
\blacken\path(957.000,315.000)(987.000,195.000)(1017.000,315.000)(987.000,279.000)(957.000,315.000)
\path(12,795)(1962,795)
\path(12,45)(1962,45)
\path(12,1695)(1962,1695)
\path(1962,1995)(1062,1995)
\path(1062,2145)(1812,2145)
\path(12,1995)(312,1995)
\path(312,2145)(12,2145)
\path(12,2295)(462,2295)
\path(462,1995)(912,1995)
\path(612,2145)(912,2145)
\path(612,2295)(1812,2295)
\path(1662,720)(3687,1395)
\path(1662,1620)(3687,2295)
\path(1662,1920)(3687,2595)
\path(1962,2145)(3387,2595)
\path(1962,2295)(3387,2745)
\path(5112,1470)(5112,945)
\blacken\path(5082.000,1065.000)(5112.000,945.000)(5142.000,1065.000)(5112.000,1029.000)(5082.000,1065.000)
\path(5112,645)(5112,195)
\blacken\path(5082.000,315.000)(5112.000,195.000)(5142.000,315.000)(5112.000,279.000)(5082.000,315.000)
\path(4137,795)(6087,795)
\path(4137,45)(6087,45)
\path(4137,1695)(6087,1695)
\path(6087,1995)(5187,1995)
\path(5187,2145)(5937,2145)
\path(4137,1995)(4437,1995)
\path(4437,2145)(4137,2145)
\path(4137,2295)(4587,2295)
\path(4587,1995)(5037,1995)
\path(4737,2145)(5037,2145)
\path(4737,2295)(5937,2295)
\path(5787,720)(7812,1395)
\path(5787,1620)(7812,2295)
\path(5787,1920)(7812,2595)
\path(6087,2145)(7512,2595)
\path(6087,2295)(7512,2745)
\path(2262,1920)(2261,1920)(2258,1919)
	(2251,1918)(2242,1915)(2231,1911)
	(2222,1906)(2214,1900)(2208,1892)
	(2203,1882)(2200,1870)(2197,1859)
	(2195,1848)(2193,1834)(2192,1819)
	(2190,1804)(2188,1787)(2187,1770)
	(2186,1753)(2184,1736)(2182,1721)
	(2181,1706)(2179,1692)(2177,1681)
	(2175,1670)(2171,1658)(2166,1648)
	(2160,1640)(2152,1634)(2143,1629)
	(2132,1625)(2123,1622)(2116,1621)
	(2113,1620)(2112,1620)
\path(7437,2295)(7436,2295)(7433,2294)
	(7426,2293)(7417,2290)(7406,2286)
	(7397,2281)(7389,2275)(7383,2267)
	(7378,2257)(7375,2245)(7372,2234)
	(7370,2223)(7368,2209)(7367,2194)
	(7365,2179)(7363,2162)(7362,2145)
	(7361,2128)(7359,2111)(7357,2096)
	(7356,2081)(7354,2067)(7352,2056)
	(7350,2045)(7346,2033)(7341,2023)
	(7335,2015)(7327,2009)(7318,2004)
	(7307,2000)(7298,1997)(7291,1996)
	(7288,1995)(7287,1995)
\end{picture}
}
\end{center}
\caption{Two isomorphic stable relative maps \lremind{isorsm}}\label{isorsm}
\end{figure}
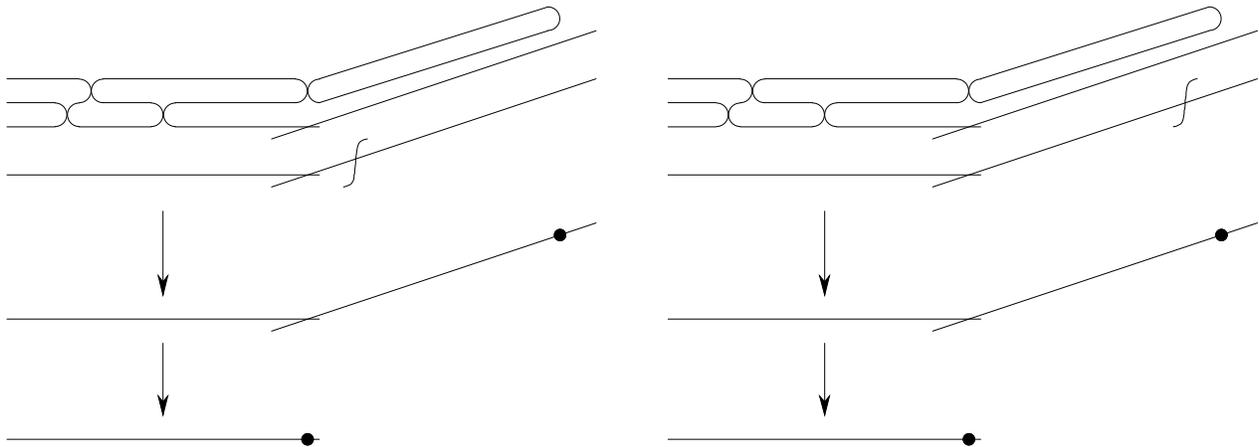

There is a compact moduli space (Deligne-Mumford stack) for stable relative
maps to $\proj^1$, denoted $\cmbar_{g,m,\al}(\proj^1, d)$.  (In order
to be more precise, I should tell you the definition of a family of
stable relative maps parametrized by an arbitrary base, but I will not
do so.)  In what follows, $m=0$, and that subscript will be omitted.
(More generally, stable relative maps may be defined with $\proj^1$ replaced by any
smooth complex projective variety, and $d$ replaced by any smooth
divisor $D$ on $X$.  The special case $D=\emptyset$ yields Kontsevich's
original space of stable maps.)

Unfortunately, the space $\cmbar_{g,\al}(\proj^1, d)$ is in general
terribly singular, and not even equidimensional.  

\noindent {\em Exercise.}  Give an example of such a moduli
space with two components of different dimensions.  (Hint:
use contracted components judiciously.)

However, it has a component which we already understand well, which
corresponds to maps from a smooth curve, which is a branched cover of
$\proj^1$.  Such curves form a moduli space $\cm_{g,\al}(\proj^1,d)$
of dimension corresponding to the ``expected number of branch points
distinct from $\infty$'', which we may calculate by the
Riemann-Hurwitz formula~\eqref{rh} to be \lremind{rh2}
\begin{equation}
\label{rh2}
r = 2g-2+n+d.\end{equation}
We have seen this formula before, \eqref{Hr}.

\noindent
{\em Exercise.}  Verify \eqref{rh2}.

These notions can be readily generalized, for example
to stable relative maps to $\proj^1$ relative to two points
(whose moduli space is denoted $\cmbar_{g,\al,\be}(\proj^1,d)$),
or to no points (otherwise known as the stable maps to $\proj^1$;
this moduli space is denoted $\cmbar_{g}(\proj^1,d)$).

\noindent {\em Exercise.}
Calculate $\dim \cm_{g,\al,\be}(\proj^1,d)$ (where $\al$ has $m$ parts
and $\be$ has $n$ parts) and $\dim \cm_{g}(\proj^1,d)$.

\epoint{Stable relative maps with possibly-disconnected source curve}
Recall that by our (non-standard) definition, nodal curves are connected.
It will be convenient, especially when discussing the degeneration
formula, to consider curves without this hypothesis.
Just as our discussion of (connected) stable curves generalized without
change to (possibly-disconnected) stable curves (see \S \ref{sic}),
our discussion of (relatively) stable maps from connected curves generalizes
without change to ``(relatively) stable maps from possibly-disconnected curves''.
Let $\cmbar_{g,\al}(\proj^1,d)^{\dis}$ be the space of stable relative
maps from possibly-disconnected curves (to $\proj^1$, of degree $d$,
etc.).  
Warning:  this is {\em not} in general the quotient of a product
of $\cmbar_{g', \al'}(\proj^1, d')$'s by a finite group.

\bpoint{The virtual fundamental class}

There is a natural homology (or Chow) class on
$\cmbar_{g,\al}(\proj^1, d)$ of dimension $r= \dim
\cm_{g,\al}(\proj^1,d)$ (cf.\ \eqref{rh2}), called the {\em virtual
  fundamental class}
$[ \cmbar_{g,\al}(\proj^1,d)]^{\virt} \in A_r(\cmbar_{g,\al}(\proj^1,d)])$, 
which is obtained from the
deformation-obstruction theory of stable relative maps, and has many
wonderful properties.  The virtual fundamental class agrees with the
actual fundamental class on the open subset  $\cm_{g,\al}(\proj^1,d)$.  The most
difficult part of dealing with the moduli space of stable relative
maps is working with the virtual fundamental class.

\noindent
{\em Aside:  relative Gromov-Witten invariants.}
In analogy with usual Gromov-Witten invariants, one can define {\em
  relative Gromov-Witten invariants} by intersecting natural
cohomology classes on the moduli space with the virtual fundamental
class.  More precisely, one multiplies (via the cup/cap product) the
cohomology classes with the virtual fundamental class, and takes the
degree of the resulting zero-cycle.  One can define $\psi$-classes and
$\la$-classes in the same way as before, and include these in the
product.  When including $\psi$-classes, the numbers are often called
{\em descendant relative invariants}; when including $\la$-classes,
the numbers are sometimes called {\em Hodge integrals} for some
reason.  For example, one can show that 
Hurwitz numbers are descendant relative
invariants of $\proj^1$.  However, this point of view turns out to be
less helpful, and we will not use the language of relative
Gromov-Witten invariants again.

The virtual fundamental class behaves well under two procedures:
degeneration and localization; we now discuss these.

\bpoint{The degeneration formula for the virtual fundamental class,
  following \cite{li2}}

We describe the degeneration formula in the case of stable maps to
$\proj^1$ relative to one point, and leave the cases of stable maps to
$\proj^1$ relative to zero or two points as straightforward
variations for the reader.  In this discussion, we will deal with
possibly-disconnected curves to simplify the exposition.

Consider the maps to $\proj^1$ relative to one point $\infty$, and
imagine deforming the target so that it breaks into two $\proj^1$'s,
meeting at a node (with $\infty$ on one of the components).  It turns
out that the virtual fundamental class behaves well under this
degeneration.  The limit can be expressed in terms of virtual
fundamental classes of spaces of stable relative maps to each
component, relative to $\infty$ (for the component containing
$\infty$), and relative to the node-branch (for both components).

Before we make this precise, we give some intuition.  Suppose we have
a branched cover $C \rightarrow \proj^1$, and we deform the target
into a union of two $\proj^1$'s, while keeping the branch points away
from the node; call the limit map $C' \rightarrow \proj^1 \cup
\proj^1$.  Clearly in the limit, away from the node, the cover looks
just the same as it did before (with the same branching).  At the
node, it turns out that the branched covers of the two components must
satisfy the kissing/predeformability condition. Say that the
branching above the node corresponds to the partition $\ga_1 + \cdots +
\ga_m$.  By our discussion about Hurwitz numbers, as we have specified
the branch points, there will be a finite number of such branched
covers --- we count branched covers of each component of $\proj^1 \cup
\proj^1$, with branching corresponding to the partition $\ga$ above
the node-branch; then we choose how to match the preimages of the
node-branch on the two components (there are $\# \Aut \ga$ such
choices).  It turns out that $\ga_1 \cdots \ga_m$ covers of the
original sort will degenerate to each branched cover of the nodal
curve of this sort.  (Notice that if we were interested in connected
curves $C$, then the inverse image of each component of $\proj^1$
would not necessarily be connected, and we would have to take some
care in gluing these curves together to get a connected union.  This
is the reason for considering possibly-disconnected components.)

Motivated by the previous paragraph, we give the degeneration formula.
Consider the degeneration of the target $\xymatrix{(\proj^1, \infty)
  \ar@{~>}[r] & \proj^1 \cup (\proj^1 , \infty)}$.  Let $(X, \infty)$ be
the general target, and let $(X', \infty)$ be the degenerated
target. Let $(X_1, a_1) \cong (\proj^1, \infty)$ denote the first
component of $X'$, where $a_1$ refers to the node-branch, and let
$(X_2, a_2, \infty) \cong (\proj^1, 0, \infty)$ denote the second
component of $X'$, where $a_2$ corresponds to the node branch.  Then
for each partition $\ga_1 + \cdots + \ga_m = d$, there is a
natural map\lremind{glue}
\begin{equation}\label{glue}
\cmbar_{g_1, \ga}(\proj^1, d)^{\dis}
\times \cmbar_{g_2, \ga, \al}(\proj^1, d)^{\dis} \rightarrow
\cmbar_{g_1+g_2 - m+1, \al}(X', d)^{\dis}
\end{equation}
obtained by gluing the points above $a_1$ to the corresponding points
above $a_2$. The image of this map can be suitably interpreted as
stable maps to $X'$, satisfying the kissing condition, which can
appear as the limit of maps to $X$.  
(We are obscuring a delicate issue here --- we have not
defined stable maps to a singular target such as $X$.)
Then Li's degeneration formula
states that the image of the product of the virtual fundamental
classes in \eqref{glue} is the limit of the virtual fundamental class
of $\cmbar_{g,\al}(\proj^1,d)^{\dis}$, multiplied by $\ga_1 \cdots
\ga_m$.

The main idea behind Li's proof is remarkably elegant, but as with any argument
involving the virtual fundamental class, the details are quite technical.

If we are interested in connected curves, then there is a
corresponding statement (that requires no additional proof): we look
at the component of the moduli space on the right side of \eqref{glue}
corresponding to maps from connected source curves, and we look at
just those components of the moduli spaces on the left side which glue
together to give connected curves.

\bpoint{Relative virtual localization \cite{thmstar}} 

The second fundamental method of manipulating virtual fundamental
classes is by means of localization.  Before discussing localization
in our Gromov-Witten-theoretic context, we first quickly review
localization in its original setting.

(A friendly introduction to equivariant cohomology is given in
\cite[Ch.~4]{mirsym}, and to localization on the space of ordinary
stable maps in \cite[Ch.~27]{mirsym}.)

Suppose $Y$ is a complex projective manifold with an action by a torus
$\C^*$.  Then the fixed point loci of the torus is the union of
smooth submanifolds, possibly of various dimensions.  Let the
components of the fixed locus be $Y_1$, $Y_2$, \dots.  The torus acts
on the normal bundle $N_i$ to $Y_i$.  Then the Atiyah-Bott localization
formula states that\lremind{localization}
\begin{equation}
\label{localization}
[Y] = \sum_{\fixed} [Y_i] / c_{\top}(N_i) =
\sum_{\fixed} [Y_i] / e(N_i),
\end{equation} in the {\em equivariant homology} of $Y$
(with appropriate terms inverted),
where $c_{\top}$ (or the Euler class $e$) of a vector bundle 
denotes
the top Chern class.  This is a wonderfully powerful fact, and to
appreciate it, you must do examples yourself.  The original paper
of Atiyah and Bott  \cite{ab} is beautifully written and remains a canonical
source.

You can cap \eqref{localization} with various cohomology classes to
get $0$-dimensional classes, and get an equality of numbers.
But you can cap \eqref{localization} with classes to get
higher-dimensional classes, and get equality in cohomology (or the
Chow ring).  One lesson I want to emphasize is that this is a powerful
thing to do.  For example, in a virtual setting, in Gromov-Witten
theory, localization is traditionally used to get equalities of
numbers.  We will also use equalities of numbers to prove the ELSV
formula \eqref{elsv}.  However, using more generally equalities of
classes will give us Theorem $\star$, and part of Faber's conjecture.

Localization was introduced to Gromov-Witten theory by Kontsevich in
his ground-breaking paper \cite{kratcurves}, in which he works on the
space of genus zero maps to projective space, where the virtual
fundamental class is the usual fundamental class (and hence there are
no ``virtual'' technicalities).  In the foundational paper \cite{vl},
Graber and Pandharipande showed that the localization formula
\eqref{localization} works ``virtually'' on the moduli space of
stable maps, where fundamental classes are replaced by virtual
fundamental classes, and normal bundles are replaced by ``virtual
normal bundles''.  They defined the virtual fundamental class of a
fixed locus, and the virtual normal bundle, and developed the
machinery to deal with such questions.

There is one pedantic point that must be made here.  The localization
formula should reasonably be expected to work in great generality.
However, we currently know it only subject to certain technical
hypotheses.  (1) The proof only works in the algebraic category.  (2)
In order to apply this machinery, the moduli space must admit 
a $\C^*$-equivariant locally closed immersion
into an orbifold.  (3) The virtual fundamental class of this fixed
locus needs to be shown to arise from the $\C^*$-fixed part of the
obstruction theory of the moduli space.  It would be very interesting,
and potentially important, to remove  hypotheses (1) and (2).

The theory of virtual localization can be applied to our relative
setting \cite{thmstar}.  (See \cite{llz} for more discussion.)  We now
describe it in the case of interest to us, of maps to $\proj^1$.
Again, in order to understand this properly, you should work out
examples yourself.

Fix a torus action on $\proj^1$
$$\si \circ [x;y] = [\si x ; y],$$ so the torus
acts with weight $1$ on the tangent space at $0$ and $-1$ on
the tangent space at $\infty$.
(The {\em weight} is the one-dimensional representation, or equivalently,
the character.)  This torus action
induces an obvious torus action on $\cmbar_{g, \al}(\proj^1, d)^{\dis}$
(and $\cmbar_{g, \al}(\proj^1, d)$).

We first determine the torus-fixed points of this action.  Suppose $C
\rightarrow T \rightarrow X$ is such a fixed map. A picture of two
fixed maps showing ``typical'' behavior is given in
Figure~\ref{fixed}.  The first has ``nothing happening above
$\infty_X$'', and the second has some ``sprouting'' of $T_i$'s.
\begin{figure}
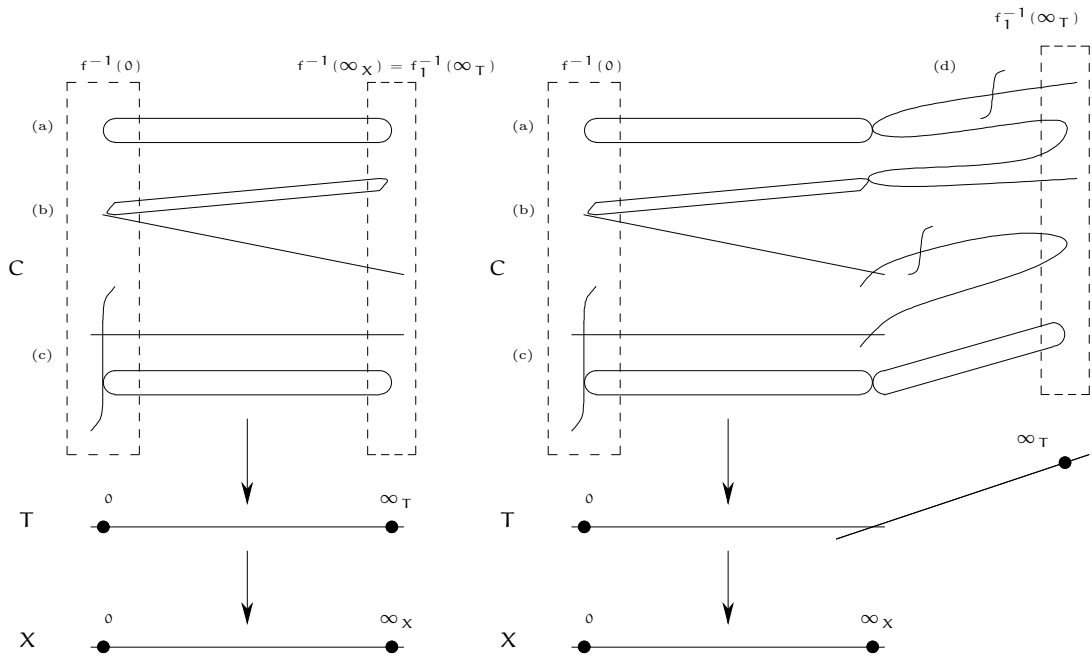

\begin{center}
\include{fixed}
\end{center}
\caption{Two examples of torus-fixed stable relative maps to $(\proj^1, \infty)$
\lremind{fixed} }\label{fixed}
\end{figure}

The map $C \rightarrow X$ must necessarily be a covering space away
from the points $0$ and $\infty$ of $X=\proj^1$.

\noindent {\em Exercise.}
Using the Riemann-Hurwitz formula, show that a surjective map $C'
\rightarrow \proj^1$ from an irreducible curve, unbranched way from
$0$ and $\infty$ must be of the form $\proj^1 \rightarrow \proj^1$,
$[x;y] \mapsto [x^a; y^a]$ for some $a$.

\point Hence the components dominating $X$ must be a union of ``trivial
covers'' of this sort.  \label{trivialcoverdef}\lremind{trivialcoverdef}

We now focus our attention on the preimage of $0$.  Any sort of
(stable) behavior above $0$ is allowed.  For example, the curve could
be smooth and branched there (Figure~\ref{fixed}(a)); or two of the
trivial covers could meet in a node (Figure~\ref{fixed}(b)); or there
could be a contracted component of $C$, intersecting various trivial
components at nodes (Figure~\ref{fixed}(c)).  (Because the
``relative'' part of the picture is at $\infty$, this discussion is
the same as the discussion for ordinary stable maps, as discussed in
\cite{vl}.)

Finally, we consider the preimage of $\infty_X$.  Possibly ``nothing
happens over $\infty$'', i.e.\ the target has not sprouted a tree
($l=0$ in our definition of stable relative maps at the start of \S
\ref{s:rvl}), and the preimage of $\infty$ consists just of $n$ smooth
points; this is the first example in
Figure~\ref{fixed}.  Otherwise, there is some ``sprouting'' of the target, and something
``nontrivial'' happens above each sprouted component $T_i$ ($i>0$), as in
Figure~\ref{fixed}(d).

\point At this point, you should draw some pictures, and convince yourself of
the following important fact: the connected components of the fixed
locus correspond to certain discrete data.  In particular, each
connected component can be interpreted as a product of three sorts
of moduli spaces:\label{atthispoint}\lremind{atthispoint}
\begin{enumerate}
\item[(A)] moduli spaces of pointed curves (corresponding to Figure~\ref{fixed}(c))
\item[(B)] (for those fixed loci where ``something happens above $\infty_X$,
i.e.\ Figure~\ref{fixed}(d)), a moduli space of maps parametrizing the
behaviour there.  This  moduli space is a variant of the space of
stable relative maps, where there is no ``rigidifying'' map to $X$.
We denote such a moduli space by
$\cmbar_{g,\al,\be}(\proj^1,d)_{\sim}$; its theory (of deformations
and obstructions and virtual fundamental classes) is essentially the
same as that for $\cmbar_{g,\al, \be}(\proj^1,d)$.
The virtual dimension of 
$\cmbar_{g,\al,\be}(\proj^1,d)_{\sim}$ is one less than that of
$\cmbar_{g,\al,\be}(\proj^1,d)$.

\item[(C)] If $\al_1 + \cdots + \al_n = d$ is the partition corresponding
to the ``trivial covers'' of $T_0$, these stable relative maps
have automorphisms $\Z_{\al_1} \times \cdots  \times \Z_{\al_n}$ corresponding
to automorphisms of these trivial covers 
(i.e.\ if one  trivial cover is of the form $[x;y] \mapsto [x^{\al_1}; y^{\al_1}]$, 
and $\zeta_{\al_1}$ is a $\al_1$th root of unity, then
$[x;y] \mapsto [\zeta_{\al_1} x ; y]$ induces an automorphism  of 
the map).  In the language of stacks, we can include a factor of 
$B\Z_{\al_1} \times \cdots \times B\Z_{\al_n}$; but the reader may prefer
to simply divide the virtual fundamental class by  $\prod \al_i$ instead.
\end{enumerate}

Each of these spaces has a natural virtual fundamental class:
the first sort has its usual fundamental class, and the 
second has its intrinsic virtual fundamental class.

The {\em relative virtual localization} formula states that
$$
[\cmbar_{g, \al}(\proj^1,d)]^{\virt} = \sum_{\fixed} [Y_i]^{\virt} 
/ e(N_i^{\virt}),
$$
in the equivariant homology of $\cmbar_{g,\al}(\proj^1,d)$ (cf.\
\eqref{localization}), with suitable terms inverted, where the virtual
fundamental classes of the fixed loci are as just described, and the
``virtual normal bundle'' will be defined now.

Fix attention now to a fixed component $Y_i$.  The virtual normal
bundle is a class in equivariant $K$-theory.  The term $1/
e(N_i^{\virt})$ can be interpreted as the product of several factors,
each ``associated'' to a part of the picture in Figure~\ref{fixed}.
We now describe these contributions.  The reader is advised to not
worry too much about the precise formulas; the most important thing is
to get a sense of the shape of the formula upon a first exposure to
these ideas.  Let $t$ be the generator of the equivariant cohomology
of a point (i.e.\ $H_T^*(\rm{pt}) = \Z[t]$).

{\em 1.}  For each irreducible component dominating $T_0$ (i.e.\ each 
trivial cover) of degree $\al_i$, we have a contribution of 
$\frac {\al_i^{\al_i}}  {\al_i! t^{\al_i}}$.

{\em 2.}  For each contracted curve above $0$ (Figure~\ref{fixed}(c))
of genus $g'$, we have a contribution of $(t^{g'} - \la_1 t^{g'-1}
+ \cdots + (-1)^{g' } \la_{g'}) / t$.  (This contribution is on
the factor $\cmbar_{g',n}$ corresponding to the contracted curve.)

{\em 3.}  For each point where
 a trivial component of degree $\al_i$
meets  a contracted curve above $0$ at a point $j$,
we have a contribution of $t /  (t/\al_i - \psi_j)$.
here, $\psi_j$ is a class on the moduli space $\cmbar_{g',n}$ corresponding to the
contracted component.

{\em 4.}  For each node above $0$ (Figure~\ref{fixed}(b)) joining
trivial covers of degrees $\al_i$ and $\al_j$, we have a contribution
of $1 / (t/\al_1 + t/ \al_2)$.

{\em 5.}  For each smooth point above $0$ (Figure~\ref{fixed}(a))
on a trivial cover of degree $\al_i$, we have a contribution of
$t/\al_i$.

At this point, if you squint and ignore the $t$'s, 
you can almost see the ELSV
formula \eqref{elsv} taking shape.

{\em 6.}  If there is a component over $\infty_X$, then we have a contribution
of $1/ (-t-\psi)$, where $\psi$ is the first Chern class of the 
line bundle corresponding to the cotangent space of $T_1$ at the point
where it meets $T_0$.

These six contributions look (and are!) complicated.  But this formula
can be judiciously used to give some powerful results, surprisingly
cheaply.  We now describe some of these.

\section{Applications of relative virtual localization}

\bpoint{Example 1:  proof of the ELSV formula \eqref{elsv}}
\label{elsvproof}\lremind{elsvproof}

As a first example, we prove the ELSV formula \eqref{elsv}.  
(This formula follows \cite{hhv}, using the simplification in the 
last section of \cite{hhv} provided by the existence of Jun Li's description
of the moduli space of stable relative maps.)
The ELSV
formula counts branched covers with specified branching over $\infty$
corresponding to $\al \vdash d$, and other fixed simple branched
points.  Hence we will consider $\cmbar_{g,\al}(\proj^1,d)$.

We next need to impose other fixed branch points.  There is a natural
Gromov-Witten-theoretic approach involving using descendant
invariants, but this turns out to be the wrong thing to do.  Instead,
we use a beautiful construction of Fantechi and Pandharipande
\cite{fp}.  Given any map from a nodal curve to $\proj^1$, we can
define a {\em branch divisor} on the target.  When the source
curve is smooth, the definition is natural (and old):
above a point $p$ corresponding to a partition $\be \vdash d$, 
the branch divisor contains $p$ with multiplicity $\sum (\be_i-1)$.
It is not hard to figure out how extend this to the case where
the source curve is not smooth above $p$.

\noindent {\em Exercise.}  Figure out what this extension should be.
(Do this so that the Riemann-Hurwitz formula remains true.)

Thus we have a map of sets $\cmbar_g(\proj^1,d) \rightarrow
\Sym^{2d+2g-2} \proj^1$.  In the case of stable relative maps, we have
a map of sets $\cmbar_{g,\al}(\proj^1,d) \rightarrow \Sym^{2d+2g-2}
\proj^1$.  As each  such stable relative map will have 
branching of at least $\sum(\al_i-1)$ above $\infty$, we can subtract
this fixed branch divisor to get a map of sets\lremind{br}
\begin{equation}
\label{br}
br: \cmbar_{g,\al}(\proj^1,d) \rightarrow \Sym^r \proj^1
\end{equation}
where $r = 2d+2g-2 - \sum(\al_i-1) = 2g-2+d+n$
(cf.\ \eqref{rh2}).

The important technical result proved by Fantechi and Pandharipande
is the following.

\tpoint{Theorem (Fantechi-Pandharipande \cite{fp})}
{\em There is a natural map of stacks $br$ as in \eqref{br}.}

We call such a map a {\em (Fantechi-Pandharipande) branch morphism}.
This morphism respects the torus action.

One can now readily verify several facts.  If the branch divisor does
not contain $p \neq \infty$ in $\proj^1$, then the corresponding map
$C \rightarrow \proj^1$ is unbranched (i.e.\ a covering space, or
\'etale) above $p$.  If the branch divisor contains $p \neq \infty$
with multiplicity $1$, then the corresponding map is simply branched
above $p$. (Recall that this means that the preimage of $p$
consists of smooth points, and the branching corresponds to the partition
$2+1+\cdots+1$.)
If the branch divisor does not contain $\infty$, i.e.\ there is no
additional branching above $\infty_X$ beyond that required by the definition of
stable relative map, then the preimage of $\infty_X$ consists
precisely of the $n$ smooth points $q_i$.  In other words, there is no
``sprouting'' of $T_i$, i.e.\ $T \cong \proj^1$. Hence
if $p_1 + \cdots + p_r$ is a general point of $\Sym^r \proj^1$, then
$br^{-1}(p_1 + \cdots + p_r) \subset  \cmbar_{g, \al}(\proj^1, d)$ 
is a finite set of cardinality equal to  the Hurwitz number $H^g_{\al}$.
This is true despite the fact that $\cmbar_{g,\al}(\proj^1, d)$
is horribly non-equidimensional --- the preimage of a general
point of $\Sym^r \proj^1$ will be contained in $\cm_{g,\al}(\proj^1,d)$, and
will not meet any other nasty components!

By turning this set-theoretic argument into something
more stack-theoretic and precise, we have that \lremind{thingy}
\begin{equation}
\label{thingy}
H^g_{\al} = \deg br^{-1}( pt ) \cap [\cmbar_{g,\al}(\proj^1, d)]^{\virt}.
\end{equation}  (For distracting
unimportant reasons, the previous paragraph's discussion is slightly
incorrect in the case where $H^g_{\al} = 1/2$, but \eqref{thingy}
is true.)

We can now calculate the right side of \eqref{thingy} using localization.
In order to do this, we need to interpret it equivariantly,
which involves choosing an equivariant lift of $br^{-1}$ of a point
in $\Sym^r \proj^1 \cong \proj^r$.  We do this by choosing our point
in $\Sym^r \proj^1$ to be the point $0$ with multiplicity $r$.
Thus all the branching (aside from that forced to be at $\infty$) 
must be at $0$.  The normal bundle to this point of $\proj^r$ 
is $r! t^r$.  Thus when we apply localization, a miracle happens.
The only fixed loci we consider are those where there is no extra branching
over $\infty$ (see the first picture in Figure~\ref{fixed}).  However, the source curve is smooth, so there is in fact only one connected
component of the fixed locus to consider, which is shown in Figure~\ref{elsvfixed}.
The moduli space in this case is $\cmbar_{g,n}$, which we
take with multiplicity $1/ \prod \al_i$ (cf.\ \S \ref{atthispoint}(C)).
Hence the Hurwitz number is the intersection on this moduli space
of the contributions to the virtual normal bundle outlined above.
\begin{figure}
\begin{center}
\setlength{\unitlength}{0.00083333in}
\begingroup\makeatletter\ifx\SetFigFont\undefined%
\gdef\SetFigFont#1#2#3#4#5{%
  \reset@font\fontsize{#1}{#2pt}%
  \fontfamily{#3}\fontseries{#4}\fontshape{#5}%
  \selectfont}%
\fi\endgroup%
{\renewcommand{\dashlinestretch}{30}
\begin{picture}(3762,2322)(0,-10)
\put(162.000,1845.000){\arc{150.000}{1.5708}{4.7124}}
\put(162.000,1095.000){\arc{150.000}{1.5708}{4.7124}}
\put(3162.000,1845.000){\arc{150.000}{4.7124}{7.8540}}
\put(3162.000,1095.000){\arc{150.000}{4.7124}{7.8540}}
\put(87,45){\blacken\ellipse{74}{74}}
\put(87,45){\ellipse{74}{74}}
\put(3237,45){\blacken\ellipse{74}{74}}
\put(3237,45){\ellipse{74}{74}}
\path(162,1920)(3162,1920)
\path(162,1770)(3162,1770)
\path(87,1470)(3237,1470)
\path(162,1170)(3162,1170)
\path(162,1020)(3162,1020)
\path(87,2070)(87,870)
\path(1587,945)(1587,195)
\blacken\path(1557.000,315.000)(1587.000,195.000)(1617.000,315.000)(1587.000,279.000)(1557.000,315.000)
\path(12,45)(3312,45)
\path(87,2070)(87,2071)(87,2078)
	(88,2092)(88,2112)(89,2133)
	(90,2153)(92,2170)(94,2185)
	(96,2197)(100,2207)(103,2218)
	(108,2227)(114,2237)(122,2248)
	(131,2260)(142,2272)(151,2283)
	(158,2291)(161,2294)(162,2295)
\path(87,870)(87,869)(87,862)
	(86,848)(86,828)(85,807)
	(84,787)(82,770)(80,755)
	(78,743)(75,732)(71,722)
	(66,713)(60,703)(52,692)
	(43,680)(32,668)(23,657)
	(16,649)(13,646)(12,645)
\put(3162,270){\makebox(0,0)[lb]{\smash{{{\SetFigFont{5}{6.0}{\rmdefault}{\mddefault}{\updefault}$\infty$}}}}}
\put(62,270){\makebox(0,0)[lb]{\smash{{{\SetFigFont{5}{6.0}{\rmdefault}{\mddefault}{\updefault}$0$}}}}}
\put(3762,45){\makebox(0,0)[lb]{\smash{{{\SetFigFont{8}{9.6}{\rmdefault}{\mddefault}{\updefault}$T=X$}}}}}
\end{picture}
}
\end{center}
\caption{The only fixed locus contributing to our calculation of the Hurwitz number
\lremind{elsvfixed}}\label{elsvfixed}
\end{figure}
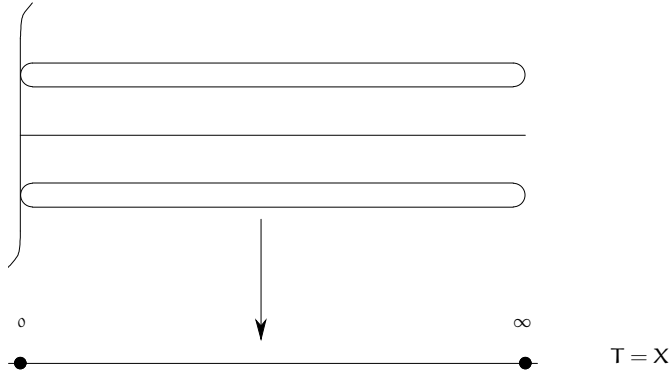

\noindent {\em Exercise.}  Verify that the contributions from {\em 1},
{\em 2}, and {\em 3} above, on the moduli space $\cmbar_{g,n}$, give the
ELSV formula \eqref{elsv}.

\bpoint{Example 2:  Proof of Theorem $\star$ (Theorem~\ref{thmstarhere})}\label{thmstarpf}\lremind{thmstarpf}

In Example 1 (\S \ref{elsvproof}), we found an equality of numbers.  Here we
will use relative virtual localization to get equality of cohomology
or Chow classes.

Fix $g$ and $n$.  We are interested in dimension $j$ (tautological)
classes on $\cmbar_{g,n}$.  In particular, we wish to show that any
such {\em tautological} class can be deformed into one supported on the locus
corresponding to curves with at least $2g-2+n-j$ genus $0$ components.
(This is just a restatement of Theorem $\star$.)
Call such a dimension $j$ class {\em good}.
Using the definition of the tautological ring in terms of $\psi$-classes,
it suffices to show that monomials in $\psi$-classes of dimension $j$
(hence degree=codimension $3g-3+n-j$) are good.

Here is one natural way of getting dimension $j$ classes.  Take any
partition $\al_1 + \cdots + \al_n = d$.  Let $r= 2g-2+n+d$ be the
virtual dimension of $\cmbar_{g,\al}(\proj^1,d)$ (i.e.\ the dimension
of the virtual fundamental class, and the dimension of
$\cm_{g,\al}(\proj^1,d)$), and suppose $r>j$.  Define the {\em Hurwitz
class} $\H^{g, \al}_j$ by 
$$
\H^{g,\al}_j := \pi_* \left( (\cap_{i=1}^{r-j} br^{-1}(p_i)) \cap
  [\cmbar_{g,\al}(\proj^1,d)]^{\virt} \right) \in A_j(\cmbar_{g,n})$$
where $\pi$ is the moduli map $\cmbar_{g,\al}(\proj^1,d) \rightarrow
\cmbar_{g,n}$ (and the $n$ points are the preimages of $\infty$), and
$p_1$, \dots, $p_{r-j}$ are generally chosen points on $\proj^1$.  We
think of this Hurwitz class informally as follows: consider branched
covers with specified branching over $\infty$.  Such covers (and their
generalization, stable relative maps) form a space of (virtual)
dimension $r$.  Fix all but $j$ branch points, hence giving a class of
dimension $j$.  Push this class to the moduli space $\cmbar_{g,n}$.

We get at this in two ways, by deformation and by localization.

\noindent {\bf 1.  Deformation.} 
(We will implicitly use Li's degeneration formula here.)  Deform the
target $\proj^1$ into a chain of $r-j$ $\proj^1$'s, each with one of
the fixed branch points $p_i$.  Then you can (and should) check that
the stabilized source curve has lots of rational components,
essentially as many as stated in Theorem $\star$.  (For example,
imagine that $r \gg 0$.  Then the $j$ ``roving'' branch points can lie
on only a small number of the $r-j$ components of the degenerated
target.  Suppose $\proj^1$ is any other component of the target, where $0$
and $\infty$ correspond to where it meets the previous and next
component in the chain.  Then the cover restricted to this $\proj^1$
can have arbitrary branching over $0$ and $\infty$, and only one other
branch point: simple branching above the $p_i$ lying on it.  This
forces the cover to be a number of trivial covers, plus one other
cover $C \rightarrow \proj^1$, where $C$ is simply branched at $p_i$,
and has one point above $0$ and two points above $\infty$, or vice
versa, forcing $C$ to be genus $0$, with three node-branches.)  Thus
any dimension $j$ Hurwitz class is {\em good}, i.e.\ satisfies the
conclusion of Theorem $\star$.

\noindent {\bf 2. Localization.}  We next use localization to express
tautological classes in terms of Hurwitz classes.  In the same way as
for the ELSV formula, we choose an equivariant lifting of
$\cap_{i=1}^{r-j} br^{-1}(p_i)$, corresponding to requiring all the
$p_i$ to go to $0$.  (Unlike the ELSV case, there are still $j$ branch
points that could go to either $0$ or $\infty$.)

We now consider what fixed components can arise.  

We have one ``main'' component that is similar to the ELSV case, where
all the $j$ ``roving'' branch points go to $0$
(Figure~\ref{elsvfixed}).  Any other component will be nontrivial
over $\infty$.  One can readily inductively show that these other
components are {\em good}, i.e.\ satisfy the conclusion of Theorem
$\star$.  (The argument is by looking at the contribution from such a
fixed locus.  The part contained in $f^{-1}(\infty_X)$ is essentially
a Hurwitz class, which we have shown is {\em good}.  The part contained in
$f^{-1}(0)$ corresponds to tautological classes on moduli spaces of
curves with smaller $2g-2+n$, which can be inductively assumed to 
be {\em good}.)

Thus we have shown that the contribution of the ``main'' component is
{\em good}.  But this contribution is straightforward to contribute:
it is (up to multiple)
the dimension $j$ component of 
$$
\frac { 1 - \la_1 + \cdots + (-1)^g \la_g}
{ ( 1 - \al_1 \psi_1) \cdots (1 - \al_n \psi_n)}
$$
(compare this to the ELSV formula \eqref{elsv}).  By expanding this
out, we find a polynomial in the $\al_i$ of degree $3g-3+n-j$ (cf.\
\eqref{Pgn} for a similar argument earlier).  We then apply the same
trick as when we computed top intersections of $\psi$-classes using
Hurwitz numbers in \S \ref{trick}: we can recover the
co-efficients in this polynomial by ``plugging in enough values''.
In other words, 
$\psi_1^{a_1} \cdots \psi_n^{a_1}$ may be obtained (modulo
{\em good} classes) as a linear combination  of Hurwitz classes.
As Hurwitz classes are themselves {\em good}, we have shown that the monomial
$\psi_1^{a_1} \cdots \psi_n^{a_1}$ is also {\em good}, completing the argument.

\section{Towards Faber's intersection number conjecture \ref{taketwo} via
relative virtual localization}

We can use the methods of the proof of Theorem $\star$ to
combinatorially describe the top intersections in the tautological
ring.  Using this, one can prove the ``vanishing'' or ``socle''
portion of the Faber-type conjecture for curves with rational tails
(and hence for $\cm_g$), and prove Faber's intersection number
conjecture for up to three points.  Details will be given in
\cite{gjv}; here we will just discuss the geometry involved.

The idea is as follows.  We are interested in the Chow ring of
$\cm_{g,n}^{rt}$, so we will work on compact moduli spaces, but
discard any classes that vanish on the locus of curves with rational
tails.  We make a series of short geometric remarks.

First, note that $R_{2g-1}(\cm^{rt}_{g,n}) \rightarrow
R_{2g-1}(\cm_{g,1})$ is an isomorphism, and $R_{2g-1}(\cm_{g,1})
\rightarrow R_{2g-1}(\cm_g)$ is a surjection.  The latter is immediate
from our definition.  The argument for the former is for example
\cite[Prop.~5.8]{thmstar}, and can be taken as an exercise for the reader
using Theorem~$\star$.  Faber showed \cite[Thm.~2]{faber} that
$R_{2g-1}(\cm_g)$ is non-trivial, so if we can show that
$R_{2g-1}(\cm_{g,1})$ is generated by a single element, then we will
have proved that $R_{2g-1}(\cm^{rt}_{g,n}) \cong \Q$ for all $n \geq 0$.

\point An extension of that argument using Theorem $\star$ shows that if we
have a Hurwitz class of dimension less than $2g-1$ (i.e.\ with fewer
than $2g-1$ ``moving branch points''), then the class is $0$ in $A_*
(\cm^{rt}_{g,n})$. \label{extension}\lremind{extension}

In order to get a hold of $R_{2g-1}(\cm^{rt}_{g,n})$, we will again
use branched covers.  Before getting into the Gromov-Witten theory, we make a
series of remarks, that may be verified by the reader, using only 
the Riemann-Hurwitz formula \eqref{rh}.

\point First, suppose we have a map $C \rightarrow \proj^1$ from a nodal
(possibly disconnected) curve, unbranched away from $0$ and $\infty$.
Then it is a union of a trivial covers (in the sense of \S \ref{trivialcoverdef}).

\point Second, suppose we have a map from a nodal curve $C$ to
$\proj^1$, with no branching away from $0$ and $\infty$ except for
simple branching over $1$, and nonsingular over $0$ and $\infty$.  Then
it is a union of trivial covers, plus one more component, that is
genus $0$, completely branched over one of $\{ 0, \infty \}$, and with
two preimages over the other.  More generally, suppose we have a map from
some curve $C$ to a chain of $\proj^1$'s, satisfying the kissing
condition, unbranched except for two smooth points $0$ and $\infty$ on
the ends of the chain, and simple branching over another point $1$.
Then the map is the union of a number of trivial covers glued
together, plus one other cover $\proj^1 \rightarrow \proj^1$ of the
component containing $1$, of the sort described in the previous
sentence.  \label{one}\lremind{one}

\point Third, if we have a map from a nodal curve $C$ to $\proj^1$, with total
branching away from $0$ and $\infty$ of degree less than $2g$, and nonsingular
over $0$ and $\infty$, then $C$ has no component of geometric genus $g$.
In the same situation, if the total branching away from $0$ and $\infty$
is exactly $g$, and $C$ has a component of geometric genus $g$, then
the cover is a disjoint union of trivial covers, and one connected curve $C'$
 of arithmetic genus $g$, where the map $C' \rightarrow \proj^1$ is contracted
to $1$ or completely 
branched over $0$ and $\infty$. \label{two}\lremind{two}

More generally, if we have a map from a curve $C$ to a chain of
$\proj^1$'s satisfying the kissing condition, with $0$ and $\infty$
points on either ends of the chain, with total branching away less
than $2g$ away from $0$, $\infty$, and the nodes, then $C$ has no
component of geometric genus $g$.  In the same situation, if the total
branching away from $0$, $\infty$, and the nodes is precisely $2g$,
then the map is the union of a number of trivial covers glued
together, plus one other cover of the sort described in the previous
paragraph.

\point The following fact is trickier.  Let $Z_{g,d}$ be the image in
$A_{2g-1}(\cm_{g,1})$ of $br^{-1}(1) \cap [
\cmbar_{g,(d),(d)}(\proj^1,d) ]^{\virt}$ (where the point in
$\cm_{g,1}$ is the preimage of $\infty$).  Then $Z_{g,d} = d^{2g}
Z_{g,1}$.  (We omit the proof, but the main idea behind this is the
Fourier-Mukai fact \cite[Lemma~2.10]{looijenga}.) 
\label{twog}\lremind{twog}

Define the Faber-Hurwitz class $\F^{g,\al}$ as the image in $A_{2g-1}(\cmbar_{g,n}^{rt})$
of 
$$
\cap_{i=1}^{r-(2g-1)} br^{-1} (p_i) \cap [
\cmbar_{g,1}(\proj^1,d)]^{\virt}$$ where the $p_i$ are general points
of $\proj^1$.  (This is the image of a Hurwitz class in
$\cmbar_{g,n}^{rt}$.)

As with the proof of Theorem $\star$, we get at this class inductively
using degeneration, and connect it to intersections of $\psi$-classes
using localization.

\bpoint{Degeneration}

Break the target into two pieces $\xymatrix{\proj^1 \ar@{~>}[r] &
  \proj^1 \cup \proj^1}$, where $\infty$ and one $p_i$ are on the
``right'' piece, and the remaining $p_i$'s are on the ``left'' piece.
The Faber-Hurwitz class breaks into various pieces; we enumerate the
possibilities.  We are interested only in components where there is a
nonsingular genus $g$ curve on one side.  We have two cases, depending
on whether this curve maps to the ``left'' or the ``right'' $\proj^1$.

\point If it maps to the left component, then all $2g-1$ ``moving'' branch
points must also map to the left component in order to get a non-zero
contribution in $A_*(\cm_{g,n}^{rt})$, by Remark~ \ref{two}.  Thus by
Remark~ \ref{one}, the cover on the right is of a particular sort,
and the cover on the left is another Faber-Hurwitz class, where
one of the branch points over $\infty$ has been replaced two,
or where two of the branch points are replaced by one.\label{left}
\lremind{left}

\point If the genus $g$ curve maps to the right component, then all $2g-1$
``moving'' branch points must map to the right component, and by
Remark~\ref{two} our contribution is a certain multiple of $Z_{g,d}$,
which by Remark~\ref{twog} is a certain multiple of $Z_{g,1}$.  The
contribution from the left is the genus $0$ Hurwitz number
$\H^0_{\al}$, for which Hurwitz gives us an attractive formula
\eqref{hformula}. \label{right} \lremind{right}

Unwinding this gives the recursion\lremind{joincutprime}
\begin{equation}\label{joincutprime}
\F^g_{\al} = \sum_{i+j=\al_k} ij H^0_{\al'} \F^g_{\al''} \binom {
  d+l(\al)-2 } { d+l(\al') -2, d+l(\al'')-1} + \sum_{\al_i + \al_j}
\F^g_{\al'} + \sum_{i=1}^{l(\al)} \al_i^{2g+1} H^0_{\al} Z_{g,1}.
\end{equation} In
this formula, the contributions from paragraph \S \ref{right} are in the first
two terms on the right side of the equation, and the contributions from \S \ref{left} are
in the last.  The second term on the right corresponds to where two
parts $\al_i$ and $\al_j$ of $\al$ are ``joined'' by the nontrivial
cover of the right $\proj^1$ to yield a new partition where $\al_i$
and $\al_j$ are replaced by $\al_i + \al_j$.  The first term on the
right corresponds to where one part $\al_k$ of $\al$ is ``cut'' into
two pieces $i$ and $j$, forcing the curve covering the left $\proj^1$
to break into two pieces, one of genus $0$ (corresponding to partition
$\al'$ and one of genus $g$ (corresponding to $\al'$).  The binomial
coefficient corresponds to the fact that the $d+l(\al) -2$ fixed
branch points $p_1$, $p_2$, \dots on the left component must be split
between these two covers.

The base case is $\F^g_{(1)} = Z_{g,1}$.  Hence we have shown that
$\F^g_{\al}$ is always a multiple of $Z_{g,1}$, and the theory of
cut-and-join type equations (developed notably by Goulden and Jackson)
can be applied to solve for $\F^g_{\al}$ (in generating function form)
quite explicitly.

\bpoint{Localization}

We now get at the Faber-Hurwitz class by localizing. As with the proof
of Theorem~$\star$, we choose a linearization on $br^{-1}(p_i)$ that
corresponds to requiring all the $p_i$ to move to $0$. We now describe
the fixed loci that contribute.  We won't worry about the precise
contribution of each fixed locus; the important thing is to see the
shape of the formula.

First note that as we have only $2g-1$ moving branch points, in any
fixed locus in the ``rational-tails'' locus, our genus $g$ component
cannot map to $\infty$, and thus must be contracted to $0$.  The fixed
locus can certainly have genus $0$ components mapping to sprouted
$T_i$ over $\infty$, as well as genus $0$ components contracted to
$0$.  

We now look at the contribution of this fixed locus, via the relative
virtual localization formula.
We will get a sum of classes glued together from various moduli spaces
appearing in the description of the fixed locus (cf.\ \S \ref{atthispoint}).
Say the contracted genus $g$ curve
meets $m$ trivial covers, of degree $\be_1$, \dots, $\be_m$
respectively.  Then the contribution from this component will
be some summand of
$$
\frac {1 - \la_1 + \cdots + (-1)^g \la_g } { (1 -\be_1 \psi_1) \cdots
(1 - \be_m \psi_m)}$$
where $\psi_i$ are the $\psi$-classes on $\cmbar_{g,m}$.
Thus the contribution from this component is visibly tautological, and 
by Remark~\ref{extension} the contribution will be zero if the dimension
of the class is less than $2g-1$.  As the total contribution
of this fixed locus is $2g-1$, any non-zero contribution must
correspond to a dimension $2g-1$ tautological class on $\cmbar_{g,m}$
glued to a dimension $0$ class on the other moduli spaces appearing
in this fixed locus.  This can be readily computed; the genus $0$
components contracted to $0$ yield binomial coefficients, any
components over $\infty_X$ turn out to yield products of {\em genus $0$
double Hurwitz numbers}, which count branched covers of $\proj^1$
by a genus $0$ curve, with specified branching $\al$ and $\be$ above two 
points, and the remaining branching fixed and simple.

Equipped with this localization formula, even without worrying
about the specific combinatorics, we may show the following.

\tpoint{Theorem}
{\em For  any $n$, and $\be \vdash d$, 
$\pi_* \psi_1^{\be_1} \cdots \psi_n^{\be_n}$ is a multiple of $Z_{d,1}$,
where $\pi$ is the forgetful map to $\cm_{g,1}$.}
\label{Zd1}\lremind{Zd1}

We have thus fully shown the ``vanishing'' (or socle) part of
Faber's conjecture for curves with rational tails. (This
may certainly be shown by other means.)  In particular,
we have completed a proof of Looijenga's Theorem~\ref{lthm}.

\bpf Call such a class an $n$-point class.  We will show that such a
class is a multiple of $\Z_{d,1}$ modulo $m$-point classes, where
$m<n$; the result then follows by induction.  As with the proof of
Theorem $\star$, we consider $\F^{\al}_g$ as $\al$ runs over all
partitions of length $n$.  Each such Faber-Hurwitz class is a multiple
of $\Z_{d,1}$ by our degeneration analysis.  By our localization
analysis, all of the fixed points for $\F^{\al}_g$ yield $m$-point
classes where $m<n$ except for one, corresponding to the picture in
Figure~\ref{elsvfixed}.  The contribution of this component is some
known multiple of a polynomial in $\al_1$, \dots, $\al_n$.  The
highest-degree coefficients of this polynomial are the $n$-point
classes, the monomials in $\psi$-classes that we seek.  By taking a
suitable linear combination of values of the polynomial (i.e.\
Faber-Hurwitz classes, modulo $m$-point classes where $m<n$), we can
obtain any co-efficient, and in particular, the leading co-efficients.
\epf

A related observation is that we have now given an explicit
combinatorial description of the monomials in $\psi$-classes, as a
multiple of our generator $Z_{g,1}$.  (In truth, we have not been
careful in this exposition in describing all the combinatorial
factors.  See \cite{gjv} for a precise description.)

This combinatorialization can be made precise as follows.
We create a generating function $\F$ for Faber-Hurwitz classes.
The join-cut equation \eqref{joincutprime} allows us to
solve for the generating function $\F$.

We make a second generating function $\W$ for the intersections $\pi_*
\psi_1^{\be_1} \cdots \psi_n^{\be_n} \la_k \in R_{2g-1}(\cm_{g,1})$
(where $\be_1 + \cdots + \be_n + k = g-2$).  Localization gives us a
description of $\F$ in terms of $\W$ (and also the genus $0$ double
Hurwitz generating function).  By inverting this relationship we can
hope to solve relatively explicitly for $\W$, and hence prove Faber's
intersection number conjecture.  Because genus $0$ double Hurwitz
number $H^0_{\al, \be}$  are only currently well-understood
when one of the partition has at most $3$ parts (see \cite{double}),
this program is not yet complete.  However, it indeed yields:

\tpoint{Theorem \cite{gjv}} {\em Faber's intersection number
  conjecture is true for up to three points.}

One might reasonably hope that this will give an elegant proof of
Faber's intersection number conjecture in full before long.

\section{Conclusion}

In the last fifteen years, there has been a surge of progress in
understanding curves and their moduli using the techniques of
Gromov-Witten theory.  Many of these techniques have been outlined
here.

Although this recent progress uses very modern machinery, it is part
of an ancient story.  Since the time of Riemann, algebraic curves have
been studied by way of branched covers of $\proj^1$.  The techniques
described here involve thinking about curves in the same way.
Gromov-Witten theory gives the added insight that we should work with
a ``compactification'' of the space of branched covers, the moduli
space of stable (relative) maps.  A priori we pay a steep price, by
working with a moduli space that is bad in all possible ways
(singular, reducible, not even equidimensional).  But it is in some
sense ``virtually smooth'', and its virtual fundamental class behaves
very well, in particular with respect to degeneration and localization.

The approaches outlined here have one thing in common:  in each case the key 
idea is direct and naive.  Then one works to develop the necessary
Gromov-Witten-theoretic tools to make the naive idea precise.

In conclusion, the story of using Gromov-Witten theory to understand
curves, and to understand curves by examining how they map into other
spaces (such as $\proj^1$), is most certainly not over, and may just
be beginning.

} % end of parskip; it started just before the introduction

\end{document}